\documentclass[twosided]{report}
\usepackage{amssymb}
\usepackage{a4} 
\usepackage{amsmath}

\usepackage{xy} 
\xyoption{all} 

\usepackage{latexsym}

\newcommand{\lie}[1] {\mathfrak{#1}}  
\newcommand{\bb}[1]{{\mathbb #1}}    

\newcommand{\bbR}{{\bb R}}
\newcommand{\bbA}{{\bb A}}
\newcommand{\bbZ}{{\bb Z}} 
\newcommand{\bbQ}{{\bb Q}}

\newcommand{\bbC}{{\bb C}}
\newcommand{\bbN}{{\bb N}}

\newcommand{\ftn}{{\rm f.t.n.}}

\newcommand{\wfn}{{\rm wfn}}
\newcommand{\vcd}{\mathop{{\rm vcd}}}

\newcommand{\rmod}{ / \,}
\newcommand{\lmod}{\backslash}
\newcommand{\GL}{{\rm GL}}

\renewcommand{\O}{{\rm O}} 
 
\newcommand{\U}{{\rm U}} 

\newcommand{\SP}{{\rm Sp}}

\newcommand{\Aff}{{\rm Aff}}
\newcommand{\Aut}{{\rm Aut}}
\newcommand{\E}{{\rm E}}

\newcommand{\End}{{\rm End}}

\newcommand{\Hom}{{\rm Hom}}
\newcommand{\Def}{{\cal D}_c}
\newcommand{\Mod}{{\cal M}_c}
\newcommand{\Cha}{{\cal C}_c}
\newcommand{\Out}{{\rm Out}}
\newcommand{\Inn}{{\rm Inn}}

\newcommand{\ad}{{\rm ad} \: }
\newcommand{\N}{{\rm N}}  
\newcommand{\Z}{{\rm Z}} 
\newcommand{\Fitt}{{\rm Fitt}} 
\newcommand{\ur}{{\rm u}} 
\newcommand{\C}{{\rm C}} 
\newcommand{\A}{{\rm A}} 

\newcommand{\bF}{{\mathbf{F}}}
\newcommand{\bC}{{\mathbf{C}}}
\newcommand{\bU}{{\mathbf{U}}}
\newcommand{\bG}{{\mathbf{G}}}
\newcommand{\bH}{{\mathbf{H}}}
\newcommand{\bL}{{\mathbf{L}}}

\newcommand{\bN}{{\mathbf{N}}}

\newcommand{\bD}{{\mathbf{D}}}
\newcommand{\bT}{{\mathbf{T}}}

\newcommand{\lu}{{\lie{u}}}
\renewcommand{\lg}{{\lie{g}}}
\renewcommand{\ln}{{\lie{n}}}
\newcommand{\la}{{\lie{a}}}
\newcommand{\lb}{{\lie{b}}}
\newcommand{\lgl}{{\lie{gl}}}

\newcommand{\id}{{\rm id}}

\renewenvironment{matrix}[1] {\left( \begin{array}{#1}}{\end{array}\right)}

\newenvironment{definition0}%
{ \par  \medskip  \noindent  {\bf Definition} \hspace{0.5em} }%
{\par \medskip } 

\newtheorem{example1}{Example}[section]

\newenvironment {example}%
{ \begin{example1} \rm }%
{ \end{example1}} 
\newenvironment {proof}%
{ \noindent {\em Proof. }}%
{\hspace*{\fill}$\Box$\par \medskip } 
\newenvironment{prf}[1]%
{ \noindent {\em #1 \hspace{0.5em}} }%
{\hspace*{\fill}$\Box$\par \medskip } 
\newenvironment {remark}%
{{\em Remark \hspace{0.5em}} }%
{\par \medskip }

\newtheorem{proposition}{Proposition}[chapter]
\newtheorem{definition1}[proposition]{Definition}
\newtheorem{theorem}[proposition]{Theorem}
\newtheorem{lemma}[proposition]{Lemma}
\newtheorem{corollary}[proposition]{Corollary}
\newenvironment{definition}%
{ \begin{definition1} \rm }%
{ \end{definition1}}

\newcommand{\rank}{{\rm  rank} \, } 
\newcommand{\ac}[1]{\overline{#1}} 

\author{Oliver Baues \thanks{e-mail: baues@math.uni-karlsruhe.de} \\
Institut f\"ur Algebra und Geometrie\\ 
Universit\"at Karlsruhe\\ 
D-76128 Karlsruhe
}

\title{Deformation Spaces for Affine Crystallographic Groups}

\begin{document}

\maketitle

\begin{abstract} 
We develop the foundations of the deformation
theory of compact complete affine space forms
and affine crystallographic groups. Using methods 
from the theory of linear algebraic groups
we show that these deformation spaces inherit
an algebraic structure from the space of
crystallographic homomorphisms.  We also study the 
properties of the action of the homotopy mapping
class groups on deformation spaces. In our context these
groups are arithmetic groups, and we construct examples 
of flat affine manifolds where every finite group of mapping
classes admits a fixed point on the deformation space. 
We also show that the existence of fixed points on the
deformation space is equivalent to the realisation of finite groups of
homotopy equivalences by finite groups of affine diffeomorphisms. 
Extending ideas of Auslander we relate the deformation spaces 
of affine space forms with solvable fundamental group to 
deformation spaces of manifolds with nilpotent fundamental group.
We give applications concerning the classification problem
for affine space forms.  \hspace{1cm}  \\

\hrule 

\hspace{1cm} \\

\noindent 
2000 Mathematics Subject Classification: Primary 57S30  22E40; \\
Secondary 22E40 22E45 20H15 58D29\\

\noindent 
Keywords: affine space form, affine crystallographic group, proper actions, 
discrete group, locally homogeneous manifold, deformation space, moduli space, mapping classes, homotopy equivalences, realisation problem,  linear algebraic group, arithmetic group, solvable group, polycyclic group, nilpotent group
\end{abstract}


\chapter*{Introduction}
This work is devoted to some aspects of the study of the 
fundamental groups of affine space forms. An affine space form
is a  {\em compact, complete, affinely flat manifold}. 
In modern terminology,  an affinely flat manifold is a
{\em locally homogeneous manifold\/} which is  modelled on the
action of the affine group $\Aff(V)$ on some  real vector space $V$.
One defines  more restricted kinds of affinely flat manifolds 
by considering any transitive subgroup $A$ of $\Aff(V)$.
For example,  if $A= E(V)$ is the group of Euclidean isometries
of $V$,  then one obtains the concept of a {\em Riemannian flat 
manifold}, and similarly one obtains geometrically interesting
cases like {\em Lorentz-flat\/} manifolds 
or {\em flat symplectic\/} manifolds.

It is well known that a compact complete affine manifold arises
as a quotient $\Gamma \lmod V$,  where $\Gamma \leq \Aff(V)$ 
is a properly discontinuously and cocompactly acting subgroup
of $\Aff(V)$. These groups are called {\em affine crystallographic groups}. 
An affine crystallographic subgroup  $\Gamma \leq \E(V)$ 
is called an {\em Euclidean\/} crystallographic group. 
The study of affine crystallographic groups has a long
history which goes back to Bieberbach's famous three
theorems and Hilbert's 18th problem on Euclidean crystallographic groups. 

Much of the recent
work on affine crystallographic groups is an effort to
generalize the satisfactory theory of Euclidean crystallographic
groups to the general case. This approach has encountered 
serious obstacles. A long standing conjecture (``Auslander's conjecture'')
states that any affine crystallographic group $\Gamma$ is a {\em virtually
polycyclic group}. Under the assumption that Auslander's conjecture
holds certain generalizations of Bieberbach's theory are possible. 
Still, for example, the question which virtually polycyclic groups 
are isomorphic to an affine crystallographic group 
remains mysterious.\footnote{Benoist \cite{Benoist} gave 
a striking example of a finitely generated nilpotent group which is
not isomorphic to an affine crystallographic group.}
\newline

In this article,  we study the deformation theory of complete affinely flat 
structures on a fixed compact manifold $M$.  Let $\Gamma$ be  a group.
An injective homomorphism $\rho: \Gamma \rightarrow \Aff(V)$, 
such that  $\rho$ gives rise to a proper action of $\Gamma$ on $V$ 
with compact quotient $\Gamma \lmod V$ is called a {\em crystallographic
homomorphism}. 
By the deformation theory of locally homogeneous spaces, developed
by Thurston,  the description of the topology of the
{\em deformation space\/} of complete affinely flat structures
on $M$  roughly amounts to an analysis of the space of
crystallographic homomorphisms of $\pi_{1}(M)$ up to
conjugacy by $\Aff(V)$.  This  approach, of course, has 
a long and fruitful history in geometry. One particular 
important example is the Teichm\"uller theory of surfaces.  

The phenomena one could expect for the topology (and
possibly geometry) of the
deformation spaces of complete  affine manifolds are, however, 
quite different from the situation of constant curvature 
geometric surfaces. 
The main theme here is 
that the deformation spaces of complete  affine manifolds
with virtually polycyclic fundamental group are of an
{\em algebraic\/} and {\em arithmetic} nature.  

We provide a conceptual framework
for the study of the set of  affine crystallographic actions of 
a given virtually polycylic group $\Gamma$. This builds 
on the strong interaction of affine crystallographic with
the theory of representations of Lie- and algebraic groups.
A basic result we obtain is that for any
virtually polycyclic group $\Gamma$ the space $\Hom_c(\Gamma, \Aff(V))$
of crystallographic embeddings into $\Aff(V)$ has a
natural structure of a real algebraic variety defined over
the rational numbers. We prove that $\Hom_c(\Gamma, \Aff(V))$
is Zariski-open in a certain closed subspace of the 
space of all homomorphisms of $\Gamma$ into
the affine group. This result may be seen as an analogue 
of a classic theorem of Weil on discrete
subgroups of Lie groups. 

The {\em deformation space\/} $\Def(\Gamma, \Aff(V))$ arises as a quotient
of $\Hom_c(\Gamma, \Aff(V))$ by the action of $\Aff(V)$. The 
algebraic structure of $\Hom_c(\Gamma, \Aff(V))$ carries over
to the deformation space which in some cases turns out to
be homeomorphic to a semi-algebraic
set, and in particular it is then a Hausdorff-space. This turns out to
be true,  for example,  for any manifold which is finitely 
covered by a torus and also for certain solvmanifolds. 
In general, the  structure of the quotient spaces $\Def(\Gamma, \Aff(V))$
seems particular hard to understand because the group actions 
which define the moduli problem are not reductive, and the usual
techniques of geometric invariant theory do not apply. \newline 

Of particular interest in this context is furthermore the action of the
mapping class group $\Out(\pi_{1}(M))$ on the deformation
space of complete flat affine structures, and its quotient, 
the moduli space of affine structures on $M$.  Compared with the
Teichm\"uller theory of conformal structures for surfaces,  
the general situation 
in our subject resembles more the case of a genus $1$ surface than the situation for higher genus. In fact, if $\Gamma$
is a polycyclic group, then $\Out(\Gamma)$ is an arithmetic
group, similar to the model situation $\Out(\bbZ^2) = \GL_{2}(\bbZ)$.\footnote{This result is however highly non-trivial, see  \cite{BauesGrunewald}}
Moreover, we show that the action of $\Out(\pi_{1}(M))$ extends to an
action of an algebraic group on the deformation space. 
In contrast to the situation for surfaces, however, the 
action of $\Out(\pi_{1}(M))$ on the deformation space
is not proper, and the moduli space may be highly singular. 

We also consider the following type of question: {\em Let $\Gamma$
be an affine crystallographic group and let $\Gamma \leq \Delta$
be a finite extension. Is $\Delta$ also isomorphic to an affine crystallographic
group?}  Questions like this are known as {\em realization problems}. 
From classic work of Dehn, Nielsen and Poincare in surface
theory it is known that there is a close interplay between such a
question and Teichm\"uller theory.   

We show in detail that, for virtually nilpotent groups, the solution of the
above realization problem depends only on the action of
the group of mapping classes for $\Gamma$ on the deformation
space. This leads us in turn to a description of the 
deformation space $\Def(\Delta, \Aff(V))$ of $\Delta$ 
as a fixed point set in  $\Def(\Gamma, \Aff(V))$ for 
the finite group of mapping classes which belongs to the
extension $\Delta$.  
We mention here that our 
contribution to the realization problem
may be understood as a suitable generalization of the
classic Burckhardt-Zassenhaus theorem on Euclidean 
crystallographic groups. This is because there exists 
a fixed point for every finite subgroup of $\GL_n(\bbZ)$
on the deformation space of Euclidean tori.
In this realm, we exhibit large clases of
infranilmanifolds where the finite group of mapping classes
have fixed points in $\Def(\Gamma, \Aff(V))$. 
Therefore the answer to the affine crystallographic
realization problem for these manifolds is positive. 

Another type of realization problem enters the stage, 
coming from Auslander's {\em nilshadow construction\/} 
for solvmanifolds. {\em Let $\Gamma$ be the fundamental
group of a solvmanifold, and let $\Theta$ be its nilshadow.
If $\Theta$ is an affine crystallographic group, is it
true that $\Gamma$ is isomorphic to an affine crystallographic
group?}
This question admits a treatment which is similar
to the solution of the problem
for finite extensions. The solution for this realization 
problem  paves the way for 
a description of the deformation spaces for affine
solvmanifolds, and 
enables us to construct 
new examples of such manifolds. 
\bigskip

\paragraph{Broader context of the article}
Let us mention briefly related work in the context of our article.
The deformation theory of locally homogeneous spaces was 
formulated and developed by Thurston, see \cite{Thurston, Goldman, Epstein}
for exposition. In similar spirit, but with different methods, is the local and
global deformation theory of complex structures on compact
manifolds, as pioneered by Kodaira and Spencer \cite{KodSpen}. 
For compact  surfaces both approaches lead to the same 
mathematical object, namely \emph{Teichm\"uller space}, which
is the analogue of our deformation space.  

Traditionally locally homogeneous spaces are modeled on
homogeneous spaces $G/H$ with $H$ compact. In this setup, 
discrete subgroups $\Gamma \leq G$ give rise to proper and,  if
$\Gamma$ is a uniform lattice, cocompact actions on $G/H$
in abundance. The case $H$ non-compact has been more or less ignored 
for quite some time, although its obvious relevance 
for relativity and Pseudo-Riemannian geometry. 
Motivated by phenomena discussed in the seminal paper \cite{Kulkarni}, 
there has been lots of recent activity to understand locally
homogeneous spaces modeled on $G/H$, where $G$ is reductive
and $H$ is non-compact. (See, \cite{Kob1, Benoist2, OhWitte, Salein} 
for important contributions.) 

Kobayashi \cite{Kob5} more generally formulated the program to 
classify the deformation and moduli spaces of proper 
actions of a discrete group $\Gamma$ 
on a homogeneous space $G/H$ with non-compact $H$.
The subject of this article with $G= \Aff(V)$ and $H= \GL(V)$,
and $\Gamma$ crystallographic is an important special case.\footnote{Margulis constructed proper non-uniform 
affine actions of finitely generated free groups on a 
three-dimensional vector space, see \cite{Margulis,Drumm}.}

Other recent activity in this program concerns the situation,
where $G$ is nilpotent and $\Gamma$ is acting properly but
not necessarily cocompact on $G/H$.  (See, for example 
\cite{Bak_Kh,Bak_Ke_Yo,KobNas,Yoshino}.)

In complex geometry,  there is related recent work 
of  Catanese et al.\ \cite{Catanese,Catanese2,Hase} 
to understand and classify the deformation of complex structures 
on compact nilmanifolds. The latter subject bears
certain similarities and interactions with the deformation 
of affine structures on these manifolds.

\paragraph{Structure of the article}
There are three main chapters. In chapter 1
we introduce and develop the basic techniques 
in the study of crystallographic actions of solvable
groups on affine space. The main result is
the description of the space of crystallographic 
homomorphism as an algebraic variety defined
over the rational numbers.
 
In chapter 2 we analyse the basic properties
of the associated deformation spaces. We
provide examples of such spaces which
have interesting properties with respect to
the action of the group of mapping classes. 
In particular, we study convexity and fixed 
point properties of deformation spaces.

Chapter 3 is technically most demanding and
concerns the realisation questions for affine
crystallographic groups, and the extension
of the action of the group of mapping classes
to an algebraic group action on the deformation
space. These constructions are carried out in the
context of nilmanifolds.

\paragraph{Acknowledgement} 
The author enjoyed many fruitful conversations
with  colleagues concerning the contents of
this article. He profited a lot from the insight
and interest in the subject shared by  Bill Goldman,  
Fritz Grunewald,  Dave Morris,  Yves Benoist.  
He wishes to thank further Herbert Abels, Richard Pink, 
Karel Dekimpe, Yoshinobu Kamishima,
Bernd Siebert, Mark Gross, Vicente Cortes,
Anders Karlsson for helpful comments and discussions. 
\newpage


\subsection*{Notational conventions}  

\paragraph{Groups} \hspace{1ex} \newline 
$\Z(G)$, the center of $G$ \\
$\C_{G}(H)$, the centralizer of $H$ in $G$ \\
$\Fitt(\Gamma)$, Fitting subgroup, the maximal nilpotent
normal subgroup of $\Gamma$\\
$H  \leq_f G$, $H$ is a subgroup of finite index in $G$\\
$[A,B] \leq G$, the commutator subgroup of $A,B \subseteq G$\\ 
$G$ is said to be {\em virtually} ${\cal P}$, if $G$ has a finite index subgroup which
satisfies ${\cal P}$\\
$G$ is called a {\em \wfn-group} if $G$ is without a non-trivial finite normal subgroup\\
\ftn-group, a finitely generated torsionfree nilpotent group

\paragraph{Algebraic sets}  \hspace{1ex} \newline  
We will frequently consider affine algebraic varieties 
defined over the rationals or reals. Our special interest
lies on the set of real points of these varieties.
We call these sets real algebraic varieties.\\ 
$\bar{M}$, the real Zariski closure of $M \subset X$, $X$ a real algebraic variety

\paragraph{Linear algebraic groups} \hspace{1ex} \newline  
A $\bbQ$-defined linear algebraic group $\bG$
is a subgroup $\bG \leq \GL_n(\bbC)$  defined
by polynomials with rational coefficients. For a
subring $R$ of $\bbC$ we put $\bG_{R} = \bG \cap \GL_n(R)$.\\
\noindent
$  g = g_s g_u$, the Jordan decomposition for $g \in \bG$. \\
$\ur(\Gamma)$, the maximal normal unipotent subgroup of $\Gamma \leq \bG$\\
$M_s = \{ m_s \mid m \in M \} $\\
$M_u= \{ m_u \mid m \in M \}$, where $M \subset \bG$ is a subset\\ 
$\ac{M}$, the Zariski-closure of $M$ in $\bG$
$\rho_u$, the map $\gamma \mapsto \rho(\gamma)_u$, where $\rho: \Gamma \rightarrow \bG$
is a homomorphism

\paragraph{Lie groups, real algebraic groups}  \hspace{1ex} \newline   
A Lie group $G$ is called real algebraic, if $G$ is (a connected component of)
the set of real points $\bG_{\bbR}$ of a linear algebraic group $\bG$\\  
$N(G)$,  nilpotent radical of G, the maximal nilpotent normal connected subgroup\\
$G^0$, identity component in the (real) Zariski-topology, if $G$ is real algebraic\\
$G_0$, identity component in the Hausdorff-topology\\
$\ac{G}$, the real Zariski-closure of a linear group $G \leq \GL_n \bbR$\\
$\ur(G)$, the unipotent radical, 
$G \subset \GL_n$ a linear group

\paragraph{Affine crystallographic groups} \hspace{1ex} \newline  
$V$, a finite dimensional real vector space\\  
$\Aff(V) = V \cdot \GL(V)$, the group of affine transformations of V\\
$A \leq \Aff(V)$, a Zariski closed subgroup of $\Aff(V)$, transitive on $V$\\
$\GL_A = A \cap \GL(V)$\\
$A_x = \{ a \in A \mid a x =x \}$, for $x \in V$\\
$\rho: \Gamma \rightarrow A$ a homomorphism,
$\rho^a(\gamma) = a \rho(\gamma) a^{-1}$

\begin{definition0} 
Let $\Gamma \leq \Aff(V)$ be a subgroup. $\Gamma$ is
called {\em properly discontinuous\/} if, for all 
compact subsets $K \subset V$, the set
$ \{   \, \gamma \in \Gamma \mid \gamma K \cap K \neq \emptyset \, \}$ 
is finite \\
$\Gamma \leq \A$ is called an {\em affine crystallographic group 
of type $A$}, if $\Gamma$ is properly discontinuous, and if the 
quotient space 
$   {_{\displaystyle \Gamma}} \,   \lmod \, V $ 
is compact.     
\end{definition0}
\noindent ACG, affine crystallographic group 



\tableofcontents

\chapter{Crystallographic Homomorphisms} \label{choms}

Let $V$ denote a finite dimensional real vector space, and  $A \leq \Aff(V)$
a Zariski closed subgroup, which acts transitively on $V$. 
A subgroup $\Gamma \leq \Aff(V)$
is called {\em affine crystallographic\/} if $\Gamma$ acts properly discontinuously
and with a compact quotient on $V$. If $\Gamma \leq A$ then $\Gamma$ is 
called an affine crystallographic group of type $A$. 
We define the {\em space of crystallographic homomorphisms\/} as
$$ \Hom_c(\Gamma, A) = 
\{ \rho: \Gamma \rightarrow A  \mid \mbox{$\rho$ is crystallographic} \} \;.$$
A homomorphism $\rho: \Gamma \rightarrow \Aff(V)$ 
is called a {\em crystallographic homomorphism\/} if $\rho$ is injective
and the image $\rho(\Gamma) \leq \Aff(V)$ is a crystallographic subgroup.
In 
this chapter we study the space of crystallographic
homomorphisms 
$$\Hom_c(\Gamma,A) \, \subset \, \Hom(\Gamma,A) $$ as a subset 
of the real algebraic variety $\Hom(\Gamma, A)$ 
of all homomorphisms of $\Gamma$ into $A$. 
We establish
that the space $\Hom_c(\Gamma,A)$ 
may be described by algebraic equalities and inequalities, and
carries itself a natural structure as a real algebraic variety.

\section{Polycyclic groups,  crystallographic groups} 
In this subsection we develop
the foundations of affine crystallographic virtually polycyclic groups.
We will need these results to build our later arguments on. Our basic references are
\cite{FriedGoldman} and \cite{GS} which are, however, written from different
points of view. We put the notion of algebraic hull for a torsionfree
polycyclic group (as developed by Mostow \cite{Mostow2}, 
see also, for example, \cite[Chapter IV]{Raghunathan}) at 
the center of our considerations. We hope that our presentation 
clarifies some aspects of the theory of polycyclic affine crystallographic
groups.\footnote{A generalisation in the context of 
affine actions on solvable Lie groups is obtained in \cite{BauesI}. 
The importance of the \emph{algebraic hull functor} for arbitrary (that is
not necessarily affine crystallographic virtually) 
polycyclic groups is now particularly evident from \cite{BauesGrunewald}.}

\paragraph{Solvable crystallographic groups}
First some further definitions and introductory remarks. 
A long standing conjecture of Louis Auslander states that affine crystallographic 
groups are virtually solvable groups. See Milnor's paper \cite{Milnor} for 
an introduction to this problem, and related results. The conjecture is verified in certain 
special cases, in particular in low dimensions. No counter example is known. 
See \cite{Abels} for details. By a result of Mostow \cite{Mostow1}, a discrete solvable subgroup of a Lie group is a polycyclic group. Hence, a 
solvable affine crystallographic group is, in particular, a polycyclic group.\footnote{There is a rich and well developed theory of polycyclic groups, as
is documented in the book \cite{Segal}.  For recent developments 
see also \cite{BauesGrunewald}.} 
Virtually polycyclic affine crystallographic
groups are therefore the objects of our study. 

We add 
that there are a few obvious mild restrictions 
for an abstract finitely generated group
$\Delta$ to be realised as a crystallographic subgroup of $\Aff(V)$: 
Let $\Delta \leq \Aff(V)$ be  crystallographic. 
Then the quotient space $\Delta \lmod V$
is a manifold if and only if $\Delta$ is torsionfree. 
Moreover, by Selberg's lemma every affine crystallographic group
contains a torsionfree normal subgroup of finite index.
Since $V$ is contractible, this implies that, 
the \emph{virtual cohomological dimension} of $\Delta$ (see \cite{Brown} for
definition) satisfies $\vcd \Delta = \dim V$.  

\begin{lemma} Let $\Delta \leq \Aff(V)$ be an ACG. Then 
every finite normal subgroup of $\Delta$ is trivial.
\end{lemma}
\begin{proof}  Let $N$ be a finite normal subgroup of $\Delta$. 
Since $N$ is finite, there exist a fixed point $x \in V$ for $N$. 
The set of fixed points for $N$ is an affine subspace $H$ of $V$, 
and it contains the orbit $\Delta x$. Since $\Delta$ acts as an 
ACG on $H$, we have $\dim H = \vcd  \Delta = \dim V$.
Therefore $H =V$, and $N= \{ 1 \}$. 
\end{proof}
A group which satisfies the conclusion of the lemma will be
called a \emph{\wfn-group}.

\subsection{Hulls, splittings and shadows}
We start by providing  some foundational material on algebraic and analytic hulls
for polycyclic groups. 
Let us first recall the definition of polycyclic groups. Namely, a group $\Gamma$
is called {\em polycyclic\/} if  $\Gamma$ admits a finite normal
series $$\Gamma= \Gamma_{0} \geq \Gamma_1 \geq \cdots \geq \Gamma_k= \{1\}$$ such that
each quotient $\Gamma_i/\Gamma_{i+1}$ is cyclic.
We let $\rank \Gamma$ denote its rank which is 
the number of 
infinite cyclic quotients $\Gamma_i/\Gamma_{i+1}$.   
(By some authors the rank of $\Gamma$ is called the Hirsch-length.)  
Then it is true that $\rank \Gamma = \vcd \Gamma$. The
rank of a virtually polycyclic group is defined as the rank of any finite
index polycyclic subgroup. 

\paragraph{The algebraic hull for a virtually polycyclic group}
Let $\bH$ be 
a linear  algebraic group. We put $\ur(\bH)$ for its maximal 
unipotent normal subgroup. The group $\ur(\bH)$ is Zariski-closed in $\bH$
and is called the {\em unipotent 
radical\/} of $\bH$.   
It is customary to  say that $\bH$ has a {\em strong unipotent radical\/} if
$\C_{\bH}( \ur(\bH) ) = \Z( \ur(\bH)) $.  

\begin{definition} \label{pchull}
A $\bbQ$-defined linear algebraic group $\bH$ with a strong unipotent radical is called an
{\em algebraic hull for $\Gamma$\/} if 
$\Gamma$ is a Zariski-dense subgroup of $\bH$,  $\Gamma \leq \bH_{ \bbQ}$, and 
$\dim \ur(\bH) = \rank  \Gamma$. 
\end{definition}
 
If $\Gamma$ is a virtually polycyclic \wfn-group then algebraic hulls for $\Gamma$ exist
and are unique up to $\bbQ$-isomorphism.  Also $\Gamma \cap \bH_{\bbZ}$ has finite index in $\Gamma$. This fact was proved by Mostow in the torsionfree polycyclic case, see also \cite[Proposition 4.40]{Raghunathan}. In section \ref{realproblems} of this article the case of finite extensions of \ftn-groups will be developed in detail. 
The more general construction for 
virtually polycyclic \wfn-groups  is given in  \cite{BauesI}. 
Further applications of the theory of algebraic hulls are developed
in  \cite{BauesGrunewald}. \\

Let us put $\bH_{\Gamma}$ for the 
algebraic hull of $\Gamma$. 
We will also need algebraic hulls over the real numbers.
\begin{definition} The group  $H_{\Gamma} = \bH_{\Gamma, \bbR}$ is called
a {\em real algebraic hull for $\Gamma$}.
\end{definition}

Note that $\Gamma$ is a discrete
subgroup in its real hull $H_{\Gamma}$ since
$\Gamma \cap \bH_{\bbZ}$ has finite index in $\Gamma$.
Let us henceforth write \ftn-group for a finitely generated torsionfree nilpotent group.   
We understand the algebraic
hull for a torsionfree polycyclic group $\Gamma$ as a substitute for 
the Malcev-completion\footnote{see \cite{Malcev,GO}} of a \ftn-group. 
In fact, if $\Gamma$ is a \ftn-group then
$\bH_{\Gamma}$ is the Malcev completion of $\Gamma$. 
Also the algebraic hull $\bH_{\Gamma}$ of a virtually polycyclic \wfn-group $\Gamma$ 
satisfies the following  rigidity property: 

\begin{proposition}  \label{pchullextension}
Let $\bH_{\Gamma}$ be an algebraic hull for\/ $\Gamma$,  
$\bG$ a $\bbQ$-defined linear algebraic group with a strong
unipotent radical.  Let $\rho: \Gamma \longrightarrow \bG$ be 
a homomorphism so that $\rho(\Gamma)$ is Zariski-dense in $\bG$.
Then $\rho$ extends uniquely to a homomorphism
$\rho_{\bH_{\Gamma}}: \bH_{\Gamma} \longrightarrow \bG$.
If $\rho(\Gamma) \leq   \bG_k$, where $k$ is a subfield
of\/ $\bbC$,  then $\rho_{\bH_{\Gamma}}$ is defined over $k$.   
  \end{proposition}
\begin{proof} We will use the diagonal argument. 
Therefore we consider  the subgroup 
$$ D= \{(\gamma, \rho(\gamma)) \mid \gamma \in \Gamma \} \,
 \leq \, \bH \times \bG \, \; .$$
Let $\pi_1, \pi_2$ denote the
projection morphisms on the factors of the product  $\bH \times \bG$.
Let $\bD$ be the Zariski-closure of $D$,  and $\bU = \ur(\bD)$ the unipotent
radical of $\bD$. Now $\bD$ is a solvable algebraic group, hence
$\bU = \bD_u$. Let $\alpha=  {\pi_1}|_{\bD}$. Since $\alpha$
is onto it follows that $\alpha$ maps  $\bU$ onto $\ur(\bH)$.
By \cite[Lemma 4.36]{Raghunathan} we have 
$\dim \bU \leq \rank \Gamma = \dim \ur(\bH)$,  and hence 
$\dim \bU = \dim \ur(\bH)$.  In particular it follows
that the restriction $\alpha:  \bU \rightarrow \ur(\bH)$ is an isomorphism. 
Thus the kernel of $\alpha$ consists only of semi-simple
elements.  Therefore every $x \in {\ker}\,  \alpha$ centralizes $\bU$, and
since $\pi_2(\bU)=\ur(\bG)$, $\pi_2(x)$ 
centralizes $\ur(\bG)$. Since 
$\bG$ has a strong unipotent radical, $x$ is in the kernel
of $\pi_2$. Hence $x =1$, and it follows that $\alpha$ 
is an isomorphism. We put $\rho_{\bH_{\Gamma}} = \pi_2 \circ \alpha^{-1}$
to get the required unique extension. 
If $\rho(\Gamma) \leq \bG_k$ then $\bD$ is $k$-defined, 
and hence also the morphism $\rho_{\bH_{\Gamma}}$. 
\end{proof}
\begin{remark} The proposition shows that the condition that the Zariski-closure
$\ac{\rho(\Gamma)}$ has
a strong unipotent radical forces the homomorphism $\rho$ to be well
behaved.  For example, $\rho$ must be unipotent on the Fitting subgroup of $\Gamma$.
See Proposition \ref{Fittu} below. 
\end{remark}
\begin{corollary} The algebraic hull\/ $\bH_{\Gamma}$ 
of\/ $\Gamma$ is unique up to $\bbQ$-isomorphism. 
In particular,  every automorphism $\phi$ of\/ $\Gamma$ extends uniquely to a
$\bbQ$-defined automorphism $\Phi$ of\/ $\bH_{\Gamma}$. 
\end{corollary}

Let us put $\Fitt(\Gamma)$
for the Fitting subgroup of $\Gamma$, that is, the 
unique maximal nilpotent normal subgroup of $\Gamma$. 
We note another property of the algebraic hull: 

\begin{proposition} \label{Fittu} 
Let $\Gamma \leq  \bH_{\Gamma}$. Then 
$\Fitt(\Gamma) \leq \ur(\bH_{\Gamma})$.
In particular, $\ur(\Gamma) = \Fitt(\Gamma)$. 
\end{proposition} 
\begin{proof} Let $F$ be the maximal 
nilpotent normal subgroup of $\bH_{\Gamma}$. 
Clearly, $F= \bF$ is a Zariski-closed subgroup.
Therefore  $\ur(\bF) = \ur(\bH_{\Gamma})$. 
Now since $\bF$ is nilpotent,  $\bF_s$ is 
a subgroup, $\bF_u = \ur(\bF)$ and 
$\bF = \bF_s \cdot\ur(\bF)$ is a
direct product of groups. 
Since 
$\bH_{\Gamma}$ has a strong unipotent radical, 
$\bF_s$ must be trivial, and it follows that $\bF= \ur(\bH_{\Gamma})$.
The Zariski-closure of $\Fitt(\Gamma)$ is a 
nilpotent normal subgroup of $\bH_{\Gamma}$,
and therefore $\Fitt(\Gamma)$ is contained in $\bF$,
hence $\Fitt(\Gamma) \leq \ur(\bH_{\Gamma})$. 
\end{proof} 

Recall that a {\em Cartan subgroup\/} of $\bH_{\Gamma}$ is
the centralizer of a maximal torus in $\bH_{\Gamma}^0$. Let
$\bC$ be a Cartan-subgroup  and $\bN = \ac{\Fitt(\Gamma)}$.
We just proved that $\bN \leq \ur(\bH_{\Gamma})$. 
\begin{lemma} \label{nilsupplement}
Any Cartan subgroup $\bC$ is nilpotent, $\bC = \bC_s \cdot \ur(\bC)$.
Moreover, $$\ur(\bH_{\Gamma}) =   \ur(\bC) \bN \; . $$
\end{lemma}  
\begin{proof} For the first statement, see \cite[11.7]{Borel}. 
Let $\bU =\ur(\bH_{\Gamma})$. 
Since $\bH_{\Gamma}$ is solvable, $(\bH_{\Gamma})_u = \bU$. 
Therefore  $\ur(\bC) =  \bC \cap \bU $. 
Since $[\, \Gamma, \Gamma] \leq \Fitt(\Gamma)$, 
$[\bH_{\Gamma}, \bH_{\Gamma}] \leq \bN$.
It is known (\cite[10.6]{Borel}) that all maximal tori are conjugate by elements of
$[\bH_{\Gamma}, \bH_{\Gamma}]$. Hence they are conjugate
by elements of $\bN$. Therefore, if $\bT$ is a maximal torus, 
$ \ur( \bH_{\Gamma}) =  \bN \C_ {\bU}(\bT) = \bN \ur(\bC)$. 
\end{proof}

\paragraph{The algebraic hull and the semisimple splitting of a solvable Lie group}
The notion of algebraic hull applies also to connected, 
simply connected solvable Lie groups $G$. We summarize, 
cf.\ \cite[Proposition 4.40]{Raghunathan}: 
\begin{proposition} There exists a linear algebraic group $\bH_ G$ with
a strong unipotent radical, so that $G \leq \bH_{\bbR}$ is a Zariski-dense
subgroup,  and $\dim \ur(\bH) = \dim G$.
\end{proposition}  

As for discrete groups, the algebraic hull $\bH_ G$ is unique 
up to isomorphism defined over $\bbR$, and it has analogous rigidity
properties. Let $N$ be the nilpotent radical of $G$, i.e., the maximal,  connected
nilpotent normal subgroup of $G$. Then it follows (with the same proof as
for Proposition \ref{Fittu}) that $N \leq \ur(\bH_ G)$. Note that 
$N$ is the connected component of $\ur(G)=  G \cap  \ur(\bH_ G)$.

\begin{definition} The group $H_G = \bH_ {G, \bbR}$ is called the 
{\em real algebraic hull\/} of $G$, and $U_G = \ur(H_G)$ 
is called the {\em unipotent shadow\/} of $G$. If  $N  = \ur(G)$ then $G$
is called {\em $u$-connected\/}. 
\end{definition} 
We remark further that $G$ is a normal subgroup of  $H_G$. In fact, 
$N \leq U_G$ is Zariski-closed in $H_G$, and $[G,G] \leq N$ implies therefore
that $[H_G, H_G] \leq N$. The existence of the algebraic hull for $G$ allows
for the following:
\begin{definition} A subset $X \leq G$ is called Zariski-dense in $G$, if $X$
is Zariski-dense in the algebraic hull of $G$.
\end{definition}
In particular, a subgroup $\Gamma \leq G$ is Zariski-dense in $G$ if it is so
in $H_G$.  We consider next the (real) 
{\em semisimple splitting\/} construction. (Compare \cite[Chapter 2,\S 3.6]{EMSII}.)
This idea goes back to Malcev's and Auslander's work on solvmanifolds.
\begin{definition} Let $M= S \cdot  U$ be a splittable,  simply connected solvable Lie group, so that
$U$ is the nilradical of $M$, and $S$ acts by semisimple automorphisms. 
$M$ is called a {\em semisimple splitting for $G$}, if $G$ is a normal subgroup of $M$,
so that $M= S \cdot G$ (semidirect product), and $M = G N$ (product of subgroups).   
\end{definition}
Note that it is said that
a group of automorphisms of $U$ is semisimple if the
corresponding induced
automorphism group of the Lie algebra acts by 
semisimple linear maps. 
Auslander (see \cite{AuslanderI}) proved that semisimple splittings
for $G$ exist and are unique. We briefly show the existence of $M$ by
realizing it inside the algebraic hull $H_G$.
\begin{proposition} \label{Msplitting}
There exists a compatible embedding of the 
semisimple splitting $M$  into the real algebraic hull $H_G$,
so that $U$ coincides with the unipotent shadow $U_G$. 
\end{proposition}
\begin{proof} Since $H_G$ is (real) algebraic and Zariski-connected, 
there exists a semidirect product decomposition $H_G = T  \cdot U_G$, 
where $T$ is a maximal torus.
We denote
$\psi: H \rightarrow U_G$ the projection map on $U_G$ which
is defined by the splitting. 
Let $C$ be a Cartan subgroup of $G$, so that $G= C N$, 
where $N \leq U_G$ is the nilradical of $G$. 
We put $S= C_s = \{ g_s \mid g \in C \}$,  so that $C \leq S \C_u$. 
$S$ is a Zariski-connected abelian subgroup of $H_G$, centralized
by $C$. By the conjugacy of maximal tori, we may assume that $S \leq T$.  
Since $H_{G} = \ac{G} \leq T \C_u N$, we conclude
that $U_G = C_u N$. 
Therefore  $U_G \leq S G$,  and also $G \leq S U_G$. 
It follows that the crossed homomorphism $\psi: G  \rightarrow U_G$ 
is surjective, in fact since $\dim U = \dim G$ it is a covering.
Since $U$ is simply
connected $\psi$ is a diffeomorphism. 
Therefore $S \cap G = \{1\}$ and we put
$M = S \cdot U_G= S \cdot G=  S C_u F = G U_G$.
Note that $S$ acts faithfully  
on the strong unipotent radical $U_G$.   
Therefore $M  \leq H_G$ is a semisimple splitting for $G$. 
\end{proof} 

\paragraph{Syndetic hulls and standard groups}
The notion of {\em syndetic hull\/} of a solvable subgroup
of a linear group is due to Fried and Goldman, cf.\  \cite[\S 1.6 ]{FriedGoldman},
and is an important tool in the study of discrete solvable groups, compare
also \cite[\S 5]{Witte}. Fried and Goldman introduced this notion in the
context of affine crystallographic groups.  
We use the slightly modified definition for the syndetic hull which is 
given in \cite{GS}. We then carry out some known 
constructions for ACGs in the algebraic hull of a general 
virtually polycyclic \wfn-group.

\begin{definition} \label{syndetic} 
Let $\Gamma$ be a polycyclic subgroup of $\GL(V)$, and
$G$ a closed, connected subgroup of $\GL(V)$. $G$ is called a
{\em syndetic hull\/} of $\Gamma$ if $\Gamma$ is a Zariski-dense 
uniform lattice in $G$, and $\dim G = \rank \Gamma$. 
\end{definition} 

We remark that a syndetic hull for $\Gamma$ is necessarily 
a connected, simply connected solvable Lie group. (One has to show that $G$ has no compact
subgroups.) 
The construction of a syndetic hull for $\Gamma$ may take place in
any (real) linear algebraic group $H$ which contains $\Gamma$,
provided $H$ satisfies certain conditions. 
\begin{definition}
A polycyclic group $\Gamma \leq H$ is 
called {\em standard in $H$\/} if $\Gamma$ is Zariski-dense, $\Gamma \leq H_{0}$,
$\Gamma' \leq \ur(H)$, and $\Gamma \rmod \ur(\Gamma)$ is torsionfree. 
$\Gamma$ is called a {\em standard polycyclic group\/} if $\Gamma$ is standard
in its real algebraic hull $H_{\Gamma}$. 
\end{definition}
Note that for $\Gamma \leq H_{\Gamma}$ to be standard it is enough to assume
that $\Gamma \leq (H_{{\Gamma}})_{0}$,  
and that $\Gamma \rmod \ur(\Gamma)$ is torsionfree. We recall
from Proposition \ref{Fittu} that for 
the embedding $\Gamma \leq H_{\Gamma}$ we have $\Fitt(\Gamma) =  \ur(\Gamma)$.
So for a standard polycyclic group $\Gamma \rmod \Fitt(\Gamma)$ is 
torsionfree.   

\begin{proposition} \label{syndeticH}
If\/  $\Gamma$ is standard then $\Gamma$ has a $u$-connected syndetic
hull $G$ in its real algebraic hull $H_{\Gamma}$.
\end{proposition}
\begin{proof} The proposition is in fact true if $\Gamma$ is a discrete standard subgroup 
in a real algebraic group $H$. Since $\Gamma \leq \bH_{\Gamma,\bbZ}$ has finite
index in $\Gamma$,  $\Gamma$  is a discrete
subgroup of $H_{\Gamma} = \bH_{\Gamma,\bbR}$. 
Now $\Gamma \leq H_{\Gamma}$ satisfies all the assumptions
which are needed to carry out the construction given in 
\cite[Proposition 4.1, Lemma 4.2]{GS}.
\end{proof}
We remark that {\em $\Gamma$ is a standard polycyclic group if and
only if\/ $\Gamma$ may be realized as a Zariski-dense lattice in 
a $u$-connected simply connected solvable Lie group $G$}. One half of
this statement is implied by the previous proposition. 
The other is the following proposition.    
\begin{proposition} \label{algebraicH}
Let $G$ be a connected simply  connected 
solvable Lie group, and\/ $\Gamma \leq G$ a Zariski-dense lattice. Then
the real algebraic hull $H_{G}$ is a real algebraic hull for $\Gamma$. 
If $G$ is $u$-connected (in $H_G$) then $\Gamma$ is standard. 
\end{proposition} 
\begin{proof} The Zariski-denseness of $\Gamma$ in $G$ implies by definition
that $\Gamma$ is Zariski-dense in $\bH_G$. Moreover,  $\rank \Gamma = \dim G$
since $\Gamma$ is cocompact.
Therefore $\bH_G$ is a $\bbR$-defined algebraic hull for $\Gamma$. 
By an application of Propositition \ref{pchullextension}, $\bH_G$ is isomorphic
over $\bbR$ to $\bH_{\Gamma}$. In particular, $H_G$ is isomorphic to $H_\Gamma$.
Since $\Gamma$ is contained in the  connected group $G \leq H_{\Gamma}$, 
$\Gamma \leq (H_{\Gamma})_0$.  

Let $N$ be the nilpotent radical of $G$.
It is known that $\Gamma /  (\Gamma \cap N) \leq A= G/N$ is torsionfree.  
Now if $G$ is $u$-connected, $N =  G \cap \ur(H_{\Gamma})$, and by Proposition \ref{Fittu} 
$\ur(\Gamma) = \Gamma \cap \ur(H_{\Gamma})  =  \Gamma \cap N$. 
We conclude that $\Gamma \rmod \ur(\Gamma) $ is torsionfree. 
In particular,  $\Gamma$ is standard. 
\end{proof} 
 
The following fact is easy to see.
\begin{proposition}  \label{hasStandard}
Let $\Delta$ be a virtually polycyclic group.
Then $\Delta$ has a normal polycyclic subgroup $\Gamma$ of finite
index which is standard. 
\end{proposition}

\paragraph{Discrete shadows}     
The semisimple splitting construction for a discrete subgroup $\Gamma$
of a connected, simply connected solvable Lie group $G$ was introduced
in the study of solvmanifolds. In particular, as Auslander
showed (see  \cite{AuslanderI}), the construction allows
to associate to $\Gamma$ certain lattices in the nilshadow of $G$. 
A semisimple  splitting for certain polycyclic
groups $\Gamma$ may be obtained with the help of {\em nilpotent supplements}.
(See for example \cite{Segal}.)
Maximal nilpotent supplements provide an analogue for the Cartan-subgroups 
in a solvable Lie-group $G$. 
The use of nilpotent supplements allows to use the Jordan-decomposition in an algebraic
group to construct the splitting. 
At this point, we only need the most elementary features of this 
construction. 
\newline  

Let $\Gamma \leq \bG$ be a torsionfree polycyclic subgroup of some linear algebraic group $\bG$. We further assume that $\Gamma = C F$,  where $F$ is a normal
subgroup of $\Gamma$, and $C$ is a nilpotent subgroup. $C$ will be called a 
{\em nilpotent supplement for $F$}. Let  $h: C \rightarrow \bG$ be
a homomorphism, where $C$ is nilpotent. 
For $\gamma \in C$, we put $h_u(\gamma) = h(\gamma)_u$.
This defines a homomorphism $h_u: C \rightarrow \bG$.
Now we associate with $\Gamma \leq \bG$ the group
$$ \Gamma^u_C=  \, < C_u,  F_u > \; . $$ 
We call $\Gamma^u_C$ the {\em  unipotent shadow of \/ $\Gamma$ in $\bG$ 
(with respect to $C$)}. 
\begin{lemma} The group
$\Gamma^u_C \leq \bG_{\bbQ}$ is a finitely generated, Zariski-dense 
subgroup of $\ur(\ac{\Gamma}) \leq\bG$.   
\end{lemma}
\begin{proof} By nilpotency,  we get from the Jordan-decomposition $C \leq C_s C_u$
and $F \leq F_s F_u$. $F_s$ is a central subgroup of $F_s F_u$ and consists of
semisimple elements. Therefore $\ac{F} = \ac{F}_s$. 
Since $C$ normalizes $F$ it normalizes $\bar{F}$. Now
$\ac{\Gamma} = \ac{C} \, \ac{F} = (\ac{C})_s  (\ac{F})_s  (\ac{C})_u (\ac{F})_u$.
Since $\ac{\Gamma}$ is solvable, $\ur(\ac{\Gamma}) =  \ac{\Gamma}_u$.
Therefore  $\ur( \ac{\Gamma}) =(\ac{C})_u (\ac{F})_u =  \ac{C_u}\,  \ac{F_u}$, and clearly 
$\Gamma^u_C=  \, < C_u,  F_u >$ is finitely generated and Zariski-dense.   
\end{proof}

Usually,  we consider $\Gamma \leq \bH_{\Gamma}$ as a subgroup of 
its algebraic hull, and $\Gamma^u_C \leq  \bH_{\Gamma}$ will be called
{\em the\/} unipotent shadow of $\Gamma$.

\begin{proposition} \label{Cshadow}
Let $\Gamma = C F$ be a torsionfree polycyclic group, where 
$F =\Fitt(\Gamma)$ and $C$ is a nilpotent supplement for $F$. Then
$$   \ur(\bH_{\Gamma})  =   (\ac{C})_u \ac{F} \; . $$
The unipotent shadow $\Gamma^u_C \leq \ur(\bH_{\Gamma})_{\bbQ}$ 
is a Zariski-dense subgroup of $\ur(\bH_{{\Gamma}})$. In particular, 
$\rank \Gamma = \rank \Gamma^U_C$.
\end{proposition}
\begin{proof} Note first that by the properties of the
Jordan decomposition $\Gamma \leq \bH_{\Gamma, \bbQ}$ 
implies that $\Gamma^u_C \leq  \ur(\bH_{\Gamma})_{\bbQ}$. 
Since $\ac{C}$ is nilpotent, $\ac{C} = S\,  \ac{C}_u$, where $S =  \ac{C}_s$
centralizes the unipotent group $\ur(\ac{C}) = (\ac{C})_u$. 
By Proposition \ref{Fittu}, $\Fitt(\Gamma) \leq \ur(\bH_{\Gamma})$. 
Therefore $\ac{F}$ is a unipotent normal subgroup of  $\ac{\Gamma}= \bH_{\Gamma}$.
By the proof of the preceding lemma,  $\ac{\Gamma} = \ac{C}_u \ac{F}_u =  \ac{C}_u \ac{F}$, 
and hence the first part of the proposition follows. Since $\Gamma_C^u$ is Zariski-dense
in $\ur(\bH_{\Gamma})$ it is known  (c.f.\ \cite[Theorem 2.10]{Raghunathan}) that 
$\rank \Gamma = \dim  \ur(\bH_{\Gamma}) \leq \rank  \Gamma_C^u$.
On the other hand for any \ftn-group $\Theta \leq  \ur(\bH_{\Gamma})_{\bbQ}$, 
$\rank \Theta \leq \dim  \ur(\bH_{\Gamma})$. 
Hence,  $\rank \Gamma = \rank \Gamma_C^u$.
 \end{proof}

Put $F= \Fitt(\Gamma)$. Let $\bG$ be a linear algebraic group. 
\begin{definition} \label{udescend} 
We say that a homomorphism 
$\rho: \Gamma \rightarrow \bG$ {\em descends to the unipotent shadow $\Gamma^u_C$\/} 
if there is a homomorphism $\rho^u: \Gamma^u_C \rightarrow \bG$ so that
$ \rho^u|_F = \rho_u$, and 
$\rho^u|_{C_u}(\gamma_u) = \rho_u(\gamma)$, for all $\gamma \in C$. 
\end{definition}
In this sense, the unipotent shadow $\Gamma^u_C \leq \bH_{\Gamma}$
has a rather strong functorial property:

\begin{proposition} \label{shadowdescent}
Let $\rho: \Gamma \rightarrow \bG$ be a homomorphism. Then 
$\rho$ descends to a  unique homomorphism 
$\rho^u: \Gamma^u_C \rightarrow \bG$. 
If $\rho(\Gamma) \leq \bG_k$, where $k$ is a subfield of $\bbC$, then
$\rho(\Gamma^u_C) \leq \bG_k$. 
\end{proposition}
\begin{proof} We may assume that $\rho(\Gamma)$ is Zariski-dense in $\bG$. 
As in the proof of Proposition \ref{pchullextension} we use the diagonal
construction for the homomorphism $\rho$. We proved there that, since
$\dim \ur(\bH_{\Gamma}) = \rank \Gamma$,  the projection map
$\alpha: \bD  \rightarrow \bH_{\Gamma}$ induces a $k$-defined isomorphism
$\alpha: \ur(\bD) \rightarrow \ur(\bH_{\Gamma})$ of the unipotent radicals. 
Since $\Gamma_C^u \leq \ur(\bH_{\Gamma})$ we can define
$\rho^u = \pi_2 \alpha^{-1}$. We have to  show that $\rho$
descends to $\rho^u$ in the above sense. Since $F \leq \ur(\bH_{\Gamma})$, 
it is enough  to verify that $\pi_2 \alpha^{-1} (\gamma_u) = \rho(\gamma)_u$, 
for all $\gamma \in \Gamma$. But, since 
$d= (\gamma, \rho(\gamma))_u = (\gamma_u, \rho(\gamma)_u) \in \ur(\bD)$
it is clear that $\alpha^{-1}(\gamma_u) = d$, 
and $\pi_2(\alpha^{-1}(\gamma_u)) = \rho(\gamma)_u$.
The rationality statement follows from well known properties
of the Jordan-decomposition. 
\end{proof}
Proposition \ref{shadowdescent} implies that 
the correspondence $\rho \longmapsto \rho^u$ defines a map
$$  s_u: \Hom(\Gamma, \bG) \longrightarrow \Hom(\Gamma^u_C,\bG) \; . $$
\begin{definition}  
We call $s_u: \Hom(\Gamma, \bG) \rightarrow \Hom(\Gamma^u_C, \bG)$ the
{\em shadow map for\/ $\Gamma$ (and $\bG$)}. 
\end{definition}

We will need the following result (compare \cite[\S7, Theorem 2]{Segal}):

\begin{proposition} \label{supplexists}
If  $\Delta$ is a polycyclic group then $\Delta$ 
has a characteristic  subgroup $\Gamma$ of finite index which admits a nilpotent
supplement $C$ for $\Fitt(\Gamma)$. 
\end{proposition}


Let $\Delta \leq \bH_{\Delta}$ be a virtually polycyclic \wfn-group, and $\Gamma \leq \Delta$
a characteristic subgroup of finite index which admits a nilpotent supplement $C$
for $\Fitt(\Gamma)$. Since $\Fitt(\Gamma) \leq F= \Fitt(\Delta)$, the group
$C \, \Fitt(\Delta)$ is of finite index in $\Delta$. $C$ is then an almost nilpotent
supplement for $\Delta$. In any case, we can associate to $\Delta$ the
group $$     \Delta_C^u  \; =  \;  < C_u,  F_u >  \, \leq \; \bH_{\Delta} \; , $$ and
we call it a {\em unipotent shadow for $\Delta$}.  

\subsection{Crystallographic groups and simply transitive groups}
In this section,  we collect, and sometimes reformulate,  basic 
results about virtually polycyclic ACGs. Our basic references are
\cite{FriedGoldman} and \cite{GS}. First, we explain the 
role of simply transitive groups. 

\subsubsection{Simply transitive hulls}

The importance of syndetic hulls in our setting comes from the following
central result,  (compare \cite{FriedGoldman},\cite[\S 4]{GS}): 

\begin{theorem} \label{sthull}
A virtually polycyclic subgroup $\Delta \subset \Aff(V)$ 
is an ACG if and only if there exists a simply transitive subgroup $G \subset \Aff(V)$
and $\Gamma \leq_f \Delta$ such that $G$ is a syndetic hull for $\Gamma$.
\end{theorem} 
The theorem describes the link between virtually polycyclic ACGs 
and certain simply transitive subgroups $G \leq \Aff(V)$.
We will call $G$ a {\em simply transitive hull\/} for 
the ACG $\Gamma$. We remark that the drawback in
considering the syndetic hull $G$ is that $G$ 
is, in general,  not uniquely determined by $\Gamma$.

\subsubsection{Extension property of simply transitive groups}
The next proposition was noticed in \cite{FriedGoldman}.
Our aim in this section is to provide a converse. 

\begin{proposition} \label{stclosure}
Let $G \leq \Aff(V)$ be a simply transitive group of
affine motions. Then the embedding of $G$ into its algebraic closure 
$\ac{G}\leq \Aff(V) $ is a real algebraic hull for $G$. 
\end{proposition} 
\begin{proof} Let $U$ denote the unipotent radical of $\ac{G}$.
Auslander \cite[\S3]{Auslander} proved that $U$
is simply transitive on $V$. Hence $\dim U = \dim V= \dim G$. 

$\C_{\ac{G}}(U)$ is an algebraic subgroup of $\ac{G}$.
Since $U$ is transitive, the elements of $\C_{\ac{G}}(U)$ act 
without fixed points on $V$.  Therefore $\C_{\ac{G}}(U)$ contains 
no semisimple elements. Hence  $\C_{\ac{G}}(U)$ is a Zariski-closed 
normal unipotent subgroup of $\ac{G}$. Hence  $\C_{\ac{G}}(U) \leq \ur(\ac{G})=U$.
So $\ac{G}$ has a strong unipotent radical.  
\end{proof}

Let us call a homomorphism $i: H \rightarrow \Aff(V)$ of real linear algebraic
groups {\em $u$-simply transitive\/} if the unipotent radical $\ur(i(H))$ acts
simply transitively on $V$. We can characterise
simply transitive groups now as follows\footnote{Auslander \cite{Auslander} observed
that the nilpotent shadow of $G$ acts simply transitively}: 
\begin{theorem} \label{stextension} 
A connected, simply connected, 
solvable Lie  subgroup $G \leq \Aff(V)$ is simply transitive 
on $V$ if and only if the embedding of $G$ into $\Aff(V)$ extends to a
$u$-simply transitive embedding of algebraic groups $H_G \rightarrow \Aff(V)$.
\end{theorem} 
\begin{proof}  The uniqueness of the real algebraic hull 
implies that the  ``only if'' part of the theorem is just 
the previous proposition. So let us 
prove that $G$ is simply transitive if  $\ac{G} \leq \Aff(V)$
is a $u$-simply transitive algebraic hull for $G$. By Proposition
\ref{Msplitting}, there exists a torus $S \leq \ac{G}$ such
that $U \leq G S$. The torus $S$ has a fixed point $x \in V$,
so that $G x \supset Ux =V$. Therefore $G$ acts transitively,
and, since $\dim G = \dim U = \dim V$, $G$ acts a fortiori 
simply transitively on $V$. 
 \end{proof}
%
%

\subsubsection{Extension property of crystallographic groups} \label{ACfoundations}
We come now to the main result of this subsection.\footnote{the result was announced in \cite[Theorem 3.2]{BauesU}} 
\begin{theorem} \label{acextension} 
Let $\Delta \leq \Aff(V)$ be a virtually polycyclic \wfn-group.
Then $\Delta$ is an ACG if and only
if the embedding of $\Delta$ into $\Aff(V)$ extends to a
 $u$-simply transitive embedding of algebraic groups $H_{\Delta} \rightarrow \Aff(V)$. 
\end{theorem}
\begin{proof} Let us assume that $\Delta$ is an ACG.   
By Theorem \ref{sthull},  $\ac{\Delta}$ contains
a simply transitive hull $G$ for a finite index subgroup $\Gamma$ of
$\Delta$. Since $G$ is a syndetic hull for  $\Gamma$, 
the algebraic closure of $\Gamma$
coincides with  $\ac{G}$, and $\rank \Gamma = \dim G$.
Therefore, by  Theorem \ref{stextension}, $\dim \ur(\ac{\Delta}) = \rank \Gamma$. 
Since $\ac{G}$ is $u$-simply transitive, $\ac{\Delta}$ is $u$-simply transitive. 
Since $\ac{\Delta}^0 = \ac{G}$ has a strong unipotent radical, we can conclude that 
$\ac{\Delta}$ is a real algebraic hull for $\Delta$.

For the converse, let us assume that $\Gamma \leq \Aff(V)$ is
a subgroup such that $\ac{\Gamma}$ is a $u$-simply transitive 
real algebraic hull for $\Gamma$. Since $\Gamma$ has a
standard subgroup of finite index, 
we may as well assume that $\Gamma$ is standard.
Now if $\Gamma$ is standard, 
there exists, by Proposition \ref{syndeticH}, a syndetic hull $G \leq \ac{\Gamma}$ for $\Gamma$.
By Proposition \ref{algebraicH}, $\ac{G}=\ac{\Gamma}$ is an algebraic hull
for $G$, and  $u$-simply transitive by assumption.
By Theorem \ref{stextension}, $G$ acts 
simply transitively on $V$. In particular,  $\Gamma$ is an ACG. 
\end{proof}

\begin{remark} The theorem, together with Proposition \ref{Fittu}, implies that
any ACG $\Gamma \leq \Aff(V)$ satisfies $\Fitt(\Gamma)= \ur(\Gamma)$. 
This is the content of Lemma C in \cite{GS}. 
\end{remark}

\section{The variety of crystallographic homomorphisms}

Let $\Gamma$ be a virtually polycyclic \wfn-group, 
and  $A \leq \Aff(V)$ a Zariski closed subgroup. The purpose
of this section is to study the space of crystallographic
homomorphisms 
$$\Hom_c(\Gamma,A) \, \subset \, \Hom(\Gamma,A) $$ as a subset 
of the real algebraic variety $\Hom(\Gamma, A)$. 
The main result of this section,  Theorem \ref{IsZariskiopen}, establishes 
that the space $\Hom_c(\Gamma,A)$ 
is described by algebraic equalities and inequalities, and
carries itself a natural structure as a real algebraic variety
which is defined over the rational numbers.
We first prove the result for finitely generated torsionfree nilpotent groups.
Our strategy is to use the unipotent shadow construction
to extend from finitely generated torsionfree nilpotent groups 
to general virtually polycyclic groups.  

\subsection{Spaces of homomorphisms} \label{topology}
Let us introduce here the spaces which we want to study. 
We start with some preliminary remarks. 

\paragraph{Topology on the space of homomorphisms}
Let $H,G$ be locally compact groups, and $R(H,G)$ the
space of all continuous  homomorphisms from $H$ to $G$. 
We equip the space $R(H, G)$ with the {\em compact open topology}, 
meaning that a fundamental system of neighbourhoods
in $R(H, G)$ is specified by all sets of the form
$$   {\cal U}_{K,U} =  \{  \, \rho \in  R(H, G) \mid \rho(K) \subset U \,  \}\:  ,$$ 
where $K\subset H$ is compact,  and $U \subset G$ is open. 

Let us assume next that $\Gamma$ is a discrete group.  We equip  $\Hom(\Gamma, G)$
with the subspace topology which is inherited from the {\em product topology\/} on
$$ {\rm Map}(\Gamma, G) =  G^{\Gamma} \; .  $$
It is easy to see that this topology coincides with the  compact open topology
on  $\Hom(\Gamma, G)= R(\Gamma, G)$. 
If $\Gamma$ is finitely generated, and 
$S = \{ \gamma_1, \ldots \gamma_n \}$ is a system of generators, then
$\Hom(\Gamma, G)$ is a subset of  $G^n$ in a natural way.  
In fact, the inclusion $$ \jmath: \Hom(\Gamma, G) \longrightarrow G^n$$ which
is  given by $\rho \mapsto (\rho(\gamma_1), \ldots, \rho(\gamma_n))$ is a
homeomorphism onto a closed subspace of $G^n$, see for example
 \cite[2.23]{MacSing}. Therefore,  if $\Gamma$ is finitely generated, 
$\Hom(\Gamma, G)$ carries also the subspace topology from $G^n$. 

If $G= \bG_{\bbR}$ is a real (linear) algebraic group, $G$ carries the locally compact
Hausdorff-topology and the Zariski-topology. So does $G^n$, and
the embedding  $\jmath$ identifies $\Hom(\Gamma, G)$ with a Zariski-closed
subset. Aside from the Hausdorff-topology, the space $\Hom(\Gamma, G)$ carries 
a natural structure of
a real algebraic variety which is independent of the embedding $\jmath$.
This follows from: 
\begin{proposition} Let $\Gamma$ be a finitely generated group, and $\bG$ 
a $\bbQ$-defined linear algebraic group. 
Then $\Hom(\Gamma, \bG)$ has a natural structure of 
affine algebraic variety defined over $\bbQ$ and 
$\Hom(\Gamma, \bG_{k}) =  \Hom(\Gamma, \bG)_k$ for 
any subfield $k \leq \bbC$. 
\end{proposition} 
See \cite{LM} for a proof. The space $\Hom(\Gamma, \bG)$
is called a representation variety. Note that $\bG$ acts 
on homomorphisms by conjugation. This turns $\Hom(\Gamma, \bG)$
into a $\bG$-variety. 

\paragraph{Crystallographic homomorphisms}
Recall that a homomorphism $\rho: \Gamma \rightarrow \Aff(V)$ is
an {\em affine crystallographic homomorphism}, if $\rho$ is 
an isomorphism onto its image $\rho(\Gamma)$ and $\rho(\Gamma)$
is a crystallographic subgroup of $\Aff(V)$.   We put
$$   \Hom_c(\Gamma, \Aff(V)) \, = \, 
\{  \, \rho: \Gamma \rightarrow \Aff(V) \mid \rho \, \mbox{ is 
crystallographic} \}  $$ for the {\em space of crystallographic homomorphisms}. 
We consider then $\Hom_c(\Gamma, \Aff(V))$ as a Hausdorff-topological space with
the subspace topology inherited from the space $\Hom(\Gamma, \Aff(V))$. 
The space $\Hom_c(\Gamma, \Aff(V))$ has a natural action of $\Aff(V)$, 
which is induced by conjugation on homomorphisms.  


\subsection{Nilpotent crystallographic groups}
We prove here that for finitely generated
torsionfree nilpotent groups
the space of crystallographic homomorphisms is a
Zariski-open subset of the space of unipotent homomorphisms. 
This result was developed in \cite{BauesV}. 


\paragraph{Nilpotent crystallographic groups} 
Let $\Gamma$ be a \ftn-group. Our aim is to describe the inclusion 
$$
  \Hom_c(\Gamma, \Aff(V)) \, \subset \, \Hom(\Gamma, \Aff(V)) \; .
$$
\begin{remark} Using the fact that rank and cohomological dimension
of $\Gamma$ coincide, it is evident that $ \Hom_c(\Gamma, \Aff(V))= \emptyset $,
for $\rank \Gamma \neq \dim V$. 
\end{remark}

Let us briefly specialise some facts of section \ref{ACfoundations}. Since
$\Gamma$ is a \ftn-group, the real algebraic hull $U_{\Gamma}$ is
just the real Malcev hull. 
The group $\Gamma$ is then a discrete, cocompact lattice in $U_{\Gamma}$.
The following extension property is a special case of Proposition
\ref{pchullextension}:

\begin{proposition} \label{uextension}
Let $\Gamma$ be a \ftn-group.  Then every homomorphism 
$\rho: \Gamma \rightarrow U$ into a unipotent  
Zariski-closed subgroup $U \leq \Aff(V)$ uniquely extends to a homomorphism 
of algebraic groups $\rho_{U_{\Gamma}}:  U_{\Gamma}  \rightarrow U$.
\end{proposition} 

From Theorem \ref{acextension} we deduce  the following characterization of 
crystallographic \ftn-groups: 

\begin{theorem} \label{Aus2}
Let $\Gamma \subset \Aff(V)$ be a \ftn-group. 
Then $\Gamma$ acts crystallographically on $V$ if and only if the algebraic
closure $\overline{\Gamma} \subset  \Aff(V)$ is a unipotent simply transitive 
(real) Malcev hull for $\Gamma$. 
\end{theorem}

The theorem will be useful together with 
\begin{lemma} \label{cishull}
Let $\Gamma \subset \Aff(V)$ be a \ftn-group, and $\rank \Gamma \leq \dim V$. 
If the algebraic
closure $U= \overline{\Gamma} \subset  \Aff(V)$ is a unipotent 
simply transitive group then $U$ is a real Malcev hull for $\Gamma$.
\end{lemma}
\begin{proof}
Let $U_{\Gamma}$ be the Malcev hull of $\Gamma$. Let 
$j: \Gamma \hookrightarrow \ac{\Gamma}=U$ be the inclusion homomorphism.
Now $\dim \ur(\ac{\Gamma}) \leq \rank \Gamma = \dim U_{\Gamma}$.
By Proposition \ref{uextension}, $j$ extends to a surjective homomorphism
$j_{U_{\Gamma}}: U_{\Gamma} \rightarrow U$ of unipotent algebraic groups. Since 
$\dim  U_{\Gamma} = \rank \Gamma = \dim V = \dim U$, 
$j_{U_{\Gamma}}$ is an isomorphism. 
\end{proof}

Before stating the main result, we introduce 
some further notation.
A function $\delta$ on a $G$-variety $V$,
where $G$ is a group,  is called a {\em relative $G$-invariant\/} if 
there exists a character $\chi$ of $G$, so that, for all $g \in G$ and $v \in V$, 
$$  \delta (g \cdot v) = \chi(g)\,  \delta(v) \; . $$ 
For any continuous function $\delta$, the set $V_{\delta}$ defined as 
$$   V_{\delta} = \{ v \in V \mid \delta(v) \neq 0 \} \;  $$  
is called a {\em special open subset\/} of $V$. \\

Let us put  
$$  \Hom_u(\Gamma, \Aff(V)) =  \{ \rho \mid \rho(\Gamma) 
     \mbox{ is unipotent } \} \, \subset \, \Hom(\Gamma,\Aff(V))  
$$  
for the variety of unipotent representations. 

\begin{lemma}  The space 
$\Hom_u(\Gamma, \Aff(V))$ is a Zariski-closed
in $\Hom(\Gamma,\Aff(V))$, and it is
defined over $\bbQ$.   
\end{lemma}
\begin{proof}
If $\Gamma \leq \GL(V)$ is nilpotent then
$\Gamma_s$ and $\Gamma_u$ are commuting subgroups,
$\Gamma \leq \Gamma_s \times \Gamma_u$
and the projection maps are homomorphisms. 
(See for example \cite[Chapter 7, Proposition 3]{Segal}.)
Therefore  $\rho(\Gamma)$ is unipotent if and only if $\rho$ is unipotent 
on a set of generators, and we conclude that
$\Hom_u(\Gamma, \Aff(V))$ is a Zariski-closed
subset of $\Hom(\Gamma,\Aff(V))$.   
\end{proof}
Theorem \ref{Aus2} implies that 
\begin{equation}    \label{cisunipotent}
\Hom_c(\Gamma, \Aff(V)) \,  \subset \,  \Hom_u(\Gamma, \Aff(V)) \; . 
\end{equation}
Surprisingly, $ \Hom_c(\Gamma, \Aff(V))$ is a special open subset of
$ \Hom_u(\Gamma, \Aff(V))$:

\begin{theorem} \label{IsZariskiopenN}  
Let $\Gamma$ be a \ftn-group which 
satisfies $\rank \Gamma = \dim V$.
Then there exists a $\bbQ$-defined polynomial function 
$\delta$ on the variety  of unipotent representations 
$\Hom_u(\Gamma, \Aff(V))$,  such that 
$$  
 \Hom_c(\Gamma, \Aff(V)) \, = \, \Hom_u(\Gamma, \Aff(V))_{ \delta} \; . \; \; 
$$ 
Moreover, $\delta$ is a relative invariant with respect to the natural
$\Aff(V)$-action on  $\Hom_u(\Gamma, \Aff(V))$. 
\end{theorem} 

\begin{remark} The character which belongs 
to the relative invariant function $\delta$ is the determinant
on $\Aff(V)$, cf.\  \cite{BauesV} .
\end{remark}
 
We need some more preparations for the proof of the theorem. 

\paragraph{The Malcev hull of a \ftn-group} 
Let  $\bU_{\Gamma}$ be the Malcev hull of $\Gamma$. This means 
that  $\bU_{\Gamma}$ is
a $\bbQ$-defined unipotent linear algebraic group, and 
$\Gamma \leq (\bU_{\Gamma})_{\bbQ}$ is a Zariski-dense subgroup. 
Let $\lie u$ be
the Lie algebra of $\bU_{\Gamma}$. The exponential map 
$$ \exp: {\lie u} \longrightarrow \bU_{\Gamma}$$  is a
polynomial map which is an isomorphism of algebraic 
varieties. The inverse map is $\log: \bU_{\Gamma} \rightarrow  {\lie u}$.
It is known that there
exists a $\bbQ$-structure on the complex Lie algebra $\lie{u}$,  
so that $\log \Gamma \leq {\lie u}_{\bbQ}$, and,  in fact,
${\lie u}_{\bbQ} = \bbQ \log \Gamma$ is the $\bbQ$-span
of $\log \Gamma$. Moreover,  $\exp$ and $\log$ are $\bbQ$-defined maps. 
(See \cite{GO} for more details.)

\begin{definition}
Let $S = \{ \gamma_1, \ldots, \gamma_n \}$ be a system of 
generators for $\Gamma$. $\Gamma$ is called a {\em Malcev-basis 
of $\Gamma$}, if  the series of
subgroups 
$$ 1  \, \leq  \; \,  < \! \gamma_1 \! > \; \leq \; \,  < \!\gamma_1, \gamma_2\! > \; \leq
\;   \cdots  \; \leq  \; \,  <\! \gamma_1, \,  \ldots, \gamma_n \!> = \Gamma$$
is a central series for $\Gamma$ with infinite cyclic factors. 
\end{definition} 
By results of Malcev
\cite{Malcev} , every
\ftn-group $\Gamma$ admits a Malcev-basis. \\

The following facts can be derived from the  
Baker-Campbell-Hausdorff formula: 
Put $({\lie u}_i)_{\bbQ} = \bbQ \log \Gamma_i$, where 
$\Gamma_i = \,  <\! \gamma_1, \,  \ldots, \gamma_i\!>$. 
The complex span ${\lie u}_i$ of $\log \Gamma_i$ is
an ideal in the Lie algebra $\lie u$ and the series 
$$ 0   \, \subset   \;   {\lie u}_1   \,  \subset   {\lie u}_2 \; 
\, \cdots  \, \subset  \;  {\lie u}$$
is a central series for $\lie{u}$ with one-dimensional factors.   
In particular, 
$${\lie u}_{\bbQ}  \, =  \, < \log \gamma_1, \, \ldots \, , \log \gamma_n >_{\bbQ} \; . $$ 

\paragraph{Simply transitive and \' etale unipotent actions}
We have to consider simply transitive unipotent actions.  
Let $U \leq \Aff(V)$ be a unipotent algebraic subgroup which
satisfies $\dim U =\dim V$.  
$U$ is called {\em \' etale on $V$} if it has an open orbit on $V$. 
By a result of Rosenlicht \cite{Rosenlicht} the orbits of the unipotent group $U$
on $V$ are all closed. Therefore $U$ is simply transitive
if and only if $U$ is \' etale. \newline 

\begin{prf}{Proof of Theorem \ref{IsZariskiopenN}}
Let  ${\lie a}(V)$ denote the Lie algebra
of  $\Aff(V)$, and let ${\lie a}_n(V)$ denote the (algebraic) subset of nilpotent elements
in  ${\lie a}(V)$. 
Recall that the exponential map 
$$ \exp: {\lie a}(V) \longrightarrow \Aff(V) $$ 
defines a polynomial map which induces a {\em polynomial\/} equivalence 
from ${\lie a}_n(V)$ onto the (algebraic) subset $\Aff_u(V)$ of unipotent elements
in  $\Aff(V)$.  On  $\Aff_u(V)$, the map $\exp$ has a well defined polynomial inverse
$$  \log:  \Aff_u(V) \longrightarrow  {\lie a}_n(V) \; . $$ 

Let $S= \{ \gamma_1, \ldots, \gamma_n \}$ be 
a Malcev-basis for $\Gamma$, where $n =\dim V$. 
If  $\rho \in \Hom_u(\Gamma, \Aff(V))$, then $\rho(\Gamma) \leq \Aff(V)$ 
is a unipotent $\ftn$-group. The Lie algebra ${\lie u}_{\rho}$ of $U_{\rho} =\ac{\rho(\Gamma)}$
is therefore contained in ${\lie a}_n(V)$. Since  ${\lie u}_{\rho}$ is a homomorphic image
of $\lie{u}_{\Gamma}$ it is spanned by the set $\{ \log \rho(\gamma_i) \mid  i=1 \ldots n \}$.
It follows from Proposition \ref{Aus2}, together with 
Lemma \ref{cishull}, that, 
$\rho \in  \Hom_c(\Gamma,\Aff(V))$ if and only if 
the algebraic closure $U_{\rho}$ 
is a simply transitive unipotent subgroup of $\Aff(V)$.  

Choose an arbitrary base point $x \in V$. The differential of the orbit map 
$A \mapsto A \cdot x$ of the action of $\Aff(V)$ on $V$ defines a linear
map $o_x: {\lie a}(V) \rightarrow V$. Therefore $U_{\rho} \leq \Aff(V)$ is \' etale 
in $x$ if and only if the restriction of $o_x$ 
to the Lie-algebra ${\lie u}_{\rho} \subset  {\lie a}(V) $ 
is an isomorphism of vector spaces. It follows a fortiori that
$U$ is simply transitive on $V$ if and only if the restriction 
 $o_x: \lie{u_{\rho}} \rightarrow V$ is an isomorphism.

Let us next  consider the linear map $\tau_x(\rho): \bbR^n  \rightarrow \bbR^n$ which is 
defined by 
$$   \tau_x(\rho): \;  (\alpha_1, \ldots , \alpha_n) \longmapsto  
         o_x(\;  \sum_{i=1}^n \alpha_i \log \rho(\gamma_i) \, ) \;.
$$
(After 
choosing an arbitrary basis in $V$, we may as well view $o_x$ as
a map with values in $\bbR^n$.)
We easily see that $\tau_x(\rho)$ is an isomorphism 
if and only if $o_x: {\lie{u}_{\rho}} \rightarrow V$ is.
We therefore define $$ \delta(\rho) = \det \tau_x(\rho) \; , $$ and,  clearly,  $\delta$ is
a polynomial function in the matrix entries of the $\rho(\gamma_i)$. Therefore
$\delta$ is a polynomial  on $\Hom_u(\Gamma, \Aff(V))$ with the
property that $\delta(\rho) \neq 0$ if and only if 
$U_{\rho} = \overline{\rho(\Gamma)} \leq \Aff(V)$
is a simply transitive subgroup. In this case, the extension
$\rho_{U_{\Gamma}}: U_{\Gamma} \rightarrow U_{\rho}$ 
is an isomorphism of algebraic groups, and  $\rho(\Gamma)$ 
is a crystallographic group. Up to now we proved  that  
$$ \Hom_c(\Gamma, \Aff(V)) =  \Hom_u(\Gamma, \Aff(V))_{\delta} \; . $$
Next we show that the function $\delta$ is independent of
$x \in V$. To see this, we remark that our previous reasoning
is also valid over the field of complex numbers. In 
particular if $\delta(\rho) \neq 0$ the corresponding 
unipotent algebraic group $\bU_{\Gamma}$ acts simply
transitively on $\bbC^n$. If $\rho$ is fixed and 
$\delta(\rho) \neq 0$ then the function
$$  x \mapsto \det \tau_x(\rho) $$
is a polynomial on $\bbC^n$ which does not vanish. 
Hence $\det \tau_x(\rho)$ is constant.
To show that $\delta(\rho) = \det \tau_x(\rho)$ is 
a relative invariant on $\Hom_c(\Gamma,\Aff(V))$ 
we remark that by direct calculation, for all $g \in \Aff(V)$, the formula 
$$  \det \tau_x(\rho^g) = (\det g) \det \tau_{g^{-1}x} (\rho) $$ 
holds. It follows therefore that
$\det \tau_x(\rho^g) = (\det g)   \det \tau_{x} (\rho)$.
\end{prf}

\begin{remark} We may also define $\tilde{\delta}(\rho)  = \delta(\rho_u)$. 
If  $\rho \in  \Hom(\Gamma, \Aff(V))$, such  that $\tilde{\delta}(\rho) \neq 0$, then, 
by the Theorem,  $\rho_u$ is in $\Hom_c(\Gamma, \Aff(V))$. In particular
$U = \ac{\rho_u(\Gamma)}$ is a simply transitive subgroup. 
But since $\rho(\Gamma)_s$
consists of semisimple elements, and is contained in the centralizer of $U$, 
it must be trivial.  Therefore $\rho =\rho_u$,
and a fortiori $\rho \in \Hom_c(\Gamma, \Aff(V))$. Hence, 
$$ \Hom_c(\Gamma, \Aff(V)) =  \Hom(\Gamma, \Aff(V))_{\tilde{\delta}} \; .$$ 
Note however that the function $\tilde{\delta}$ is {\em not\/} continuous on
$\Hom(\Gamma, \Aff(V))$. 
\end{remark}

\noindent The next example illustrates the theorem in the
simplest possible situation: 
\begin{example} Let $V = \bbR$. Then 
\begin{align*}   
\Hom(\bbZ, \Aff(\bbR)) &  =  \Aff(V) =  \left\{ 
\begin{matrix}{cc} 
\epsilon & v \\   0 & 1 
\end{matrix} \mid \epsilon \neq 0 \right\}  \,  \leq \, \GL(2, \bbR) \\
\Hom_{u}(\bbZ, \Aff(\bbR)) &  =   V = \left\{ 
\begin{matrix}{cc} 
1 & v \\   0 & 1 
\end{matrix}  \right\}  \\
\intertext{and} 
  \Hom_c(\bbZ, \Aff(\bbR)) &  = \; \left\{ 
\begin{matrix}{cc} 
1  & v \\   0 & 1 
\end{matrix} \mid  v \neq 0 \right\} \; \; . 
\end{align*}
\end{example}

\subsection{Virtually polycyclic crystallographic groups}
After some more preparations we come to the
statement and proof of  our main result on the space of
crystallographic homomorphisms, Theorem \ref{IsZariskiopen}.
 
\paragraph{The crystallographic shadow map} 
Let $\Gamma$ be a virtually  polycyclic \wfn-group. 
The main step in our description of the space 
$\Hom_c(\Gamma,A)$  is to characterise the 
elemens $\rho \in \Hom_c(\Gamma,A)$  by properties which may be verified 
on the unipotent shadow of $\Gamma$. 
Let us assume here that $\Gamma$ is polycyclic, and 
$\Gamma=C F \leq \Aff(V)$ admits a
nilpotent supplement for $F = \Fitt(\Gamma)$. Recall then
from Proposition  \ref{Cshadow} that there is a certain
unipotent group $\Gamma^u_C \leq \Aff(V)$, the unipotent
shadow of $\Gamma$ in $\Aff(V)$,  which is associated with $\Gamma$. 

\begin{proposition} \label{crystshadow}
Let\/ $\Gamma \leq \Aff(V)$  be a torsionfree polycyclic group which satisfies
$\rank \Gamma= \dim V$. Moreover, we 
assume that\/ $\Gamma = C F$ admits a nilpotent supplement\/
$C$ for\/ $F = \Fitt(\Gamma)$. Then $\Gamma$ is an ACG
if and only if the unipotent shadow\/ $\Gamma^u_C$ of\/ $\Gamma$ 
in $\Aff(V)$ is an ACG. 
\end{proposition}
\begin{proof}  
Let us first assume that $\Gamma$ is an ACG.
Then, by Theorem \ref{acextension}, the algebraic closure 
$\ac{\Gamma}$ is a $u$-simply transitive
real algebraic hull for $\Gamma$. By Proposition \ref{Cshadow}, 
$ \Gamma^u_C  = C_u F$
is a lattice in $U= \ur(\ac{\Gamma})$. Since $U$
is simply transitive, $\Gamma^u_C$  is
an ACG. 

Now let us assume that $ \Gamma^u_C  = C_u F_u$ is an 
unipotent ACG. It follows that $U =  \ac{ \Gamma^u_C} = \ac{C_u} \ac{F_u}$
is a unipotent simply transitive subgroup of $\ac{\Gamma}$. 
Since $\ac{\Gamma}$ is solvable, $U \leq \ur(\ac{\Gamma})$. 
By a lemma of Mostow (cf.\ \cite[4.36]{Raghunathan})
it is known that $\dim \ur(\ac{\Gamma}) \leq \rank \Gamma= \dim V$,
and hence it follows that $\ur(\ac{\Gamma}) =U$. 
So $\ac{\Gamma}$ is in fact $u$-simply transitive, and 
$\dim \ur(\ac{\Gamma}) =  \rank \Gamma$. 
The centralizer $\C_{\ac{\Gamma}}(U)$ 
of the simply transitive normal subgroup $U$
is a unipotent subgroup of $\ac{\Gamma}$. 
Therefore $\C_{\ac{\Gamma}}(U) \leq \ur(\ac{\Gamma}) = U$.
Hence $\ac{\Gamma}$ has a strong unipotent
radical, and is a $u$-simply transitive algebraic hull for $\Gamma$. 
By Theorem \ref{acextension},  $\Gamma$ is an ACG.   
\end{proof}

Recall from Definition \ref{udescend} that,  for $\Gamma = C F$ as above, 
a homomorphism $\rho: \Gamma \rightarrow \bG$ is said to
descend to  the unipotent shadow $\Gamma^u_C \leq H_{\Gamma}$ if
there is a homomorphism $\rho^u: \Gamma^u_C \rightarrow \bG$ so that
$ \rho^u|_F = \rho_u$, and 
$\rho^u|_{C_u}(\gamma_u) = \rho_u(\gamma)$, for all $\gamma \in C$.

\begin{proposition} \label{cshadowdescent}
Let\/ $\Gamma$  be a torsionfree  polycyclic group which satisfies 
$\rank \Gamma = \dim V$. We 
assume also that\/ $\Gamma = C F$ admits a nilpotent supplement\/
$C$ for\/ $F = \Fitt(\Gamma)$.  Let $\rho: \Gamma \rightarrow \Aff(V)$
be a homomorphism. Then $\rho$ is
a crystallographic homomorphism if and only if $\rho$ descends to a
crystallographic homomorphism $\rho^u$ on the unipotent  shadow $\Gamma^u_C$
of\/ $\Gamma$. 
\end{proposition}
\begin{proof}
Let us first assume that $\rho$ is crystallographic. 
By Theorem \ref{acextension},  $\rho$ extends to a $u$-simply transitive
morphism of (real) algebraic groups 
$\rho_{H_{\Gamma}}: H_{\Gamma} \hookrightarrow \Aff(V)$. 
As in the proof of Proposition \ref{shadowdescent}, $\rho$ descends 
to the restriction of  $\rho_{H_{\Gamma}}$ to 
$\Gamma^u_C \leq \ur(H_{\Gamma})$.
Since $\ur(\rho(H_{\Gamma}))$ is simply transitive, this homomorphism is crystallographic. 

Conversely,  let us assume that $\rho$ descends to a crystallographic
homomorphism $\rho^u$ on  $\Gamma^u_C \leq \ur(H_{\Gamma})$. Then,  
since $\rho^u(\Gamma^u_C)$ is a unipotent shadow for $\rho(\Gamma)$ in
$\Aff(V)$, Proposition \ref{crystshadow} implies that  $\rho(\Gamma)$ 
is crystallographic on $V$.
In particular, $ \rank  \rho(\Gamma)  =\dim V =\rank \Gamma$. It follows that
the kernel of $\rho$ is finite. Since $\Gamma$ is torsionfree, $\rho$
is injective, hence  $\rho$ is a  crystallographic homomorphism.  
\end{proof} 

Proposition \ref{cshadowdescent} implies that 
the correspondence $\rho \longmapsto \rho^u$ defines a map
$$  s_u: \Hom_c(\Gamma, A) \longrightarrow \Hom_c(\Gamma^u_C,A) \; . $$
\begin{definition}  
We call $s_u: \Hom_c(\Gamma, A) \rightarrow \Hom_c(\Gamma^u_C,A)$ the
{\em crystallographic shadow map for\/ $\Gamma$}. 
\end{definition}
In fact, the  crystallographic shadow map for $\Gamma$ is just the
restriction of the shadow map $s_u: \Hom(\Gamma, A) \rightarrow \Hom(\Gamma^u_C,A)$,
see Proposition \ref{shadowdescent}. 

\paragraph{Adapted systems of generators}
Let $\Gamma= CF \leq \bH_{\Gamma}$ be a torsionfree polycyclic group 
which admits a nilpotent supplement $C$ for $F= \Fitt(\Gamma)$. 
\begin{definition} 
Let $S = \{ \gamma_i \in \Gamma$,  $i=1,  \ldots,  l \}$, be a system of
generators  for $\Gamma$. $S$ is called 
{\em adapted to the supplement $C$}, if $C = \, < \! \gamma_i , \, i=1,  \ldots,  s \!>$, and 
$F= \, < \! \gamma_i , \, i=s+1,  \ldots , l \! > $. 
\end{definition} 
Let  $S$ be an adapted system of generators for
$\Gamma = C \, \Fitt(\Gamma)$. Then it follows that
$$ S_u \, = \, \{ \,  \theta_i = (\gamma_i)_u \in \ur(H_{\Gamma}) \, \}$$
is a set of generators for $\Gamma^u_C$. Assume now also that
$\rank \Gamma = \dim V$. 

\begin{proposition} \label{adapted}
Let $\rho: \Gamma \rightarrow \Aff(V)$
be a homomorphism, $S$ an adapted system of generators. Then $\rho$ is 
a crystallographic homomorphism if and only if the assignment 
$$   \theta_i \, \longmapsto \, \rho(\gamma_i)_u \; \; , \; \theta_i \in S_u$$
defines a crystallographic homomorphism 
$\rho^u \in \Hom_c(\Gamma^u_C, \Aff(V))$. 
\end{proposition}
\begin{proof} By Proposition \ref{shadowdescent}, there exists
a unique homomorphism $$ \rho^u  \, \in \, \Hom(\Gamma_C^u, \Aff(V)) \; , $$
so that  $\rho$ descends to $\rho^u$, and  moreover  
$\rho^u(\theta_i) = \rho(\gamma_i)_u$.   
The proposition  then follows from Proposition \ref{cshadowdescent}.
\end{proof} 
\begin{remark} By the remark following Theorem \ref{acextension}, every
crystallographic group has a unipotent Fitting subgroup.  
Astonishingly enough, we do not have to {\em assume\/} this necessary condition on 
$\rho$ to ensure that $\rho$ is crystallographic. 
In fact, the proof of Proposition \ref{crystshadow} shows that 
if $\rho^u$ is crystallographic then also  
$\rho(\gamma) =\rho(\gamma)_u$,  for all $\gamma \in \Fitt(\Gamma)$.
\end{remark}

\paragraph{The crystallographic restriction}
A nilpotent crystallographic group $\Gamma$ is necessarily unipotent. 
This lead us to consider the space of unipotent homomorphisms 
$\Hom_u(\Gamma, \Aff(V))$. For a general virtually polycyclic \wfn-group $\Gamma$,
a corresponding restriction arises from the internal structure
of the algebraic hull. \newline 

For $g \in H_{\Gamma}$,  we let $c_U(g) \in \Aut(U_{\Gamma})$ denote 
conjugation with $g$ on $U_{\Gamma}=\ur(H_{\Gamma})$. Taking the differential 
of $c_U(g)$ in the identity defines a representation
$$  \alpha_U:  H_{\Gamma}    \longrightarrow   \Aut(\lu_{\Gamma}) $$
on the Lie algebra $\lu_{\Gamma}$ of $U_{\Gamma}$. Since $H_{\Gamma}$
has a strong unipotent radical, the kernel of $\alpha_U$ is contained
in $U$. For $\gamma \in H_{\Gamma}$, we define the characteristic 
polynomial 
$$      \chi_{U}(\gamma,T) \; = \; \det\left(T {\rm id}- \alpha_U(\gamma)\right) \;  , $$
and, correspondingly for $g \in \Aff(V)$, we let $\chi(g , T)$ denote 
the characteristic polynomial of $g$ as an element of $\GL(V \oplus \bbR)$. 
The crystallographic restrictions for $\Gamma$ are described by the following:
\begin{lemma} Let $\Gamma$ be virtually polycyclic \wfn-group. 
If $\rho: \Gamma \rightarrow \Aff(V)$ is a
crystallographic homomorphism then, for all $\gamma \in \Gamma$,
\begin{equation} \label{cr}
      \chi(\, \rho(\gamma) , T \, ) \;  = \;  \chi_U (\gamma, T) \;  (T-1) \; .
\end{equation} 
\end{lemma}
\begin{proof} Since $\rho$ is crystallographic there exists a 
$u$-simply transitive homomorphism $\rho_{\Gamma}: H_{\Gamma} \rightarrow \Aff(V)$
which extends $\rho$, so that $U= \rho_{\Gamma}(\U_{\Gamma})$ is 
a simply transitive subgroup of $\Aff(V)$. $\rho(H_{\Gamma})$ acts by
conjugation on $U$, as well as on the Lie algebra $\lu \subset \End(V \oplus \bbR)$ 
of $U$. Let $\phi_{\Gamma}:  \lu_{\Gamma} \rightarrow \lu$ be the differential
of the homomorphism $\rho_{\Gamma}$. The isomorphism  $\phi_{\Gamma}$ is equivariant
with respect to conjugation, i.e.,  $\phi_{\Gamma}$ satisfies, for all $X \in \lu_{\Gamma}$, 
$g \in H_{\Gamma}$, 
$$           \phi_{\Gamma}( \alpha(g) X )  =    \rho(g)   \phi_{\Gamma}(X) \rho(g)^{-1}
\; . $$
Since $\chi_U(\gamma, T) =  \chi_U(\gamma_s,T)$ we are allowed to 
assume that $\gamma = \gamma_s$ is a semisimple element. We may also
assume  that $\rho(\gamma_s)= \rho(\gamma)_s \in \GL(V)$. Let $o_0: \lu \rightarrow V$
be the differential of the orbit map of $U$ in $0$ (the evaluation map in $0$). 
An immediate calculation shows that the isomorphism  
$o_0 \, \phi_{\Gamma}: \lu_{\Gamma} \rightarrow V$
is equivariant with respect to the action of $\gamma_s$, i.e., for all $X \in \lu_{\Gamma}$,
$$    o_0\,  \phi_{\Gamma} \, ( \alpha(\gamma_s ) X) \, = \,   \rho(\gamma_s) \, o_0 \, \phi_{\Gamma}(X) $$
Hence, $\chi_U(\gamma_s,T) =   \chi_{\GL(V)} (\rho(\gamma)_s ,T)$, and the lemma follows. 
\end{proof} 
We define now a certain Zariski-closed subspace of $\Hom(\Gamma,A)$.
Namely, we put
$$   \Hom_{\chi} (\Gamma, A) \;  = \; \left\{ \, \rho \mid 
\chi(\rho(\gamma),T) = \chi_U(\gamma,T) (T-1)) 
\mbox{ , for all } \gamma \in \Gamma \, \right\}   \; $$
for the subvariety of homomorphisms of type $A$ which satisfy the
crystallographic restrictions.  
If $\Gamma$ is nilpotent, then $\Hom_{\chi} (\Gamma, A) = \Hom_u(\Gamma,A)$.
From the lemma we infer that 
$$ \Hom_{c} (\Gamma, A) \, \subset \, \Hom_{\chi} (\Gamma, A) \; . $$ 
We will show next that $\Hom_{c} (\Gamma,A)$ is 
Zariski-open in $\Hom_{\chi} (\Gamma,A)$.

\paragraph{The variety of crystallographic homomorphisms} 
We are ready now to generalize Theorem \ref{IsZariskiopenN} from nilpotent groups
to virtually polycyclic groups. Also 
we consider the situation now for any Zariski-closed subgroup 
$A \leq \Aff(V)$. If $\Gamma$ is a polycyclic finite index subgroup of $\Delta$,
we let $  \Hom_{\chi} (\Delta, A) \subset \Hom(\Delta, A)$ be the
Zariski-closed subset of homomorphisms which satisfy the
crystallographic restrictions (\ref{cr}) on $\Gamma$.  

\begin{theorem} \label{IsZariskiopen}  
Let $\Delta$ be a virtually polycyclic \wfn-group which satisfies 
$\rank \Delta = \dim V$.
Then there exists a unique $\bbQ$-defined polynomial function $\delta$ on 
the real algebraic variety\/
$\Hom_{\chi}(\Delta, A)$, such that 
$$  
 \Hom_c(\Delta, A) \; = \; \Hom_{\chi}(\Delta, A)_{ \delta} \; . \; \; 
$$ 
Moreover, $\delta$ is a relative invariant
for the conjugation action  of $A$ on $\Hom_{\chi}(\Delta, A)$.
\end{theorem} 
\begin{proof} 
Since $\Hom(\Delta,A)$ is 
Zariski-closed in $\Hom(\Delta, \Aff(V))$, it is
clearly enough to prove the result in the case $A= \Aff(V)$.
Observe that $\rho(\Delta)$ is a crystallographic group if
and only if a finite index subgroup of $\rho(\Delta)$ is a crystallographic
group. Since $\Delta$ is a \wfn-group, $\rho$ is crystallographic (in particular
injective) if and only if it is so on a finite index subgroup $\Gamma$. Since the
restriction map $\Hom(\Delta, \Aff(V)) \rightarrow \Hom(\Gamma, \Aff(V))$
is algebraic, it is therefore enough to show the theorem for $\Gamma$. 
By Proposition \ref{supplexists},  there exists a finite index subgroup $\Gamma$ 
of $\Delta$ which is torsionfree and admits a nilpotent supplement 
$C$ for $\Fitt(\Gamma)$. We now show the result for $\Gamma$.  

By Theorem \ref{IsZariskiopenN}, there exists a polynomial function 
$\delta_u$ on $\Hom_u(\Gamma_C^u, \Aff(V))$
which is a relative invariant and tests 
crystallography for  $\Hom(\Gamma_C^u, \Aff(V))$. 
The shadow map for $\Gamma$ is a map 
$$ s_u: \Hom(\Gamma, \Aff(V)) \longrightarrow
\Hom_u(\Gamma_C^u, \Aff(V)) \; . $$ 
Therefore (compare Proposition \ref{adapted}) the function
$\delta = \delta_u \circ s_u$
tests crystallography on $\Hom(\Gamma, \Aff(V))$, and
in particular on $\Hom_{\chi}(\Gamma, \Aff(V))$.   

We now show that $s_u$ is continuous on 
$\Hom_{\chi}(\Gamma, \Aff(V))$. Let $S= \{ \gamma_i \}$
be an adapted system of generators for $\Gamma$.
By Proposition
\ref{adapted}, the shadow map $s_u$ is obtained
by computing $\rho(\gamma_i)_u$, for each generator
$\theta_i = (\gamma_i)_u$ of the unipotent 
shadow $\Gamma_C^u$. Now,  for $g \in \GL(W)$, 
we may compute the unipotent part of the multiplicative 
Jordan decomposition of $g$ by the 
formula 
$$   g_u   = I + g_s^{-1} g_n \, , $$
where $g = g_s +g _n$ is the additive Jordan-decomposition
in $\End(W)$. Moreover $g_s = P(g)$, $g_n= Q(g)$, where
$P(T),Q(T) \in \bbQ[T]$ are certain polynomials which
depend only on the characteristic polynomial of $g \in \GL(W)$.  
If $\rho \in \Hom_{\chi}(\Gamma, \Aff(V))$ satisfies   
the crystallographic restriction then $\rho(\gamma_i)$
has characteristic polynomial $\chi_U(\gamma,T) \, (T-1)$ 
which does not depend on $\rho$. Therefore, there
exist polynomials $P_i(T)$, $Q_i(T)$ such that, for every 
$\rho \in \Hom_{\chi}(\Gamma, \Aff(V))$, 
$$   \rho(\gamma_i)_u  =    I +  P_i(\rho(\gamma_i) )^{-1} Q_i(\rho(\gamma_i))\; .$$ 
This means that $s_u$ is a $\bbQ$-defined 
polynomial map on $\Hom_{\chi}(\Gamma, \Aff(V))$, 
and so is $\delta = \delta_u \circ s_u$.
It is easy to see that
$s_u$ is $\Aff(V)$-equivariant. Hence, $\delta$ is a 
$\bbQ$-defined polynomial function which is a relative invariant, 
and  the theorem is proved. 
\end{proof} 

\begin{remark} Let $\Delta$ be as in Theorem \ref{IsZariskiopen}. 
The theorem implies that the existence problem
for crystallographic subgroups of $\Aff(V)$ which are isomorphic to $\Delta$ 
may be reduced to the existence problem for real solutions 
of a certain system of algebraic equations with rational coefficients.   
Therefore, by
the Tarski-Seidenberg Theorem this decision problem admits an effective
solution. 
\end{remark} 

In particular, Theorem \ref{IsZariskiopen} shows that 
$\Hom_c(\Delta, A)$ is a Zariski-open subset of\/  
$\Hom_{\chi}(\Delta, A)$. 
The theorem therefore also establishes a certain {\em rigidity property\/} for affine 
crystallographic homomorphisms: 
{\em If  $\rho \in \Hom_c(\Gamma, \Aff(V))$ then any nearby 
$\rho' \in  \Hom_{\chi}(\Gamma, \Aff(V))$ is in $\Hom_c(\Gamma, \Aff(V))$.}
In his celebrated paper \cite{Weil},  A.\ Weil proved the following\footnote{Abels (see \cite{Abels}) gave generalizations 
of Weil's theorem in the general context of proper actions on homogeneous 
spaces of Lie groups.}: \\

Let $\Gamma \leq G$ be
a discrete cocompact subgroup of the Lie group $G$. Then the 
space $$ R(\Gamma,G) = \{\rho: \Gamma \rightarrow G \mid \rho 
\mbox{ is injective and } \rho(\Gamma)
\mbox{ is discrete and cocompact} \} $$
is open in the space of all homomorphisms from $\Gamma$ to $G$.  
(The space $R(\Gamma,G)$ is now called the Weil space.)\\

\noindent 
Our  Theorem \ref{IsZariskiopen} thus may be interpreted as an analogue 
of the result of Weil in 
the specific situation of crystallographic affine actions for 
virtually polycyclic group. \\ 

We state some first consequences now. The main applications, however, 
concern the structure and topology of the deformation
spaces of affine crystallographic actions. This will be
described in the following chapter of this article. 

\begin{corollary} Let $\Gamma$ be a virtually polycyclic group. Then
$\Hom_c(\Gamma,A)$ is a locally closed subset with respect to
the Zariski-topology on $\Hom(\Gamma, A)$. 
\end{corollary}

As a special open subset of an affine real algebraic variety 
the set of crystallographic homomorphisms has a natural algebraic structure itself: 
\begin{corollary} \label{Astructure} 
$\Hom_c(\Gamma, \Aff(V))$ 
has the structure of a $\bbQ$-defined real affine algebraic variety.
In particular, $\Hom_c(\Gamma, \Aff(V))$ has only finitely many connected
components also in the Hausdorff-topology. 
\end{corollary} 

Another important application of the reasoning in
the proof of Theorem \ref{IsZariskiopen} is:
\begin{corollary} \label{shadowmapisc} 
Let $\Theta$ be a unipotent shadow for $\Gamma$. 
Then  the shadowmap 
$$ s_u: \Hom_c(\Gamma, A) \longrightarrow \Hom_c(\Theta,A) $$ 
is a morphism of real algebraic varieties defined over $\bbQ$. 
In particular, $s_u$ is a continuous map. 
\end{corollary}  
\begin{proof} In fact, in the course of the proof of Theorem \ref{IsZariskiopen}
we proved that the shadow map 
$$   s_u:  \Hom_{\chi}(\Gamma,A) \longrightarrow \Hom_u(\Theta,A)  $$
is given by certain rational polynomials with respect to an adapted system of
generators of $\Gamma$. Therefore $s_u$ is algebraic also on 
 $\Hom_c(\Gamma, A) \subset \Hom_{\chi}(\Gamma,A)$.
\end{proof}

\chapter{Deformation Spaces} \label{defspaces}
In this chapter we study properties of the 
deformation spaces of affine crystallographic actions
of virtually polycyclic groups. These spaces are
defined in a purely algebraic way as quotients
of spaces of homomorphisms with an appropiate topology. 
If $M$ is a fixed compact smooth manifold then the deformation space 
$\Def(M,A)$ of complete affine 
$A$-structures on $M$ is the space of such structures up to diffeomorphism,
equipped with the $C^{\infty}$-topology from the space of
developing maps. If $\Gamma$ is the fundamental group of
$M$ the {\em deformation theorem\/} of Thurston links the spaces
$\Def(M,A)$ and $\Def(\Gamma,A)$ via a continuous 
and open map, the {\em holonomy map\/}
$$ hol: \Def(M,A) \longrightarrow \Def(\Gamma,A) \;  . $$
Note that the theory of Thurston covers the general situation 
of locally homogeneous manifolds $M$,  modelled on a homogeneous
space $G/H$, see \cite{ThurstonL}, 
and, in particular, \cite{Epstein,Goldman} for further explanation
and proofs. In \cite{BauesI} implications are discussed
for the deformations of compact complete afffine manifolds. 

In this realm, Kobayashi \cite{Kob5},
more generally formulates the program to determine the deformation
spaces for \emph{proper}, not necessarily cocompact, actions of a 
group $\Gamma$ on a homogeneous
space $G/H$. The situation $G= \Aff(V)$, and $G/H = V$ is
an important special case in this program.\footnote{See \cite{Bak_Kh, Bak_Ke_Yo,KobNas,Yoshino} for recent contributions on this problem.} 

If $\Gamma$ is a torsionfree virtually polycyclic ACG of type $A$ then the 
quotient space $M= \Gamma \lmod V$ is a compact complete affine
manifold which admits an affine atlas with coordinate changes in the 
group $A$. It is also a smooth
aspherical compact manifold with fundamental group $\Gamma$. 
Results of the previous chapter (cf. Proposition\ref{syndeticH})
imply that  the smooth manifold $M$ falls into 
the class of infrasolvmanifolds. 
In this class of smooth manifolds
the diffeomorphism type of $M$ is determined by the 
fundamental group $\Gamma$ alone. This is proved in
\cite{BauesI}.\footnote{A result in \cite{FriedGoldman}
shows that any two compact complete affine manifolds 
with virtually solvable fundamental group are (polynomially) 
diffeomorphic.} 
This \emph{smooth rigidity} for compact complete affine 
manifolds implies that  the holononomy map is a homeomorphism,
compare \cite{BauesV}. Therefore,  the algebraic viewpoint
on the deformation space $\Def(M,A)$ captures the whole
picture.  In this sense, we allow ourselves at some points
to speak about the deformation spaces of affine
structures on certain manifolds, although we only
give results on the affine crystallographic 
actions of their fundamental groups.

\section{The geometry of the deformation space}
Let $\Gamma$ be a group, and let $\Hom_c(\Gamma,A)$ be
the space of crystallographic homomorphisms of type $A$. The group $A$
acts by conjugation on the space of homomorphisms. Our principal object
of interest is the deformation space of crystallographic homomorphisms
of type $A$ which is defined as
$$    \Def(\Gamma,A) \;  =  \; 
\Hom_c(\Gamma,A) \,  \rmod \! {_{\displaystyle A}} \; . $$
In this situation, the group $\Aut(\Gamma)$ of automorphisms of $\Gamma$ acts
freely on  $\Hom_c(\Gamma,A)$ commuting with the action of $A$. 
The quotient topological space 
$$    \Cha(\Gamma,A) \; =  \;  
{_{\displaystyle \Aut(\Gamma)}}  \lmod \,   \Hom_c(\Gamma,A) $$
is the space of crystallographic subgroups of $A$ which are
isomorphic to $\Gamma$. This space is sometimes called the 
{\em Chabauty space}. The group $A$ acts by conjugation on 
$\Cha(\Gamma,A)$, and the action of
the group $\Aut(\Gamma)$ on $\Hom_c(\Gamma,A)$ factorizes to an
action of the outer automorphism group 
$\Out(\Gamma) =\Aut(\Gamma) / {\rm Inn}(\Gamma)$ 
on $\Def(\Gamma,A)$. The quotient space 
$$    \Mod(\Gamma,A)  \; = \; 
{_{\displaystyle \Out(\Gamma)}}  \lmod \,    \Def(\Gamma,A) $$
is called the {\em moduli space}. The situation is described in the following 
commutative diagram of maps and spaces:
\begin{equation}
\xymatrix{    & \Hom_c(\Gamma,A)  \ar@/_/[dl] \ar@/^/[ddr]  &  &    \\
  \Cha(\Gamma,A) \ar@/_/[rdd] &  &  & \\
   & & \Def(\Gamma,A) \ar@/^/[ld]  & \\
   & \Mod(\Gamma,A) & & 
} \; . 
\end{equation} 
Basic problems in the deformation theory for $\Gamma$ are to understand the 
topology and geometry of the spaces involved in this diagram.\newline

We restrict ourselves here to the case that $\Gamma$ is a 
virtually polycyclic group. Recall from chapter \ref{choms}
that then the affine crystallographic Weil space $\Hom_c(\Gamma,A)$
carries a Hausdorff-topology, and also a Zariski-topology which is 
induced from $\Hom(\Gamma,A)$. This is the essence of our 
Theorem \ref{IsZariskiopen}. The next step is to learn more
about the remaining spaces. Then the following interesting 
picture emerges:

\paragraph{Algebraic and arithmetic nature of deformation spaces}
The deformation space arises 
as a quotient of the real points of the $\bbQ$-defined algebraic 
variety $\Hom_c(\Gamma,A)$ by a compatible $\bbQ$-defined algebraic action of
$A$.  The moduli space $\Mod(\Gamma,A)$ 
arises as a quotient of the deformation space by a compatible
$\bbQ$-defined algebraic action of the \emph{arithmetic group}\footnote{see \cite{Serre}
for a discussion of the notion of arithmetic group}  $\Out(\Gamma)$. 
In fact, as is proved recently in \cite{BauesGrunewald}, the outer automorphism
group of a virtually polycyclic group is an arithmetic group. 

Moreover, the action of $\Aut(\Gamma)$ on $\Hom_{c}(\Gamma,A)$ 
extends to an algebraic $\bbQ$-defined action of the \emph{algebraic automorphism group} 
$\Aut_{a}(H_{\Gamma})$, and
this action factors over the algebraic inner group $$\Out_{a}(H_{\Gamma}) = 
\Aut_{a}(H_{\Gamma})/ {\rm Inn} (H_{\Gamma})$$ 
to an action on $\Def(\Gamma,A)$, which extends the $\Out(\Gamma)$-action
on $\Def(\Gamma,A)$.  These facts are announced in \cite{BauesU}. In chapter 3 of 
this article we will give a proof for (virtually) \ftn-groups $\Gamma$.
(The proof for the general fact needs more of the machinery 
which relates the groups $\Aut(\Gamma)$ and $\Aut_{a}(H_{\Gamma})$,
as is developed in \cite{BauesI, BauesGrunewald}.)

\paragraph{Topology of deformation spaces}

The deformation space, as well as the
moduli space have nice topological properties in some well known
geometric situations. For example they are Hausdorff-manifolds in
some cases. This is in particular true for the deformation spaces of
constant curvature Riemannian metrics, and the related Teichm\"uller theory
of complex analytic structures on surfaces, see for example \cite{Ratcliffe}.
Such properties were suspected to be no longer true 
for the deformation spaces of affine crystallographic groups, see \cite{Goldman}. 

To the contrary it was
shown in \cite{BauesG} that the deformation space $\Def(\bbZ^2,\Aff(\bbR^2))$
is a Hausdorff-space, in fact, homeomorphic to $\bbR^2$. Whereas the moduli space 
$$\Mod(\bbZ^2,\bbR^2) =  {_{\displaystyle \GL_{2}(\bbZ)}}  \lmod \,  \bbR^2$$ 
is highly non Hausdorff. 

Note that via the holonomy map the space $\Def(\bbZ^2,\Aff(\bbR^2))$ is homeomorphic to the  deformation space of complete affine structures on 
the two torus $T^2$. A geometric interpretation of the coordinates of 
$\Def(\bbZ^2,\Aff(\bbR^2))$
in terms of periods of developing maps for affine structures on $T^2$
is given in \cite{BauesGoldman}. 

This example also nicely illustrates the
\emph{arithmetic and algebraic nature} of the deformation space.
For an interpretation of rational points in $\bbR^2$ in terms of 
geometric properties of the corresponding affine structures on $T^2$, 
see \cite{BauesG}. It is shown there that rationality of an affine structure 
is equivalent to the existence of closed geodesics.\\

In this chapter we concentrate our 
effort to a further study of the deformation spaces $\Def(\Gamma,A)$. 
We will show  that in some cases the deformation spaces
for affine crystallographic groups are in fact Hausdorff and
homeomorphic to a semi-algebraic set.
Let us mention first the
following basic result on the separation properties 
of deformation spaces:

\begin{theorem} Let $\Gamma$ be a 
virtually polycyclic group, and $A \leq \Aff(V)$ a Zariski-closed
subgroup. Then the deformation space $\Def(\Gamma,A)$ is a $T1$-topological 
space\footnote{meaning that every point is closed} with 
finitely many connected components.
\end{theorem}
\begin{proof}
By Corollary \ref{Astructure}, the space $\Hom_{c}(\Gamma, A)$ has
only finitely many components.  Therefore, the same holds for
the quotient space $\Def(\Gamma,A)$. 

To see that points in $\Def(\Gamma,A)$ are closed, we shall 
need some of the machinery, which will be developed in the sequel.
We can argue as follows: By the realisation theorem, Theorem 
\ref{shadowembedding}, it is enough to prove separation properties
for deformation spaces of \ftn-groups. If $\Gamma$ is a \ftn-group then, by 
Proposition \ref{Matsushima},  the deformation space $\Def(\Gamma,A)$
may be represented as the quotient space of a real algebraic
variety by a unipotent action of the real Malcev hull $U_{\Gamma}$.
By the Rosenlicht theorem \cite{Rosenlicht} every orbit of $U_{\Gamma}$ is
closed. Therefore, the points in $\Def(\Gamma,A)$ are closed.
\end{proof}

\paragraph{Geometric notions for deformation spaces}
Geometric notions are introduced by the natural group actions 
on our spaces. Let $\rho \in \Hom_c(\Gamma, A)$ be a crystallographic
homomorphism. An important concept concerns the dimension of the
deformation space. 

\begin{definition}
The action $\rho$ is called {\em (locally) $A$-rigid\/} 
if the orbit $A \rho$ is an open set in $\Hom_c(\Gamma, A)$. 
\end{definition}
Local rigidity of $\rho$ is equivalent to the fact that 
the point $[\rho \,]$ constitutes a connected
component of the deformation space $\Def(\Gamma,A)$. 
The action $\rho$ is {\em rigid\/} if the deformation
space is a point. It is well known that rigidity occurs
for certain geometric structures, for example
hyperbolic structures in dimension $n>2$. 
We expect that local rigidity fails for most 
affine crystallographic actions.
\newline  

Another important notion is
\begin{definition} \label{convex}
The deformation space $\Def(\Gamma,A)$ 
is called {\em convex\/} if every finite subgroup of $\Out(\Gamma)$
has a fixed point. A special case occurs if the group $\Out(\Gamma)$ has 
fixed points. Then we will call $\Def(\Gamma,A)$ {\em fixed pointed}.    
\end{definition} 

The convexity of the deformation space of $\Gamma$ is in particular 
important because of the role it plays in 
the realization problem for finite extensions, 
see chapter \ref{realizations}.
In Theorem  \ref{convexftngroups} below we show 
that natural classes of \ftn-groups, and as a consequence many  
\ftn-by finite and torsionfree polycyclic groups
have fixed-pointed or convex deformation spaces.  


\section{Models for deformation spaces} 
The purpose of this section is to provide 
some examples of deformation spaces for 
affine crystallographic groups which are 
Hausdorff topological spaces. We first give
some details on the structure of deformation   
spaces for tori. The results
of chapter \ref{realizations} allow
then to construct more examples of
deformation spaces which are Hausdorff.
In particular we show that the deformation spaces
of complete affine structures of certain
three-manifolds are Hausdorff spaces.

\subsection{Deformation spaces for affine tori}  \label{atori}
In \cite{BauesV} it was proved that the deformation
spaces $\Def(\bbZ^n, \Aff(V))$ of affine tori are 
homeomorphic to a semi-algebraic set in some Euclidean space. 
Here we want to study 
more closely these deformation spaces and 
the structure 
of the algebraic varieties associated to them. \newline

Theorem \ref{IsZariskiopen} asserts that 
the space $ \Hom_c(\bbZ^n, A)$ has a natural 
structure as a real algebraic variety. 
Here we want to describe the 
algebraic variety $ \Hom_c(\bbZ^n, A)$, that is, 
the space of crystallographic homomorphisms of type $A$, 
and the associated deformation space $\Def(\bbZ^n,A)$ in
more detail. A basic result on  these deformation spaces is:
\begin{theorem} \label{tsemialg}
The deformation space $\Def(\bbZ^n, A)$ is 
homeomorphic to a semi algebraic
set, and in particular $\Def(\bbZ^n, A)$ is a Hausdorff space.
\end{theorem}
\begin{proof} The proof of \cite[Corollary 2.9]{BauesV} for $A = \Aff(V)$
carries over almost verbatim to a more general $A$ in the case 
that $\GL_A= A \cap \GL(V)$ is a real reductive group. 
Let
us briefly recall the argument. Since $\bbZ^n$ is abelian
the centralizer of every cystallographic action is a
simply transitive subgroup of $A$, and 
$\Def(\bbZ^n, A)$ is homeomorphic to the orbit space 
$\Hom_c(\bbZ^n, A)/ \GL_A$. 
It follows from Lemma \ref{centralizer} that $\GL_A$ acts 
freely on $\Hom_c(\bbZ^n, A)$, and in particular all orbits
are closed.
It is known that the  space of closed orbits of the reductive 
group $\GL_A$ is homeomorphic to a semi-algebraic set. 
For a more general $A$ the theorem will follow from
Proposition \ref{fproduct} below.  
\end{proof} 

Using a fixed choice of basis we embed $\bbZ^n$ as
a lattice in $V$. As a real algebraic group, $V$ identifies
then with the real Malcev hull of $\bbZ^n$. As shown in section \ref{Groupactions}
the automorphism group of $V$, that is the group $\GL(V)$,  
acts freely on the space $\Hom_c(\bbZ^n, A)$. For $g \in \GL(V)$ 
the action is given by
$$  \rho \longmapsto \rho \, g = \bar{\rho} \, g^{-1} \, j \; ,  $$
where $\bar{\rho}: V \rightarrow A$ is the unique simply
transitive representation of $V$ which extends $\rho$,
and $j: \bbZ^n \rightarrow V$ is the fixed embedding of $\bbZ^n$ in $V$.
This action of $\GL(V)$ commutes with the conjugation action of $A$. 
The quotient space 
$$G_{st}(V,A) \;  = \; \;     {_{\displaystyle \GL(V)}}  \lmod \,  \Hom_c(\bbZ^n, A)$$
is the space of simply transitive abelian subgroups of $A$.
We describe the 
structure of the algebraic variety $\Hom_c(\bbZ^n, A)$
as follows: 

\begin{proposition} \label{structure} 
The space $G_{st}(V,A)$ is homeomorphic to
a real algebraic variety and there is an isomorphism of algebraic 
varieties $$ \Hom_c(\Gamma,A) = \GL(V)\times  G_{st}(V,A) \; ,  $$
such that the action of $\GL(V)$ on  $\Hom_c(\Gamma,A)$ is
given by right multiplication on the first factor. The action
of $\GL_A$ by conjugation on $ \Hom_c(\Gamma,A)$ is given
by its action on  $G_{st}(V,A)$ and by left multiplication on
$\GL(V)$.  
\end{proposition}
\begin{proof} 
Each $\rho \in \Hom_c(\bbZ^n,A)$ extends to 
a homomorphism $U_{\bbZ^n}= V  \rightarrow A$.
The translational components of its derivative define
a linear isomorphism $\bar t(\rho) : V \rightarrow V$, since $V$
acts simply transitively.  
Also it is easy to see that 
the $\GL(V)$-  (change of basis) action 
on $\Hom_c(\bbZ^n,A)$ corresponds to (transposed) right
multiplication on $\GL(V)$, i.e., 
$\bar t (\rho \, g) = \bar t(\rho) (g^{-1})^t$.
In particular, each element
of $G_{st}(V,A)$ has a unique representative 
$\rho_{id} \in \Hom_c(\bbZ^n,A)$
which satisfies $\bar t(\rho_{id}) = id_V$, and this 
representative is computed by the formula
$\rho_{id} = \rho \,  \bar t(\rho)^{-1}$.  
Therefore $G_{st}(V,A)$ identifies with a Zariski-closed
subspace $\Hom_{id}(\bbZ^n,A)$ of $\Hom_c(\bbZ^n,A)$,
and the natural map 
$$ \GL(V) \times \Hom_{id}(\bbZ^n,A) \longrightarrow  \Hom_c(\bbZ^n,A)$$ 
given by $(g,\rho) \mapsto \rho \, g$ is an algebraic isomorphism. 
The equivariance statements are easy to verify. Hence 
the proposition follows.
\end{proof}

We let $\GL_A$ act on $\GL(V)$ by left multiplication and put
$$  X_A \; = \;    {_{\displaystyle \GL_A}}\! \lmod \, \GL(V)    \; .$$
The following result reveals the structure of the 
deformation spaces $\Def(\bbZ^n, A)$ and the nature of the 
$\GL(V)$-action on  $\Def(\bbZ^n, A)$ more precisely: 
\begin{proposition} \label{fproduct} 
Let $F= {G}_{st}(V,A)$ be the algebraic 
variety of simply transitive abelian subgroups of $A$ with the natural
action of $\GL(A)$. 
Then 
$$  \Def(\bbZ^n, A) \; = \; \GL(V) \underset{\GL_A}{\times}  F$$
is homeomorphic to a fiber product, 
and in particular a bundle over the homogeneous
space $X_A$ with fiber $F$. The $\GL(V)$-action on $\Def(\bbZ^n, A)$
corresponds to the natural $\GL(V)$-action on the fiber product which is
induced by right multiplication of  $\GL(V)$ on itself. 
\end{proposition} 
\begin{proof} We already remarked in the proof of Theorem \ref{tsemialg}
that $\Def(\bbZ^n, A)$ is homeomorphic to the quotient by $\GL_A$.
Therefore the proposition follows from Proposition \ref{structure}. 
\end{proof}

As a consequence we see that $\Def(\bbZ^n,\Aff(V))$ is actually 
homeomorphic to a real algebraic variety. 
\begin{corollary} The deformation space $\Def(\bbZ^n,\Aff(V))$ is
homeomorphic to the real algebraic variety ${G}_{st}(V,\Aff(V))$.  
\end{corollary}
\begin{proof}  The group $\GL_{\Aff(V)}(V)= \GL(V)$ acts by conjugation
on $\Hom_c(\bbZ^n,\Aff(V))$, so that $t(\rho^g) = g t(\rho)$. Therefore
each orbit intersects 
$\Hom_{id}(\bbZ^n,\Aff(V))$ 
precisely once. The corollary now follows from Proposition \ref{fproduct}.  
\end{proof}

\begin{example} The deformation space $\Def(\bbZ^2, \Aff(\bbR^2))$
is homeomorphic to $\bbR^2$ and the actions of $\GL(\bbR^2)$ and 
$\GL_2(\bbZ)$ on $\Def(\bbZ^2, \Aff(\bbR^2))$ correspond to
the canonical linear actions. The fixed point 
$o \in \Def(\bbZ^2, \Aff(\bbR^2))$ corresponds to
the natural action of $\bbZ^2$ by translations. 
(See \cite{BauesG})
\end{example}

We next cover a few more special cases:
\begin{example}
Let $q$ be an inner product on $V$, $\O(q)$ the
group of linear $q$-isometries, and $\A(q)$ the group
of affine isometries for $q$. It is easy to see and well known
that ${G}_{st}(V,A(q))$ consists of the groups of
translations only. (A proof may be found in \cite{BauesCortes})
It follows that  $ \Def(\bbZ^n, \A(q))$ is homeomorphic to
$X_{A(q)}$.
\end{example}

The symplectic case is already more interesting. 
Let $V$ be a vector space with a nondegenerate
alternating product $\omega$. We let  $\SP(\omega)$ denote the
group of linear $\omega$-isometries of $V$, 
and $\A(\omega)$ the group of affine isometries for $\omega$.
Let ${\cal U}_k$ denote the tautological 
vector bundle over the Grassmanian 
of $k$-dimensional isotropic subspaces of the 
symplectic vector space $(V, \omega)$. It is a
natural $\SP(\omega)$-variety.
Let  $2n = \dim V$. 
The following is proved in \cite{BauesCortes}: 
\begin{proposition} The variety ${G}_{st}(V,A(\omega))$
admits a $\SP(\omega)$-invariant stratification 
$$     {G}_{st}(V,A(\omega)) = \bigcup_{k=0}^n  {\cal G}_{st}(V,A(\omega))_k \;. $$
Each stratum is an open subbundle of the third symmetric power $S^3 {\cal U}_k$ of
${\cal U}_k$ . 
\end{proposition} 
Together with Proposition \ref{fproduct} we obtain a model  for 
the space $\Def(\bbZ^n, A(\omega))$. \newline

\begin{remark}  
As is observed in \cite{BauesCortes} the
study of simply transitive abelian symplectic affine groups plays a role
in the construction of flat models for
{\em special K\"ahler geometry}, a particular
geometry which arises in supersymmetric quantum
field theory.  
\end{remark}

\subsection{Semi-algebraic deformation spaces}
We proved up to now that the spaces $\Def(\bbZ^n, A)$ 
are homeomorphic to semi algebraic sets. The results
of chapter \ref{realizations} allow us to 
further generalize this result. A finite effective extension
group of $\bbZ^n$ is traditionally called a {\em Bieberbach
group}.

\begin{corollary} \label{Bieberbach} 
Let $\Delta$ be a Bieberbach group. 
Then the deformation space $\Def(\Delta, A)$ is homeomorphic
to a semi algebraic set and is a Hausdorff space. 
\end{corollary}  
    
This implies corresponding results on the deformation
spaces of affine space forms which are finitely covered
by a torus. More details and the 
proofs will be given in section \ref{applications1}.
A further result is:
\begin{corollary} \label{abshadow} 
Let $\Gamma$ be a virtually torsionfree polycyclic group such that
the unipotent shadow of \/ $\Gamma$ is abelian.
Then the deformation space $\Def(\Gamma, A)$ is homeomorphic
to a semi algebraic set and is a Hausdorff space. 
\end{corollary}  

The proof will be given in section \ref{srealizations}. 

\paragraph{Deformation spaces of complete affine three manifolds}
The following question was raised by Bill Goldman\footnote{private communication}:
 \newline

{\em Let $M$ be a closed 3-manifold which is a torus bundle over $S^1$.
Is the deformation space of complete affine structures on $M$ Hausdorff?}
\newline

As is proved in \cite{FriedGoldman} if $M$ is a closed complete affine 3-manifold 
then $M$ is finitely covered by a torus bundle over $S^1$.
The previous results imply a partial answer to
the above question:
\begin{theorem} Let $M$ be a closed complete 
affine three manifold, and $A \leq \Aff(3)$ 
a Zariski-closed subgroup. 
\begin{itemize}
\item[i)] If $M$ is finitely covered by a torus, then the deformation
space $\Def(M, A)$ is Hausdorff.
\item[ii)] If $M$ is a torus bundle over $S^1$ where the attaching map
is hyperbolic with a positive trace then $\Def(M, A)$ is
Hausdorff.
\end{itemize}
\end{theorem}
\begin{proof} 
The first claim follows from Corollary \ref{Bieberbach}
above, keeping in mind the nontrivial fact that the deformation space of complete
affine structures on $M$
is homeomorphic to the deformation space of crystallographic
actions of the the fundamental group.

For the second claim, we remark that if the condition on
the bundle is satisfied, the fundamental group $\pi_1(M)$ of $M$ is
a Zariski-dense lattice in the 3-dimensional Lie group ${\rm Sol}$.  
Therefore, the unipotent shadow of $\Gamma= \pi_1(M)$ is abelian,
and Corollary \ref{abshadow} applies. 
\end{proof}

\begin{remark} It would be further interesting to understand the separation properties of
$\ftn$-groups which are lattices in the 3-dimensional Heisenberg group. 
\end{remark}


\section{Convexity properties of deformation spaces} \label{convexamples}
The purpose of this section is to provide 
some examples of deformation spaces for 
affine crystallographic groups which are 
convex in the sense of Definition \ref{convex}.

\paragraph{Strong convexity}
Let $\Gamma$ be a \ftn-group   
and let $U_{\Gamma}$ be a real Malcev hull for $\Gamma$.
The automorphism group of $U_{\Gamma}$ naturally 
acts on the space $\Hom_c(\Gamma,A)$ 
and on the deformation space $\Def(\Gamma,A)$.
Moreover, there is an induced  action of 
$\Out(U_{\Gamma}) = \Aut(U_{\Gamma)} / {\rm Inn}(U_{\Gamma})$ on $\Def(\Gamma,A)$
which extends the natural $\Out(\Gamma)$-action. 
This is explained in section \ref{Groupactions}.\\

We may as well interpret the fixed point properties of the  $\Aut(U_{\Gamma})$-action 
as a {\em convexity property\/} of the deformation space $\Def(\Gamma,A)$.
For a \ftn-group $\Gamma$, we will henceforth 
understand Definition \ref{convex} in the following stronger sense: 

\begin{definition} \label{convex2} 
Let $\Gamma$ be a \ftn-group.
The space $\Def(\Gamma,A)$ is called 
{\em strongly convex\/} if every reductive subgroup of $\Aut(U_{\Gamma})$
has a fixed point in  $\Def(\Gamma,A)$. 
If the group $\Aut(U_{\Gamma})$ has a fixed
point, $\Def(\Gamma,A)$
will be called {\em fixed pointed}. 
\end{definition}

Note, if $\Def(\Gamma,A)$ is fixed pointed in the sense
of Definition \ref{convex2} then every finite subgroup
of $\Out(\Gamma)$ has a fixed point on $\Def(\Gamma,A)$.
In fact, every finite subgroup $\mu \leq \Out(\Gamma)$ may
be lifted to a finite subgroup of $\Aut(U_{\Gamma})$.\footnote{Using Theorem \ref{Mostow}}

\begin{example} Let $\Gamma = \bbZ^n \leq V$ be a lattice,
so that $U_{\Gamma} = V$. The space $\Def(\bbZ^n,\Aff(V))$ has a natural fixed 
point for all automorphisms. The fixed point for $\Out(U_{\bbZ^n}) = \GL(V)$
is given by the natural action of $\bbZ^n$ on $V$ via translations.
(Compare also Proposition \ref{fproduct}) 
\end{example} 

For the moment being, we restrict our interest to the principal 
case $A=\Aff(V)$.
Our aim here is to construct examples of groups with 
fixed pointed deformation spaces $\Def(\Gamma, \Aff(V))$, 
and more generally with convex deformation spaces. 
The first step is to consider \ftn-groups.

\paragraph{Conditions for strong convexity}
Let $\Gamma$ be an \ftn-group. Recall
that the {\em nilpotency class\/} of $\Gamma$ is
the length of a shortest central series for $\Gamma$.
(For example, if $\Gamma$ is abelian then $\Gamma$ is of class 1.) 
By abuse of language we 
say that {\em $\Gamma$ admits an invariant grading\/} 
if the Lie algebra of the Malcev hull $U_{\Gamma}$ has an invariant grading. 
(See Definition \ref{filtrations}.) Let us consider now the
following conditions for $\Gamma$:
\begin{itemize}
\item[i)]  $\Gamma$ is of nilpotency class $\leq 2$,       
\item[ii)] $\Gamma$ is of nilpotency class $\leq 3$,  
\item[iii)] $\Gamma$ admits a positive invariant grading.  
\end{itemize}

It is known from constructions given by Scheuneman that if one of the conditions
i)-iii) is satisfied $\Gamma$ admits affine crystallographic
actions. This will be explained further below.\\

We will prove below:
\begin{theorem} \label{convexftngroups}
Let $\Gamma$ be a \ftn-group. If $\Gamma$ satisfies
condition i) then the deformation space 
$\Def(\Gamma, \Aff(V))$ is fixed pointed. 
If one of the conditions ii) or iii) 
is satisfied by $\Gamma$ then
$\Def(\Gamma, \Aff(V))$ is strongly 
convex. 
\end{theorem}

The proof of the theorem will show that is possible 
to construct affine crystallographic actions which
are fixed by reductive groups if the algebraic 
structure of the \ftn-group $\Gamma$ is not too complicated.
In particular if the rank of $\Gamma$ is small, 
for example, $\rank \Gamma \leq 5$,
conditions ii) or iii) are satisfied. 
This follows from known classification results
for nilpotent Lie algebras.   
On the other hand, it is known that for a ``sufficiently generic''
nilpotent group connected reductive groups 
of automorphisms do not exist at
all, compare \cite{Goze}. 
This leaves us with the open problem: \newline 

{\em Does there exist an affine crystallographic \ftn-group $\Gamma$ with
a non-convex deformation space?}\newline 

The answer to this question is in particular important
in the light of the solution to the realization problems
in chapter \ref{realizations}. For example, by Theorem 
\ref{Affinerealization1} the convexity 
of $\Def(\Gamma,A)$ implies that $\Gamma$ 
satisfies the realization property for finite extensions. 
In particular, every
finite effective extension $\Delta$ of $\Gamma$ admits
an affine crystallographic action. Therefore the question
above is more or less equivalent to the
question: \newline 

{\em Does there exist an affine crystallographic \ftn-group 
$\Gamma$ with a finite effective extension $\Delta$ which is
not isomorphic to an affine crystallographic group?}
\newline 

By Theorem \ref{inheritance1}
the convexity of the deformation space is inherited 
to finite extensions of $\Gamma$. 
Therefore Theorem \ref{convexftngroups} also yields: 

\begin{corollary} \label{convexftnbfgroups}
Let $\Delta$ be a finite extension of an \ftn-group, $\Gamma =\Fitt(\Delta)$.
Then the deformation space $\Def(\Delta, \Aff(V))$ 
is fixed pointed if 
condition i) is satisfied by $\Gamma$.
The deformation space 
$\Def(\Delta, \Aff(V))$ is convex if the group 
$\Gamma$ satisfies one of ii) or iii). 
\end{corollary}

We remark that the corollary in particular establishes the existence 
of affine crystallographic actions for $\Delta$. (See 
also section \ref{existence})
By Theorem \ref{inheritance2}
the convexity of the deformation space is also inherited 
from unipotent shadows. This yields: 

\begin{corollary} \label{convexpcgroups}
Let $\Gamma$ be a torsionfree polycyclic group, $\Theta$ a
unipotent shadow for $\Gamma$. Then the deformation space 
$\Def(\Gamma, \Aff(V))$ is fixed pointed if
condition i) is satisfied by $\Theta$. 
$\Def(\Gamma, \Aff(V))$ is convex if the shadow  
$\Theta$ satisfies one of ii) or iii).
\end{corollary}

The previous corollaries are an immediate consequence
of Theorem \ref{convexftngroups} and the above mentioned
inheritance results of chapter \ref{realizations}.
Theorem \ref{convexftngroups} follows from the
existence of certain well understood simply
transitive unipotent actions constructed
by Scheuneman, see below. \newline  

\begin{prf}{Proof of Theorem \ref{convexftngroups}}
For a simply transitive subgroup $U \leq \Aff(V)$
we let  $\Aut_{\Aff(V)}(U) \leq \Aut(U)$ be the
image of the normalizer of $U$. (cf.\ section \ref{Groupactions}).

Let $\Gamma$ be of nilpotency class $\leq 2$, and $U_{\Gamma}$
its real Malcev hull. By Proposition
\ref{2step}, $U_{\Gamma}$ embeds into $\Aff(V)$ so
that the image $U \leq \Aff(V)$ is a simply transitive
subgroup with the property that $\Aut_{\Aff(V)}(U) = \Aut(U)$.
This representation corresponds to an element $\rho \in \Hom_{c}(\Gamma, \Aff(V))$.

Recall now that $\Aut(U_{\Gamma})$ acts on the
deformation space $\Def(\Gamma, \Aff(V))$, extending the
action of $\Aut(\Gamma)$. 
By Proposition \ref{Autstabilizer}, the subgroup 
$\Aut_{\Aff(V)}(U) \leq \Aut(U)$ 
corresponds to a subgroup of the stabiliser $\Aut(U_{\Gamma})_{[\rho]}$ of
$[\rho] \in \Def(\Gamma, \Aff(V))$.
Therefore, $\Aut_{\Aff(V)}(U) = \Aut(U)$ implies
that $\Aut(U_{\Gamma})_{[\rho]} = \Aut(U_{\Gamma})$.
This means that $[\rho]$ is a fixed point for $\Aut(U_{\Gamma})$.

The second statement is proved similarly. If $\Gamma$ satisfies
condition ii) or condition iii), Scheuneman's examples
(cf.\ Theorem \ref{3step}, Theorem \ref{graded}) provide 
us with crystallographic homomorphisms $\rho$ for $\Gamma$ 
so that $\Aut_{\Aff(V)}(U)$ contains a Levi subgroup of $\Aut(U)$.
This shows that $\Aut(U_{\Gamma})_{[\rho]}$ contains a
Levi subgroup $L$ of $\Aut(U_{\Gamma})$. Therefore $[\rho]$
is a fixed point for $L$. Hence any reductive subgroup
of $\Aut(U_{\Gamma})$ has a fixed point in the
$\Aut(U_{\Gamma})$-orbit of $[\rho \,]$. In particular,
this holds for finite subgroups of $\Aut(U_{\Gamma})$.
\end{prf}

\subsection{Scheuneman's examples}
Scheuneman (cf.\ \cite{Scheune1, 
Scheune2}) constructed 
unipotent simply transitive actions on affine space
to give new examples of compact complete affine manifolds.
We analyze here the automorphism groups of these actions.
Some of the examples we study were
known already to Elie Cartan (cf.\ \cite{Helgason}). 
Scheuneman's examples provide the only known
general method to establish the existence of affine crystallographic
actions on large classes of (reasonably well behaved) torsionfree
nilpotent groups.\newline 
 
To present the constructions,  we need to establish first the
infinitesimal picture of simply transitive actions.

\subsubsection{Simply transitive affine actions of Lie algebras} 
Let ${\lie g}$ be a Lie algebra, $\phi: \lg \rightarrow \lgl(V)$
a representation of $\lg$ on a vector space $V$. A map
$D: \lg \rightarrow V$ is called a derivation for $\phi$ if
$$ D([X,Y]) = \phi_X DY - \phi_Y D X \; \; , \;  \mbox{ for all $X,Y \in \lg$.} $$
$D$ will be called {\em nonsingular\/} if $D$ is an isomorphism
of vector spaces. We call the  pair $(\phi,D)$ an {\em affine
representation on $V$}.
Let  $\la(V) = \lgl(V) \oplus V$ be the Lie algebra of $\Aff(V)$.  
Every Lie algebra homomorphism $\varphi: \lg \rightarrow \la(V)$ is,
with respect to the splitting of $\la(V)$, a sum $\varphi = \phi + D$, such that
$(\phi,D)$ is an affine representation. \\

 Let $U$ denote a simply connected nilpotent Lie group with Lie algebra
 $\lu$. The next proposition is well known (compare \cite{Kim}):

\begin{proposition} 
An affine representation $(\phi,D): \lu \rightarrow \la(V)$ 
is the differential of a unipotent simply transitive representation
of\/ $U$ if and only if $\phi$ is a representation of nilpotent linear
operators and $D$ is a nonsingular derivation for $\phi$.      
\end{proposition} 

If $U$ is a unipotent simply transitive subgroup of $\Aff(V)$,  
we call the Lie algebra $\lu \subset \la(V)$ 
a {\em simply transitive subalgebra\/} of $\la(V)$.

\paragraph{Normalisers of simply transitive actions}
Let $U$ be a unipotent (real) linear algebraic group with
Lie algebra $\lu$. Recall that the automorphism group $\Aut(U)$ 
of $U$ identifies with the group $\Aut(\lu)$ of Lie algebra
automorphisms of $\lu$ by taking
the differential in the identity. (This identification also defines
the natural algebraic group structure on $\Aut(U)$.)

Let $U \leq A$ be a subgroup of the affine group.
We define (see also section \ref{Groupactions}) 
the group $$ \Aut_A(U) \leq \Aut(U) $$ of
$A$-automorphisms of $U$ as the natural image of 
$N_{\GL_A}(U)$ in $\Aut(U)$, and $$ \Aut_A(\lu) \leq \Aut(\lu)$$
as the corresponding group of Lie algebra automorphisms of $\lu$.

Note that $\Aut_A(\lu)$ may also be computed as the image of $N_{\GL_A}(\lu)$ in 
$\Aut(\lu)$, since $\lu$ is a subalgebra of $\la$, where $\GL_A$ acts
by conjugation on  $\la$. 
Now let $\varphi: \lu \rightarrow \la$  be a simply transitive
representation of $\lu$, $\varphi = (\phi, D)$. We put 
$$ \Aut_A(\lu,\varphi) = \Aut_A(\varphi(\lu))^{\varphi^{-1}} \; . $$
To compute $\Aut_A(\lu,\varphi)$ it 
is convenient to identify $V$ with the underlying vector space of $\lu$,
i.e., we consider simply transitive representations $\varphi$ of 
$\lu$ on itself. For elements $D,g \in \GL(\lu)$ let us  define
$ (D,g): = g^{-1} D g D^{-1}$. The group $\Aut_A(\lu)$ 
may now  be characterized in terms of $\phi$ as follows: 
 
\begin{proposition} \label{autocharact} Let $A \leq \Aff(\lu)$ be a
transitive subgroup. Let  $\varphi: \lu \rightarrow \la(\lu)$ be a
a unipotent simply transitive representation of $\lu$ on itself. 
Let $g \in \Aut(\lu)$. Then $g \in \Aut_A(\lu,\varphi)$ if and only if 
$$ D g D^{-1} \in N_{GL_A}(\phi(\lu))$$ 
and 
$$  (D,g) \,  \phi_X \,  (D,g)^{-1} = g^{-1} \phi_{g X}\, g \;\;, \;  
\mbox{ for all $X \in \lg$ .}  $$
\end{proposition} 
\begin{proof} Let $h \in N_{GL_A}(\phi(\lu))$. If $\varphi = (\phi, D)$ then
the conjugate representation is $\varphi^h = (\phi^h,h D)$. 
The image $\Phi_h$ of $h$ in $\Aut_A(\lu,\varphi)$ is defined 
by the condition that 
$$  \varphi \;  \Phi_h \; = \; \varphi^h  \; .$$ 
This is equivalent to $\phi^h = \phi \,  \Phi_h$, and 
$h D = D \Phi_h$. By the latter condition it follows that   $$ D \Phi_h D^{-1}  = h \in 
N_{GL_A}(\phi(\lu)) \; .$$  The
proposition is immediate with the first
condition.  
\end{proof}

This gives rise to the following 
\begin{definition} \label{compatible}
Let $\phi: \lu \rightarrow \lgl(\lu)$ be a representation
of Lie algebras, and  $g \in \Aut(\lu)$. Then $g$ is called {\em compatible\/} with  $\phi$  
if $$     g^{-1} \phi_{g X} \, g = \phi_X\;\;, \;  
\mbox{ for all $X \in \lg$ .}  $$
\end{definition}

We remark that a particular special case is the 
adjoint representation of $\lu$ on itself, $\phi= \ad$.
By its very definition, $\Aut(\lu)$ is compatible with $\ad$. 

\subsubsection{Basic examples of unipotent simply transitive actions}
We come now to a first family of examples 
which was known to Cartan already. These examples are fully 
invariant and may be seen as a generalization of the natural 
simply transitive representation of an abelian Lie algebra. 

\begin{definition} We call a simply transitive representation $\varphi$ of $\lu$ 
{\em invariant\/} if $\Aut(\lu,\varphi) = \Aut(\lu)$. We call $\varphi$ 
invariant by the subgroup $H \leq \Aut(\lu)$ if $H \leq \Aut(\lu, \varphi)$.
\end{definition}

\begin{proposition} \label{2step}
Let $\lu$ be a 2-step nilpotent Lie algebra. Then $\lu$
has a natural invariant simply transitive representation.  
\end{proposition}   
\begin{proof} We take
$\phi= {1 \over 2} \, \ad$, and $D= \id_{\lg}$, where $\ad$ denotes
the adjoint representation. 
\end{proof}

This construction was generalized by Scheuneman (\cite{Scheune1}) 
as follows:
\begin{proposition} \label{derivation} 
Let $\lu$ be a nilpotent Lie algebra, which admits a
nonsingular derivation $D$.  Then $\lu$ has a simply transitive affine
representation $\varphi = ( \ad,D)$. The representation $\varphi$ 
is invariant by the centralizer of $D$ in $\Aut(U)$. 
\end{proposition}
\begin{proof} Since $\ad$ is compatible with $\Aut(\lu)$, the 
centralizer $$ \C_{Aut(\lu)}(D) = \{ g \in \Aut(\lu) \mid (g,D) = id_{\lu} \}$$
is contained in $\Aut(\lu, \varphi)$, by Proposition \ref{autocharact}.
\end{proof} 

\paragraph{Construction of derivations}
The proposition raises the particular problem to study the centralizers
of certain elements in $\Aut(\lu)$. 
A semisimple element 
$g \in \Aut(\lu)$ is called expanding if all eigenvalues have absolute value
$>1$. We construct now examples of Lie algebras 
which contain expanding elements in the center of the
Levi subgroups of $\Aut(\lu)$. (We call  a maximal reductive
subgroup of the linear algebraic group a Levi subgroup, cf.\ 
Theorem \ref{Mostow}.)\\

For preparation we have to introduce
a new concept. 

\begin{definition} \label{filtrations} 
A (positive) {\em filtration} on a Lie algebra $\lg$ is
a nested sequence of subspaces 
$$   \lg= \lg_1 \supseteq \lg_2 \supseteq \lg_3 \supseteq \cdots $$
such that $[ \lg_i , \lg_j ] \subseteq \lg_{i+j}$. The filtration is
called {\em invariant\/} if it is preserved by $\Aut(\lg)$. 
\end{definition}

For each positive filtration on $\lg$ there is an associated graded
Lie algebra 
$$\ac{\lg} = \bigoplus_{i=1, \ldots}\; \lg_i \rmod \lg_{i+1} $$ where
the Lie product on $\ac{\lg}$ is defined by setting, for $x \in  \lg_i$, 
$y \in  \lg_j$,  
$$   [\ac{x}, \ac{y} ] = \ac{[ x, y ]} \in \lg_{i+j} \rmod \lg_{i+j+1} \; . $$

A grading  $$ \lg = \oplus_{l=1, \ldots, k} V_l $$ is called 
invariant by a subgroup $G \leq \Aut(\lg)$ if every $g \in G$
preserves the decomposition.
If $\lg$ has a positive grading  $\lg = \oplus_{l=1, \ldots, k} V_l$
then the nested sequence of ideals $\lg_i = \oplus_{l=i, \ldots, k} V_l$ 
defines an associated positive
filtration. 

Note that every positive grading
is preserved by a one parameter group $l_{\lambda}$ of expanding
automorphisms in $\Aut(\lg)$. In fact, for $\lambda>1$,  $l_{\lambda}$ is given 
by  $$ l_{\lambda}(v)= \lambda^i v  \; , \;   v \in V_i \; . $$

\begin{proposition} Let $\lu$ be a nilpotent Lie algebra which 
has a positive grading with an invariant associated filtration.
Let $L$ be a Levi-subgroup
of the linear algebraic group $\Aut(\lu)$. Then $L$ contains a one parameter group of 
expanding automorphisms in its center.
\end{proposition}
\begin{proof}
Since the filtration $\lu_i$ of $\lu$ associated with the 
grading is invariant, 
$L$ acts on the factor spaces $\lu_{i} \rmod \lu_{i+1}$. 
Therefore $L$ acts by automorphisms on the
associated graded Lie algebra $\ac{u}$ such that the grading
is preserved. Since $L$ is reductive the action of $L$ on $\ac{u}$ 
is also faithful. 

Now the fact that $\lu$ is graded with respect to the filtration $\lu_i$
implies that there is an isomorphism of filtered Lie algebras 
$\pi: \lu \rightarrow \ac{\lu}$.
Hence it follows that $\Aut(\lu)$ has a Levi subgroup 
which preserves a grading compatible which is compatible 
with the original filtration. 
Since all Levi subgroups are conjugate in $\Aut(\lu)$, every Levi subgroup 
$L$ preserves a corresponding positive grading of $\lu$. 
The one parameter family $l_{\lambda}$ of expanding 
automorphisms which belongs to this positive grading commutes with $L$.
\end{proof} 

As a consequence we obtain
\begin{theorem} \label{graded}   
Let $u$ be a nilpotent Lie algebra which 
has an invariant positive grading, and let $L$ be a Levi-subgroup
of $\Aut(U)$. Then there exists  an affine simply transitive
representation which is invariant by $L$.
\end{theorem}

\begin{remark} It is easily checked from certain lists of 
nilpotent Lie algebras that every nilpotent Lie algebra of
dimension less or equal to five has an invariant positive
grading. 
\end{remark}

\paragraph{Scheuneman representations for 3-step  nilpotent Lie algebras}
A little more refined construction which builds on 
the previous cases was given by Scheuneman in \cite{Scheune2}.
The construction is compelling because it works for every
3-step nilpotent Lie algebra. 

\begin{theorem} \label{3step} 
Let $u$ be a 3-step nilpotent Lie algebra, 
and let $L$ be a Levi-subgroup of $\Aut(U)$. Then there exists 
an affine simply transitive representation which is invariant by $L$.
\end{theorem}
\begin{proof} Following \cite{Scheune2}, we consider the filtration 
$$\lu= \lu_1 \supseteq \lu_2 \supseteq \lu_3 \supseteq \{ 0 \} $$ 
which is given by the descending central series, i.e., $\lu_2 = [ \lu, \lu ]$,
$\lu_3 = [ \lu, \lu_2]$. We choose a compatible
decomposition of $\lu = \oplus_{l=1, \ldots, 3}\,  V_l$ as 
a direct sum of vector spaces, $V_3 = \lu_3$. We define 
a diagonal element $g= g_{\alpha,\beta,\gamma} \in \GL(\lu)$ 
which preserves the grading by $g v = \alpha v$, for $v \in V_1$, 
$g v = \beta v$, for $v \in V_2$, and $g v = \gamma v$, for $v \in V_3$. 
Correspondingly, we define $D= D_{r,s,t}$. It is easy to 
see (compare \cite{Scheune2}) that there exist $\alpha,\beta,\gamma>0$, 
and $r,s,t>0$ so that $D$ is a derivation for the
representation $\phi = \ad^g$ of $\lu$ on itself. The pair $(\phi,D)$
defines then a simply transitive representation $\varphi$ of $\lu$. 
Since the filtration
we chose is invariant, we may assume that the spaces $V_i$ are
invariant by automorphisms of a Levi subgroup $L$ of $\Aut(\lu)$.
In this case $g$ commutes with $L$, and the representation $\phi$ is
compatible with $L$ in the sense of Definition \ref{compatible}.
Since also $D$ commutes with $L$, 
it follows that $L \leq \Aut(\lu,\varphi)$, by Proposition \ref{autocharact}. 
\end{proof}

\begin{remark} When applied to the 2-step nilpotent case
the method recovers just the natural invariant simply
transitive action.
\end{remark} 


\chapter{Geometric Realization Problems} \label{realizations}
\section{Nielsen's realization problem} \label{Nielsen}
Nielsen considered in {\cite{Nielsen}} the following question which
is now called the {\em Nielsen realization problem}. Namely, 
{\em can any finite group $G$ of isotopy classes of
self-homeomorphisms of a surface be realized as a finite group 
of self-homeomorphisms?} 
As it turned out, Nielsen's question is equivalent to the following: 
{\em Let $\Delta$ be a finite effective extension of the fundamental group 
$\pi_1(S)$ of a surface $S$. 
Is it true, that $\Delta$ acts discontinuously on the plane?}

Geometry gives a positive answer to this question. 
For example, every effective finite extension of
$\bbZ^2$ acts as an Euclidean crystallographic group. These are
the well known so called wall-paper groups. Therefore the
realization problem for Euclidean surfaces has a positive solution. 
For hyperbolic
surfaces the problem was settled not long ago. 
Kerckhoff's  \cite{Kerck} celebrated 
fixed point theorem states:  
Any finite group of mapping classes 
of a closed orientable surface $M_g$ of genus $g \geq2 $ 
has a fixed point in the Teichm\"uller space of hyperbolic structures.
A particular consequence of the fixed point theorem is:
Let $\Delta$ be a finite effective extension of $\pi_1(M_g)$,
then $\Delta$ is isomorphic to a discontinuous group of hyperbolic motions.
These results give a beautiful example of the significance of geometric methods 
in topology. \newline 

Recall, that a homomorphism
$\alpha: F \rightarrow \Out(\Gamma)$ 
from a group $F$ to the outer automorphism
group of a group $\Gamma$ is called an {\em abstract kernel for $\Gamma$}.
Every group $F$ of homeomorphisms of a manifold 
naturally defines an abstract kernel $\alpha: F \rightarrow \pi_1(M)$.
If this is the case, we say that $\alpha$ is realized by a group of homeomorphisms.
The {\em topological realization problem\/} may be posed    
as follows: {\em Let $M$ be a manifold and let $\alpha$ be a finite
abstract kernel for $\pi_1(M)$. Is it possible to realize $\alpha$
by a group of homeomorphisms of $M$?} 

Under suitable topological conditions on the manifold $M$, 
for example if $M$ is aspherical, a necessary condition for 
$\alpha$ to be realizable as a group of homeomorphisms is the 
existence of a finite {\em effective\/} 
extension group $\Delta$ for $\Gamma=\pi_1(M)$ which induces $\alpha$.
The group $\Delta$ is then an {\em algebraic realization\/} of the kernel $\alpha$. 
As illustrated above in the case of surfaces, 
geometry may help to realize $\Delta$ as properly discontinuous group
of homeomorphisms of the universal cover $\tilde{M}$ of $M$. 
Still, {\em geometric realization questions\/} are 
interesting in their own right. An answer to such a question may show for
example how well a particular geometry 
is adapted to the topology of a manifold. \newline

A classic example which extends beyond the surface case is the following: 
Bieberbach proved in 1912 his famous theorems on the structure of
the Euclidean crystallographic groups.  
Later Burckhardt and Zassenhaus showed that Euclidean crystallographic 
groups are characterized by their algebraic properties. In particular,
their remark implies (see \cite{Zieschang}[\S 25]) a positive answer to the 
geometric realization problem
for {\em Euclidean crystallographic groups}: 

\begin{theorem} \label{Euclrealization}  
Let $\Gamma$ be isomorphic to an Euclidean crystallographic group. If
$\Delta$ is a finite effective extension of\/ $\Gamma$ then $\Delta$ 
is also isomorphic to an Euclidean crystallographic group. 
\end{theorem} 

Since, by Bieberbach's second theorem, flat Euclidean manifolds
are determined by their fundamental group up to affine 
diffeomorphism this implies: 

\begin{theorem} \label{Euclrealization2} 
Let $M$ be compact Euclidean space form, and $\alpha$ an injective abstract kernel
for $\pi_1(M)$. Then the following are equivalent: 
\begin{itemize}
\item[i)] $\alpha$ has an algebraic realization $\Delta$,
\item[ii)]$\alpha$ can be realized by a
group of affine diffeomorphisms of $M$.   
\end{itemize}
\end{theorem} 


\paragraph{The realization problem for affine crystallographic groups}
In our work we are concerned with the realization 
problem for affine crystallographic groups.  
We restrict ourselves here to the case that $\Gamma$ is a virtually nilpotent ACG.\footnote{Most of our results here
generalise to virtually polycyclic ACGs in the context of infrasolvmanifolds.
See \cite{BauesU,BauesI} for some of the general constructions.}

Every affine space form $M$ with fundamental group isomorphic 
to $\Gamma$ is finitely covered by a nilmanifold. (A nilmanifold
is diffeomorphic to the quotient of a nilpotent Lie group by a closed
subgroup. $M$ itself has an infranilmanifold structure.) 
For certain aspherical manifolds the topological realization question 
was solved by Lee and Raymond (cf.\ \cite{LeeRaymond1}) with 
methods from the theory of Seifert fiber spaces. 
One particular result is: 

\begin{theorem} Let $M$ be an infranilmanifold, $\alpha$ an abstract kernel
for $\pi_1(M)$. If there exists an algebraic realization of $\alpha$ then
$\alpha$ may be realized as a group of homeomorphisms of $M$. 
\end{theorem}

However, the following {\em geometric question\/} remains: 
{\em Which of the finite quotient manifolds of an affine
space form $M$ are again an affine space form?}  
Our basic result shows that if $M$ has a virtually nilpotent
fundamental group then the solution to the geometric realization problem 
is a question about the fixed points 
of abstract kernels on the deformation space. 

\begin{theorem} \label{Affinerealization1}
Let $\Lambda$ be virtually nilpotent and 
isomorphic to an affine crystallographic group of type $A$. Then a
finite effective extension $\Delta$ of\/ $\Lambda$ is isomorphic 
to an affine crystallographic group of type $A$ if and only if the kernel 
associated to the extension has a fixed point in 
the deformation space $\Def(\Lambda,A)$. 
\end{theorem}

We will also consider the question for which crystallographic groups
the analogue of Theorem \ref{Euclrealization} holds. These are
precisely those affine crystallographic groups $\Gamma$ which
satisfy the following {\em realization property\/}.
\begin{definition}
We say that $\Gamma$ has the affine realization property if
every finite effective extension $\Delta$
of $\Gamma$ is isomorphic to an affine crystallographic group. 
\end{definition} 

We remark that if the fundamental group $\pi_1(M)$ of an affine space
form $M$ satisfies the realization property then,
by some standard reasoning\footnote{using that 
compact complete affine manifolds with isomorphic fundamental group are
polynomially diffeomorphic, see \cite{FriedGoldman}}, the following generalization 
of Theorem \ref{Euclrealization2} holds: 

\begin{theorem} \label{Affinerealization2}
Let $M$ be a compact affine space form with a virtually nilpotent
fundamental group $\pi_1(M)$, and $\alpha$ an injective abstract kernel
for $\pi_1(M)$. If $\pi_1(M)$ has the realization property 
then the following are equivalent: 
\begin{itemize}
\item[i)] $\alpha$ has an algebraic realization $\Delta$.
\item[ii)]$\alpha$ can be realized by a
group of polynomial diffeomorphisms of $M$.   
\end{itemize}
\end{theorem} 

Theorem \ref{Affinerealization1} implies that $\Gamma$ satisfies
the realization property 
if the deformation space $\Def(\Gamma,\Aff(V))$ is convex
in the sense of Definition \ref{convex}. In particular, 
compare Corollary \ref{nilexamples}, 
the examples of chapter \ref{defspaces} show that some natural
classes of virtually nilpotent crystallographic 
groups satisfy the realization property. 
 
\begin{corollary} Let $\Gamma$ be a virtually nilpotent 
affine crystallographic group.
If the deformation space  $\Def(\Gamma, \Aff(V))$ is convex 
then $\Gamma$ has the realization property. 
\end{corollary}
This fact should provide us with enough motivation to study the
deformation spaces for affine space forms. 




%


\section{The Realization of finite extensions for \ftn-groups} \label{realproblems}

Let $\Gamma$ be a group which is isomorphic to an ACG of type $A$, 
and let  $\beta: F \longrightarrow \Out(\Gamma)$ be an abstract kernel for $\Gamma$. 
As a subgroup of  $\Out(\Gamma)$,  $\beta(F)$ acts on the deformation 
space $\Def(\Gamma,A)$.
Let us assume that there exists a finite normal extension $\Gamma \leq_f \Delta$ which 
realizes the abstract kernel $\beta$. This means that $\beta$ coincides with
the natural homomorphism $F= \Delta/\Gamma \longrightarrow \Out(\Gamma)$ which is
associated to the extension.  
The following  is a fundamental observation: 

\begin{proposition} \label{fundamental}
Let $\beta: F  \longrightarrow \Out(\Gamma)$
be the kernel associated to the normal extension  $\Gamma \leq_f \Delta$.
Assume that there exists a crystallographic homomorphism 
$\rho_{\Delta} \in \Hom_c(\Delta,A)$. 
Let  $\rho$ denote the restriction of $ \rho_{\Delta}$ to $\Gamma$.  Then 
$ \rho$ is in $\Hom_c(\Gamma,A)$, 
and $[\rho\,] \in \Def(\, \Gamma,A)$ is a fixed point for the action of
$F$ on $\Def(\, \Gamma,A)$.  
\end{proposition}
\begin{proof}
Clearly, the restriction of $\rho_{\Delta}$ to $\Gamma$ remains
crystallographic. For $g \in F $, we have 
$\beta(g)[ \rho \,  ] = [ \rho^{\phi_g} ]$, where 
$\phi_g \in \Aut(\Gamma)$ represents $\beta(g) \in \Out(\Gamma)$.  
But, since $\Delta$ realizes the kernel $\beta$,  there exists $\delta \in \Delta$
such that $\beta(g)$ is represented by conjugation with $\delta$. Therefore,
we may assume that $\phi_g(\gamma) = \delta \gamma \delta^{-1}$, 
for all $\gamma \in \Gamma$. It follows that 
$ \rho^{\phi_g} (\gamma) =   \rho_{\Delta}(\delta) \rho(\gamma)  \rho_{\Delta}(\delta)^{-1}$.
Since  $\rho_{\Delta}(\delta) \in A$, we conclude 
$\beta(g)[ \rho \,  ]=  [\rho^{\phi_g}] = [\rho\,]$.
\end{proof}

From now on let us assume that $\Gamma$ is a \ftn-group.
\begin{definition}
A normal extension group $\Delta$ of $\Gamma$ is called
{\em effective\/} if the associated kernel 
$\beta: \Delta/\Gamma \longrightarrow \Out(\Gamma)$ is
an injective homomorphism. 
\end{definition}
We are going to  show that the converse of Proposition \ref{fundamental}
holds for all finite effective extensions $\Delta$ of $\Gamma$.
Namely, the existence of fixed points for $\beta$
on the deformation space of $\Gamma$ 
is the only obstruction to realize $\Delta$ 
as an affine crystallographic group.

\begin{theorem} \label{realizationA}  
Let  $\Gamma$ be a \ftn-group which is isomorphic to an
affine crystallographic group of type $A$, and let $\Gamma \leq_f \Delta$
be an effective finite extension with
associated  abstract kernel $\beta: F \longrightarrow \Out(\Gamma)$.
Then $\Delta$  may be realized as an affine crystallographic group of type $A$
if and only if the kernel $\beta$ has a fixed point in $\Def(\Gamma, A)$.
\end{theorem} 

\begin{remark} In the course of the proof of Theorem \ref{realizationA}, we
will show a stronger statement: Every $\rho \in \Hom_c(\Gamma,A)$ which is a fixed point
for $\beta$ may be extended to a homomorphism 
$\rho_{\Delta} \in \Hom_c(\Delta, A)$.  
\end{remark}

Theorem \ref{Affinerealization1} follows from Theorem  \ref{realizationA}.
The proof of Theorem  \ref{Affinerealization1} 
will be given below in section \ref{applications}.
The purpose of the remainder of this section is to prove Theorem  \ref{realizationA}. 
We will proceed as follows: \newline 

The first step  is,  to consider the related problem of 
realizing a finite extension $\Delta$ of $\Gamma$ as a Zariski dense subgroup
of  a linear algebraic group. The solution to this problem   (Proposition \ref{algreal})
leads us to the construction of an algebraic hull for $\Delta$ which
satisfies functorial properties with respect to the representations
of $\Delta$. 
As one application of the construction
we deduce a splitting result for the extension $\Gamma \leq \Delta$ (Proposition \ref{splitting}). 
Another application is a certain 
characterization of the finite effective nilpotent extensions of  $\Gamma$.
(Corollary \ref{nilextensions})

The second step is the proof of the affine crystallographic realization for the
finite effective extensions $\Delta$ of $\Gamma$. This
breaks into two distinct parts: The case that $\Delta$ itself  is a  \ftn-group, and
the case that the extension $\Gamma \leq \Delta$ splits. \newline 

We will use frequently certain results from  
the theory of \ftn-groups.  For  a general reference on the relevant
aspects of the theory, see  \cite[Chapter 2]{EMSII}, \cite[Chapter II]{Raghunathan}, and 
in particular the book \cite{Segal}. We start with some preparatory
material.  

\subsection{The algebraic hull of an extension}
\paragraph{Automorphisms and the Malcev hull}
Let $\Gamma$ be a \ftn-group.
It is known 
that there exists a $\bbQ$-defined unipotent algebraic group $\bU$
(called the Malcev completion of $\Gamma$),
and an embedding $j: \Gamma \rightarrow \bU_{\bbQ}$ such
that $\Gamma$ is  a Zariski-dense 
subgroup in $\bU$.
Every
automorphism of $\Gamma$ extends to a unique $\bbQ$-defined 
automorphism of $\bU$. In fact, the automorphism group $\Aut(\bU)$ of $\bU$ is a
$\bbQ$-defined linear algebraic group, in such a way that 
$\Aut(\bU)_{\bbQ} = \Aut(\bU_{\bbQ})$,
and the extension defines
an injective homomorphism   
$$ \epsilon_j: \Aut(\Gamma) \longrightarrow  \Aut(\bU)_{\bbQ} \; $$
which satisfies $\epsilon_j(\phi) \,  j   = j  \, \phi $, for all $\phi \in \Aut(\Gamma)$.  
By factorization there
exists a homomorphism $o_j: \Aut(\Gamma)   \rightarrow \Out(\bU)$
such that the diagram 
$$ \xymatrix{
\Aut(\Gamma)  \ar[d]  \ar[r]^{\epsilon_j}  & \Aut(\bU)   \ar[d] \\
\Out(\Gamma)  \ar[r]^{o_j}  & \Out(\bU) }
$$
is commutative. Note that $o_j$ is, in general, not injective. 
The preimage of its kernel is 
$\epsilon_j^{-1}(\Inn(\bU))$. 
Let us therefore denote  
$$ \Inn_{\bbQ}(\Gamma)= \epsilon_j^{-1}(\Inn(\bU)) \text{
and } \Out_{\bbQ}(\Gamma)= \Aut(\Gamma) / \Inn_{\bbQ}(\Gamma) \; . $$

\begin{definition} \label{radicably}
Let $\beta: F \rightarrow \Out(\Gamma)$ be an abstract kernel.
$\beta$ is called {\em radicably effective} if $\beta$ induces an embedding 
$F \hookrightarrow \Out_{\bbQ}(\Gamma)$. $\beta$ is called {\em radicably trivial} if 
$\beta$ induces the trivial homomorphism to $\Out_{\bbQ}(\Gamma)$.
\end{definition} 

Let $\Gamma \leq \Delta$ be a finite
extension with $\Gamma$ normal in $\Delta$, 
and let $c: \Delta \rightarrow \Aut(\Gamma)$ be the homomorphism
which is induced by conjugation. We put $\Gamma^r =  c^{-1}(\Inn_{\bbQ}(\Gamma))$.
So $\Gamma \leq \Delta$ is radicably
effective if and only if $\Gamma^r =\Gamma$.  
We remark that if  the extension is finite and effective $\Gamma^r$ 
is a nilpotent normal subgroup  in $\Delta$. In fact, it follows 
(Corollary \ref{rkernel}) that $\Gamma^r$ is the Fitting subgroup of $\Delta$.
\newline 

Later we will need the following lemma:
\begin{lemma} \label{invariance} 
Let $\mu \leq \Aut(\bU_{\bbQ})$ be a finite subgroup,
and $\, \Gamma \leq \bU_{\bbQ}$ a finitely generated subgroup.
Then there exists $\Gamma_f \leq \bU_{\bbQ}$ such that $\Gamma \leq_f \Gamma_f$,
and $\, \Gamma_f$ is normalized by $\mu$.
\end{lemma}
\begin{proof} Let $\lie{u}$ be the Lie algebra of $\bU$.  $\lie{u}$ is defined over $\bbQ$, and
the exponential map $\exp: \lie{u} \rightarrow \bU$ identifies
$\lie{u}_{\bbQ}$ with $\bU_{\bbQ}$.  The group $\mu$ acts as
a subgroup of automorphisms in $\GL( \lie{u}_{\bbQ})$
on $\lie{u}_{\bbQ}$, and to show that $\mu$ normalizes
$\Gamma \leq \bU_{\bbQ}$ is equivalent to show that $\mu$ normalizes
the set $\log \Gamma \subseteq \lie{u}_{\bbQ}$. 
Replacing, if necessary, $\Gamma$ with a finite extension
group, we may assume (\cite{Segal}[\S 6 B, Theorem 3])
that $\Gamma$ is a lattice subgroup, such that 
$\Lambda= \log \Gamma  = \bbZ \log \Gamma$ is a full lattice in $\lie{u}_{\bbQ}$.
Since $\mu \leq\GL( \lie{u}_{\bbQ})$ is  finite there exists an integer $m$  
sucht that,  for all $g \in \mu$, $g \, {1 \over m} \Lambda \subseteq  {1 \over m} \Lambda$.
Let $\Gamma^{1\over m}$ be the group which is generated by the set 
$\{ \exp u \mid u \in {1 \over m} \Lambda \}$. 
Now $\Gamma_f = \Gamma^{1\over m}$ is a finite extension of $\Gamma$, 
and is normalized by $\mu$.
\end{proof} 

\paragraph{The Malcev extension functor}
As before let the $\bbQ$-defined algebraic 
group $\bU$ denote the Malcev completion of $\Gamma$.  
We need to consider the extension functor ({\em Malcev-rigidity\/}) from the 
unipotent representations of $\Gamma$ to representations 
of $\bU$. Remark, that this is 
just a special case of Propostion \ref{pchullextension}.

\begin{proposition} \label{Malcev} 
Let $\rho: \Gamma \rightarrow \bH$ be
a homomorphism of $\Gamma$ to a unipotent $\bbQ$-defined linear algebraic group
$\bH$. Then there exists a unique morphism of algebraic groups 
$$ \rho_{\bU}: \bU \longrightarrow \bH $$
which extends $\rho$. If $\rho(\Gamma) \leq   \bH_k$, where $k$ is a subfield
of\/ $\bbC$,  then $\rho_{\bU}$ is defined over $k$.   
\end{proposition} 
\begin{proof} Consider the subgroup 
$D= \{(\gamma, \rho(\gamma)) \mid \gamma \in \Gamma \}$ of
the product $\bU \times \bH$. Let $\pi_1, \pi_2$ denote the
projection morphisms on the factors of the product. Let $\bD$ be
the Zariski-closure of $D$. By \cite[Theorem 2.10]{Raghunathan} we have
$\dim \bD = \rank \Gamma = \dim \bU$, and it follows
that $\pi_1:  \bD \rightarrow \bU$ is an isomorphism. Hence 
$\rho_{\bU}   = \pi_2 \circ \pi_1^{-1}$ is the unique extension.
If $\rho(\Gamma) \leq \bH_k$ then $\bD$ is $k$-defined, 
and hence also the morphism $\rho_{\bU}$. 
\end{proof} 

At some point later we consider the spaces 
$\Hom(\Gamma, \bH_{\bbR})$ as topological 
spaces equipped with the compact-open topology. 
(Compare section \ref{topology}) Now let $\Gamma^*$ be a \ftn-group which is a finite
extension of $\Gamma$. There exists a unique embedding
of $\Gamma^*$ into $\bU_{\bbQ}$ which is the identity on $\Gamma$,
see \cite{Segal}[\S 6A]. We note:

\begin{corollary} 
\label{Extfunctor} Let $\bH$ be a unipotent $\bbQ$-defined linear algebraic group.
Then every $\rho \in\Hom(\Gamma, \bH_{\bbR})$ has a unique 
extension $\rho_{\Gamma^*} \in \Hom(\Gamma^*,\bH_{\bbR})$.
The extension functor induces a homeomorphism
$$ e: \Hom(\Gamma, \bH_{\bbR}) \longrightarrow  \Hom(\Gamma^*, \bH_{\bbR}) \; .$$
\end{corollary}
\begin{proof} We use the embedding $\Gamma^* \hookrightarrow  \bU_{\bbQ}$.
Malcev rigidity implies that  every $\rho \in  \Hom(\Gamma, \bH_{\bbR})$ 
extends to a homomorphism $\rho_{\bU}: \bU_{\bbR} \rightarrow \bH_{\bbR}$.
The restriction of $\rho_{\bU}$ to $\Gamma^*$ is the  unique extension 
$ e(\rho)= \rho_{\Gamma*} \in  \Hom(\Gamma^*, \bH_{\bbR})$. 
The map $e$ is then clearly a bijection. 

Let $G$ be a discrete group.  For any  $V \subset  \bH_{\bbR}$
which is open in the Hausdorff-topology on  $\bH_{\bbR}$, and any $\gamma \in G$
we put $$ {\cal U}_G (\gamma,V) \, = \, \{ \rho \in \Hom(G, \bH_{\bbR}) 
 \mid \rho(\gamma) \in V \} \; . $$ 
The compact open topology on 
$\Hom(G,  \bH_{\bbR})$ is then generated by all sets $ {\cal U}_G (\gamma,V)$.
It  is immediate that the restriction map $e^{-1}$ is continuous. 
To  prove that $e$ is continuous we take the following approach: 
We can assume that $\Gamma^* =  \Gamma^{1\over m}$, for some $m \in \bbN$.  
By \cite{Segal}[\S 6B, Proposition 2] there exists $s \in \bbN$ such 
that for every $\tilde{\gamma} \in  \Gamma^*$,  $\tilde{\gamma}^s \in \Gamma$.
We put  
$V^s = \{ v^s \mid  v\in V \}$, $\gamma = \tilde{\gamma}^s . $
Now $$ e^{-1} (\, {\cal U}_{{\Gamma^*}}(\tilde{\gamma},V)) \, 
= \{ \rho \mid \rho(\gamma) \in V^s \} 
=  \,  {\cal U}_{\Gamma} (\gamma, V^s) \; .$$ Hence,  $e$ is continuous.    
\end{proof} 


\paragraph{The algebraic hull of an effective extension}
Let $\Gamma$ be a \ftn-group. Then $\Gamma \leq \bU_{\bbQ}$,  where
$\bU$ is the Malcev completion of $\Gamma$. In fact, we may assume that 
$\Gamma \leq \bU_{\bbZ}$, i.e., $\Gamma$  is contained in the subgroup  
$\bU_{\bbZ}$ of integer points of $\bU$. Now let $\Delta$ be a finite normal
extension group of $\Gamma$, i.e., $\Gamma$ is a normal subgroup $|\Delta/\Gamma| < \infty$.  
We can also realize the group $\Delta$ as a
subgroup of a linear algebraic group (in fact, the construction realizes $\Delta$
as an arithmetic group,  as follows from \cite{GP}):  

\begin{proposition} \label{algreal} 
There exists a $\bbQ$-defined linear algebraic group $I(\bU,\Delta)$
with  $\bU$ its component of identity, and an embedding
$\psi: \Delta \rightarrow I(\bU,\Delta)_{\bbQ}$ which is the identity on $\Gamma$, such that
$I(\bU,\Delta)= \psi(\Delta)\bU$ and $\psi(\Delta) \cap \bU= \Gamma$.
\end{proposition} 
\begin{proof} We may assume that 
$\Gamma \leq \bU_{\bbZ}$, i.e., $\Gamma$  is contained in the subgroup  
$\bU_{\bbZ}$ of integer points of $\bU$. 
Let $\Delta= \Gamma r_1  \cup \, \cdots \, \cup \, \Gamma r_m$ be
a decomposition of $\Delta$ into left cosets. 
By Malcev-rigidity, conjugation with $r_i$ on $\Gamma$ extends to
$\bbQ$-defined rational homomorphisms $f_i$ of $\bU$. 
A straightforward application of the construction used to 
prove Proposition 2.2 in \cite{GP} implies the result. 
\end{proof}
In the case that the extension $\Gamma \leq \Delta$ is
effective we may further refine
Proposition \ref{algreal}. Let $\Fitt(\Delta)$ be the 
Fitting subgroup of $\Delta$, i.e., its maximal nilpotent normal
subgroup. We show that  $\Fitt(\Delta)$ is a \ftn-group and 
embeds into $\bU_{\bbQ}$.  
\begin{proposition} \label{refreal}
If the finite extension $\Gamma \leq \Delta$ is effective then 
there exists a $\bbQ$-defined linear algebraic group $I^*(\bU,\Delta)$ with $\bU$ its
component of identity, and an embedding $\psi: \Delta \rightarrow I^*(\bU,\Delta)_{\bbQ}$
which is the identity on $\Gamma$, such that $I^*(\bU,\Delta)= \psi(\Delta)\bU$ and
$\psi(\Delta) \cap \bU = \psi(\Fitt(\Delta))$. Moreover the centralizer of $\bU$ in 
$I^*(\bU,\Delta)$ is contained in $\bU$.
\end{proposition}
\begin{proof} By Proposition \ref{algreal} we may assume that
 $\Delta$ is a subgroup of $\bG_{\bbQ}$, where $\bG$ is a linear algebraic group
with $\ur(\bG)=\bU$,
and $\Gamma \leq \bU$. 
Let $\Gamma^*= \Fitt(\Delta)$. Since $\Gamma$ is a subgroup
of finite index in $\Gamma^*$, the group 
$\Gamma^*_s =\{ \gamma_s \mid \gamma \in \Gamma^* \}$ is finite.
Since $\Gamma^*$ is normal in $\Delta$, $\Gamma^*_s$ is normalized
by $\Gamma$. The centralizer of $\Gamma^*_s$ is a Zariski-closed 
subgroup of $\bG$ which contains a finite index subgroup of $\Gamma$. 
Therefore the centralizer $\C_{\bG}(\Gamma^*_s)$ contains $\bU$, in particular 
$\Gamma^*_s$ centralizes $\Gamma$.  
Consider the homomorphism $\psi_u: \Gamma^* \rightarrow \bU_{\bbQ}$ given
by $\gamma \mapsto \gamma_u$. Since the extension
$\Gamma \leq \Gamma^*$ is effective, the homomorphism $\psi_u$ is
injective. Therefore $\Gamma^*$ is a \ftn-group, and embeds as a
subgroup of $\bU_{\bbQ}$ containing $\Gamma$. We obtain $I^*(\bU,\Delta)$
by applying Proposition \ref{algreal} to the extension $\Gamma^* \leq 
\Delta$. 

The centralizer of $\bU$ in $I^*(\bU,\Delta)$ splits
as a direct product $\C_{I^*(\bU,\Delta)}(\bU) = H \times \Z(\bU)$ 
where $H$ is a finite group of semisimple elements. It follows
that the set $X_{\Delta} = \{ \gamma \in \Delta \mid \gamma_s \in H \}$
is a normal subgroup of $\Delta$, and the map $\psi_u: \gamma \mapsto \gamma_u$
is a homomorphism on  $X_{\Delta}$. Since the extension $\Gamma \leq \Delta$
is effective and $\Gamma \leq \bU$, $\psi_u$ is injective on  $X_{\Delta}$. 
Therefore $X_{\Delta}$ is nilpotent
and hence contained in $\Fitt(\Delta)$. 
By construction $\Fitt(\Delta)$ is unipotent in  $I^*(\bU,\Delta)$
and hence no $\gamma_s \neq 1$ 
centralizes $\bU$. Since $I^*(\bU,\Delta)= \Delta \bU$, $H$ is
trivial.
\end{proof}

Let $c: \Delta \rightarrow \Aut(\Gamma)$ be the homomorphism
which is induced by conjugation, and 
$\Gamma^r =  c^{-1}(\Inn_{\bbQ}(\Gamma))$. Then
there are finite extensions
$$  \Gamma \,  \leq \, \Gamma^r \,  \leq  \, \Delta \;  . $$
Recall that the extension 
$\Gamma \leq \Delta$ is called radicably effective (see  Definition \ref{radicably}) 
if and only if $\Gamma^r =\Gamma$.
\begin{corollary} \label{rkernel} Let $\Gamma$ be a \ftn-group and
$\Gamma \leq \Delta$ an effective finite extension. Then $\Gamma^r = \Fitt(\Delta)$, and in
particular, the extension $\Fitt(\Delta) \leq \Delta$ is radicably effective. 
\end{corollary}
\begin{proof}
We consider $\Gamma^r \leq  I^*(\bU,\Delta)$.
Since $c(\Gamma^r) \leq \Inn_{\bbQ}(\Gamma)$ consists of unipotent automorphisms the
semisimple parts of the elements of $\Gamma_r$ centralize $\bU$. 
Hence, by Proposition  \ref{refreal},  $\Gamma^r \leq U$ and therefore 
$\Gamma^r = \Delta \cap \bU =\Fitt(\Delta)$. 

Since $\Gamma \leq \Fitt(\Delta)$, the extension
$\Fitt(\Delta) \leq \Delta$ is effective. 
Since $\Fitt(\Delta)$ is a \ftn-group, 
the first part
of the corollary implies that $\Fitt(\Delta)^r = \Fitt(\Delta)$. 
\end{proof}

\begin{corollary} \label{nilextensions} 
Let $\Gamma$ be a \ftn-group and
$\Gamma \leq \Delta$ an effective finite extension. 
Then $\Delta$ is
a nilpotent group if and only if the extension is radicably
trivial. If $\Delta$ is nilpotent then $\Delta$ is a \ftn-group.  
\end{corollary}

The following lemma shows 
that the class of finite effective extensions of \ftn-groups
is closed under the operation of taking effective extensions. 
 
\begin{lemma} \label{eclosed}
Let $\Lambda$ be a finite effective extension
group of the \ftn-group $\Theta$, and $\Delta$ a finite
effective extension of $\Lambda$. Let $\Gamma= \Fitt(\Delta)$.
Then $\Gamma$ is a \ftn-group and 
the extension $\Gamma \leq \Delta$ is effective. 
\end{lemma}
\begin{proof}
We may assume that $\Theta = \Fitt(\Lambda)$ and 
embed $\Delta$ as a subgroup of a linear algebraic group 
$I(\bU,\Delta)$ such that $\bU \cap \Delta = \Theta$. 
It follows that $\Z_{\Delta}(\Theta)_s=  \Z_{\Delta}(\bU)_s$ 
is a finite subgroup of $I(\bU,\Delta)$ which is normalized by
the Zariski-closure $\ac{\Lambda}$ of $\Lambda$. 
Since $\Lambda$ is a normal subgroup of $\Delta$,
$\ac{\Lambda}$ is also normalized by $\Z_{\Delta}(\Theta)_s$.
Therefore 
$C= [\Z_{\Delta}(\Theta)_s, \Lambda] \leq \ac{\Lambda} \cap \Z_{\Delta}(\Theta)_s$.
Since $C \leq \ac{\Lambda} $ centralizes $\bU$, Proposition \ref{refreal} implies 
that $C$ is unipotent. Hence $C$ is trivial,
and $\Z_{\Delta}(\Theta)_s$ centralizes $\Lambda$. 
Since $\Delta$ is an effective extension of $\Lambda$ the map
$\gamma \mapsto \gamma_u$ is an injective homomorphism on 
$\Z_{\Delta}(\Theta)$. Therefore $\Z_{\Delta}(\Theta)$ is torsionfree, 
and also $\tilde{\Gamma}=\Z_{\Delta}(\Theta)\Theta$ is a $\ftn$-group.
Since the extension $\tilde{\Gamma} \leq \Delta$ is effective,
$\Delta$ is a finite effective extension of a \ftn-group.
The lemma follows.
\end{proof}

\begin{definition}
Let $\Delta$ be an  effective finite  extension of some \ftn-group.
A $\bbQ$-defined linear algebraic group $\bU_{\Delta}$ is called an
{\em algebraic hull for $\Delta$\/}, if the following hold: 
$\Delta$ is a Zariski-dense subgroup of $\bU_{\Delta}$,
so that $\Delta \leq \bU_{{\Delta}, \bbQ}$, 
$\Fitt(\Delta)$ is contained in the unipotent radical $\ur(\bU_{\Delta})$,
and  $\dim \ur(\bU_{\Delta}) = \rank  \Fitt(\Delta)$.\footnote{See Definition \ref{pchull}, for generalisation in the more general context
of finite extensions of polycyclic groups.}
\end{definition}

\begin{corollary} Every $\Delta$ as above admits an algebraic hull, and the
algebraic hull is unique up to $\bbQ$-isomorphism.\footnote{In fact, every 
virtually polycyclic \wfn-group has an algebraic hull. See \cite{BauesI}, and
chapter 1 of this article.}
\end{corollary}

The existence of the algebraic hull for $\Delta$ is proved by Proposition
\ref{refreal}. In fact $\bU_{\Delta} = I^*(\bU,\Delta)$, where $\bU$ is
the Malcev completion of $\Fitt(\Gamma)$. The uniqueness of the algebraic
hull follows from the following functorial property.

\begin{proposition} \label{hullextension}
Let $\bU_{\Delta}$ be an algebraic hull for $\Delta$,  
$\bG$ a $\bbQ$-defined linear algebraic group. Then every
homomorphism $\rho: \Delta \rightarrow \bG_{\bbQ}$, 
such that $\rho(\Fitt(\Delta))$ is 
contained in a $\bbQ$-defined unipotent subgroup $\bH$  of $\bG$,
extends uniquely to a $\bbQ$-defined homomorphism
$\rho_{\bU_{\Delta}}: \bU_{\Delta} \rightarrow \bG$.  
\end{proposition}
\begin{proof} 
Let $\Gamma = \Fitt(\Delta)$, and let $\bU$ be the unipotent
radical of $ \bU_{\Delta}$.  
By the Malcev extension property 
the induced homomorphism $\rho: \Gamma  \rightarrow \bH_{\bbQ}$
extends to a $\bbQ$-defined homomorphism 
$\rho_{\bU}:  \bU \rightarrow \bH$. 
Now every $g \in \bU_{\Delta}$ may be written as a product 
$g = u \, \delta$, where $u \in \bU$, $ \delta \in \Delta$. Therefore
any extension $\rho_{\bU_{\Delta}}$ of $\rho$ must satisfy 
$$ \rho_{\, \bU_{\Delta}}(g) = \rho_{\bU} (u) \rho(\delta) \; . $$
We claim that 
this expression also defines a homomorphism from $\bU_{\Delta}$ to
$\bG$.  
We have to verify that the expression is well defined on $\bU_{\Delta}$: 
Assume therefore that  $u \delta = u'  \delta'$. 
Then $ \delta (\delta')^{-1}   = u' u^{-1} \in \bU \cap \Delta$. 
Since $\bU \cap \Delta = \Gamma$ it follows that $\rho_{\bU}(u') \rho_{\bU}(u^{-1}) =  \rho(u'u^{-1}) =
\rho(\delta)   \rho(\delta')^{-1}$. Hence, $\rho_{\bU} (u) \rho(\delta)= \rho_{\bU} (u') \rho(\delta')$.

Similarly, we note that, by Zariski-denseness of $\Gamma$ in $\bU$,  
for all $\delta \in \Delta$, $u \in \bU$, the identitity 
$$\rho_{\bU} (\, \delta  u \delta^{-1} ) = \rho(\delta) \rho_{\bU}(u) \rho( \delta^{-1})$$
holds. This shows that the expression for $\rho_{\, \bU_{\Delta}}$ defines 
a homomorphism. We see  that $\rho_{\, \bU_{\Delta}}$ 
is a $\bbQ$-defined morphism on $\bU_{\Delta}$, by computing $\rho_{\, \bU_{\Delta}}$
on the product of varieties $\bU_{\Delta} = \mu \times \bU$, where $\mu \leq \bU_{\Delta}$ is
a finite subgroup of $\bU_{\Delta}$.  
\end{proof}

In summary,  the algebraic hull for $\Delta$ satisfies the same rigidity 
properties as the Malcev-completion for $\Gamma$ does. The proposition 
implies in particular: 
\begin{corollary} \label{hullrigidity} 
Every automorphism of $\Delta$ extends to a unique $\bbQ$-defined
automorphism of $\bU_{\Delta}$.   
\end{corollary}


\paragraph{Splitting extensions}
We use the algebraic hull of $\Delta$ 
to show that the extension $\Gamma \leq \Delta$ splits in a finite extension.  
Let us first introduce some terminology. 

Let  $\tilde{\Delta}$ be a group, $\tilde{\Gamma} \leq \tilde{\Delta}$ a normal subgroup,
and assume there is a subgroup $\mu \leq  \tilde{\Delta}$ such that 
$\tilde{\Delta}= \tilde{\Gamma} \cdot  \mu$ is a semidirect product. 
Let us put $F= \Delta/\Gamma$.  
The normal extension $\Gamma \leq \Delta$ is said to {\em split\/} (in 
the extension $\tilde{\Gamma} \leq \tilde{\Delta}$)
if there exists  an embedding $\psi: \Delta \longrightarrow \tilde{\Delta}$ 
and  a commutative diagram of homomorphisms
\begin{equation} \label{splitdiagram}   
\xymatrix{
1   \ar[r] & \Gamma  \ar[d]  \ar[r] & \Delta   \ar[d]^{\psi}  \ar[r] & 
F  \ar[d]^{{\rm id}_F} \ar[r] & 1 \\
1   \ar[r] & \tilde{\Gamma}  \ar[r] & \tilde{\Delta}   \ar[r] & F  \ar[r] & 1} \; .
\end{equation}

Now consider $\Delta$ as a subgroup of its algebraic hull $\bU_{\Delta}$. 
We recall that $\bU_{\Delta} = \Delta \bU$, where $\bU$ is the unipotent radical of $\bU_{\Delta}$.
It is known (c.f.\ Theorem \ref{Mostow}) that 
that there exists a finite subgroup $\mu \leq \bU_{\Delta,\bbQ}$ such that
$\bU_{\Delta, \bbQ} = \mu \cdot \bU_{\bbQ}$  is a semi-direct product. 
This implies, in fact,  that the extension $\Gamma \leq \Delta$ splits in
the extension $\bU_{\bbQ} \leq  \bU_{\Delta, \bbQ}$.  

\begin{proposition} \label{splitting}
There exists a finite normal extension  $\Gamma_f \leq \Delta_f$,
where $\Gamma_f$ is a \ftn-group, and an embedding $\Delta \leq_f \Delta_f$,  
such that the extension $\Gamma \leq \Delta$ splits 
in the extension $\Gamma_f \leq \Delta_f$.
\end{proposition}
\begin{proof} 
Let  $\Delta= \Gamma r_1  \cup \, \cdots \, \cup \, \Gamma r_m$  
be a decomposition of $\Delta$ in left cosets. It follows from the splitting
of $\bU_{\Delta, \bbQ}$ that $r_i  = u_i g_i$ with $u_i \in \bU_{\bbQ}$, $g_i \in \mu$, $i=1, \ldots, m$. 
Let $\tilde{\Gamma}$ be the subgroup of $\bU_{\bbQ}$ generated by the $u_i$ and 
$\Gamma$.  Since $u_i \in \bU_{\bbQ}$, $\tilde{\Gamma}$ is a finite extension of $\Gamma$. 
Moreover, since the finite group $\mu$ normalizes $\bU_{\bbQ}$, Lemma \ref{invariance}
implies that there exists a subgroup $\Gamma_f \leq \bU_{\bbQ}$ which
is  a finite extension of $\Gamma$ and which is  normalized by $\mu$. A finitely generated
subgroup of $\bU_{\bbQ}$ is a \ftn-group, hence $\Gamma_f$ is a \ftn-group.
Let $\Delta_f = \Gamma_f \, \mu$ then it follows from the construction that
$\Delta$ is contained in $\Delta_f$. 
\end{proof}

The splitting group $\Delta_f$ inherits the functorial properties of $\bU_{\Delta}$. 
\begin{proposition} \label{fextension}
Let $\rho: \Delta \rightarrow \bG_{\bbQ}$ be a homomorphism of $\Delta$
to a $\bbQ$-defined linear algebraic group $\bG$ such that $\rho(\Gamma)$ is 
contained in $\bH_{\bbQ}$, where $\bH$ is a $\bbQ$-defined  
unipotent subgroup of\/ $\bG$. Then $\rho$ extends to a homomorphism
$\rho_f: \Delta_f \rightarrow \bG_{\bbQ}$.  
\end{proposition}
\begin{proof} By Proposition \ref{hullextension}, $\rho$ extends to a homomorphism
$\rho_{\bU_{\Delta}}: \bU_{\Delta} \rightarrow \bG$. 
Now by construction of $\Delta_f$ we have 
inclusions $\Delta \leq \Delta_f \leq \bU_{\Delta}$, and the restriction
of $\rho_{\bU_{\Delta}}$ to  $\Delta_f$ induces the required extension $\rho_f$. 
\end{proof} 


\subsection{Realization as affine crystallographic group}
We are now dealing with the proof  of 
Theorem \ref{realizationA}. Let $\Delta$, $\Gamma$
satisfy the assumptions of Theorem \ref{realizationA}, i.e.,
$\Gamma$ is an \ftn-group, $\Gamma \leq_f\Delta$ is an effective normal
extension with associated kernel $\beta$.  We assume that 
$\Gamma$ is isomorphic to an ACG of type $A$. \newline 

Let us first consider the special case, where the extension $\Gamma \leq \Delta$
is radicably trivial.  Proposition \ref{extaction} implies that, in this
case, the kernel of the extension acts trivially on $\Def(\Gamma,A)$. 
Also, by Proposition \ref{refreal}, $\Delta$ is a \ftn-group. 
Therefore let us consider now a  (not necessarily normal) finite extension
 $\Gamma \leq_f \Gamma^*$, where $\Gamma^*$ is a \ftn-group.
The following result is an immediate consequence of 
Proposition \ref{Extfunctor}. 

\begin{proposition} \label{nilrealization}
Every $\rho \in\Hom_c(\Gamma,A)$ has an 
extension $\rho_{\, \Gamma^*} \in \Hom_c(\Gamma^*,A)$, and this
extension is unique. 
\end{proposition}

Now let us consider a tower of finite extensions $\Gamma \leq \Gamma^*  \leq \Delta$, 
where the extensions $\Gamma \leq \Delta$, $\Gamma^* \leq \Delta$ are normal. 
For later purposes we  note: 
\begin{lemma} \label{fixlift} 
If\/  $[\rho \, ] \in \Def(\Gamma,A)$ is
a fixed point for $\Delta \rmod \Gamma$
then $[\rho_{ \, \Gamma^*}] \in \Def(\Gamma^*,A)$ is a fixed
point for $\Delta \rmod  \Gamma^*$.   
\end{lemma}
\begin{proof}
Let $c_g: \Gamma^* \rightarrow \Gamma^*$
denote conjugation with $g \in \Delta$. Since $[\rho \,] \in \Def(\Gamma,A)$
is a fixed point for $\Delta\rmod \Gamma$ there exists $a \in A$, such that,
for all $\gamma \in \Gamma$, 
$$\rho_{ \, \Gamma^*}(c_g (\gamma)) =  
\rho \, (c_g(\gamma)) = \rho^a(\gamma) =   \rho_{ \, \Gamma^*\,}^a\!(\gamma) \; .$$ 
By the uniqueness of extensions
we must have $\rho_{ \, \Gamma^*} \circ c_g = \rho_{ \, \Gamma^*\,}^a$.
Therefore $[\rho_{ \, \Gamma^*}]$ is a fixed point for $\Delta \rmod  \Gamma^*$.
\end{proof} 

We turn now to the realization of radicably effective extensions. 
We start by considering split extensions $\Delta= F \cdot \Gamma$,
where $F \leq \Delta$ is a finite subgroup with $F \cap \Gamma= \{1\}$.
Recall from section  \ref{Groupactions}, Proposition \ref{Gstabilizer}, 
that there exists a natural homomorphism   
$c_{\rho}: \N_{A}(\, \rho(\Gamma))  \longrightarrow  \Aut(\Gamma)_{[\rho \, ]}$
which is onto.

\begin{proposition} \label{splitF} 
Let $F \leq \Aut(\Gamma)$ be a finite subgroup.
If\/  $[ \rho \, ] \in \Def(\Gamma,A)$ is a fixed point for $F$, then there exists
an embedding
$$    \iota: \, F  \longrightarrow   \N_{A}(\, \rho(\Gamma)) $$ 
which satisfies  $c_{\rho}(\,  \iota(g)) =g $, for all $g \in F$.  Any two 
embeddings $\iota, \iota': F \rightarrow  \N_{A}(\, \rho(\Gamma))$ with
the latter property  are conjugate by an element of $\C_A(\overline{\rho(\Gamma)})$. 
\end{proposition}
\begin{proof} We will identify $\Aut(\Gamma)$ with a subgroup of $\Aut(U_{\Gamma})$, 
where $U_{\Gamma}$ is a real hull of $\Gamma$, as in section  \ref{Groupactions}. 
From the assumption that $[\rho\,]$ is a fixed point we deduce that  
$F \leq \Aut(\Gamma)_{[\rho \, ]}$. 
We are going to apply some of the facts which are collected in section \ref{Groupactions}.

The group  $U= \overline{\rho(\Gamma)}$ 
is a simply transitive subgroup of $A$. 
By Proposition \ref{Ustabilizer} and Proposition \ref{Gstabilizer},
the map $c: N_A(U) \rightarrow \Aut(U_{\Gamma})_{[\rho]}$
is onto with kernel $\C_A(U)$, and restricts to a surjective map
$c_{\rho}: N_A(\rho(\Gamma)) \rightarrow \Aut({\Gamma})_{[\rho]}$.
Now let $H= c^{-1}(F)$. $H$ is a Zariski-closed subgroup of $N_A(U)$.
Since $F \leq  \Aut({\Gamma})_{[\rho]}$ it follows that $H \leq \N_A(\rho(\Gamma))$.
$H$ contains $\C_A(U)$ as a normal  subgroup of finite index. In fact,
$H \rmod \C_A(U)$ is isomorphic to $F$. Since $\C_A(U)$ is unipotent, $\C_A(U)$ is the 
unipotent radical of $H$. Now by splitting of algebraic groups there
exists a subgroup $\mu \leq H$ such that $H = \mu \cdot\C_A(U)$. 
Since the restriction of $c$ to $\mu$
is an isomorphism onto $F$, we can set $\iota = {c}^{-1}$. 
We thus proved the existence of $\iota$.

We prove now the conjugacy statement.  Let us remark first that
every homomorphism $\iota: F  \rightarrow \N_{A}(\, \rho(\Gamma)) $
which satisfies the assumption maps $F$ into the group $H$. 
Moreover  $\iota$ is
uniquely determined by its image $\iota(F) \leq H$. Since both
$\iota(F)$ and $\iota'(F)$ are Levi-subgroups in $H$, $\iota'(F)$ is
conjugated to $\iota(F)$ by an element of $\C_A(U)$.  
Hence also the homomorphisms $\iota$ and $\iota'$ are conjugate. 
\end{proof}

\begin{corollary} \label{srealization}
Let $F \leq \Aut(\Gamma)$ be a finite subgroup, and
$\Delta = F \cdot \Gamma$ be the corresponding split extension
of $\Gamma$.  Then $\rho \in \Hom_c(\Gamma,A)$ extends 
to a homomorphism $\rho_{\Delta} \in \Hom_c(\Delta,A)$, 
if and only if\/ $[ \rho \, ] \in \Def(\Gamma,A)$ is a fixed point for $F$.  
Any two extensions $\rho_{\Delta}, \rho'_{\Delta} \in \Hom_c(\Delta,A)$ of $\rho$
are conjugate by an element of $ \C_A(\overline{\rho(\Gamma)})$.
\end{corollary}
\begin{proof} By Proposition \ref{splitF}  there exists an embedding 
 $\iota: F \rightarrow \N_{A}(\, \rho(\Gamma))$
such that $c_{\rho}(\,  \iota(g)) =g $, for all $g \in F$.
For $g \in F$, $\gamma \in \Gamma$
we can define $$\rho_{\Delta}(g \,  \gamma) = \iota(g) \rho(\gamma) \; ,$$
to get the required extension $\rho_{\Delta}$. 
Now,  if $\rho_{\Delta}$ is any extension of $\rho$ to $\Delta$, then
$\iota(g) = \rho_{\Delta}(g)$ defines a homomorphism 
$\iota: F \rightarrow  \N_{A}(\, \rho(\Gamma))$ which satisfies
$c_{\rho}(\,  \iota(g)) =g $, for all $g \in F$. Therefore the conjugacy
statement of Proposition \ref{splitF} implies that any two 
extensions of  $\rho$ to $\Delta$ are conjugate by an element 
$u \in \C_A(\overline{\rho(\Gamma)})$. 
\end{proof} 

We are now ready to finish the proof of Theorem \ref{realizationA}. \newline

\noindent \begin{prf}{Proof of Theorem  \ref{realizationA}} 
Let $\rho \in \Hom_c(\Gamma,A)$ be such that $[\rho] \in \Def(\Gamma,A)$
is a fixed point for $\beta$. Let $\Gamma^* = \Fitt(\Delta)$.
By Proposition \ref{nilrealization}, $\rho$ extends to 
a unique homomorphism $\tilde{\rho} \in  \Hom_c(\Gamma^*,A)$.
By Lemma \ref{fixlift},  
$[\tilde{\rho} \, ] \in \Def(\Gamma^*, A)$  
is a fixed point for the kernel of the extension $\Gamma^* \leq \Delta$. 

So we may assume now that $\Gamma = \Fitt(\Delta)$, or equivalently
that the extension $\Gamma \leq \Delta$ is radicably
effective. Let $F = \Delta \rmod \Gamma$. We apply Proposition \ref{splitting} 
to embed $\Delta$ in an extension group $\Delta_f = \mu \cdot \Gamma_f$ 
which splits. 
 $\Gamma$ embeds as a subgroup of
finite index in the \ftn-group $\Gamma_f$, and $\mu$ is isomorphic to $F$. 
By Proposition  \ref{nilrealization}, 
$\rho$ uniquely extends to $\rho_f  \in \Hom_c(\Gamma_f,A)$.
We remark that,  by a similar argument as given in the proof of Lemma \ref{fixlift}, 
$[\rho_f] \in \Def(\Gamma_f, A)$ is a fixed point for $\mu$.
Moreover, since  $\Gamma \leq \Delta$ is radicably
effective, the extension $\Gamma_f \leq  \Delta_f$ is a (radicably)
effective extension. So we may view $\mu$ as a subgroup of $\Aut(\Gamma_f)$.
By Corollary \ref{srealization}, 
the affine crystallographic representation of $\Gamma_f$ extends to 
$\Delta_f$. Since $\Delta$ is a subgroup of $\Delta_f$ we have, in particular,
extended the original representation $\rho$ of $\Gamma$ to
a representation $\rho_{\Delta}$ of $\Delta$. 
\end{prf}

 \paragraph{Appendix: Splitting of algebraic groups} \label{Algsplitting}
In the preceding proofs we made crucial use of the existence of Levi subgroups
in linear algebraic groups and the fact that the Levi subgroups are all
conjugated, a result which is originally due to Mostow. 
Let us recall the precise statement of the result, cf.\ \cite{BS}[5.1]:
\begin{theorem} \label{Mostow} 
Let $\bG$ be a $\bbQ$-defined linear algebraic group. 
Then there exists a reductive $\bbQ$-defined subgroup $\bL$ of $\bG$
such that $\bG= \bL \ur(\bG)$ and $\bL \cap \ur(\bG) = \{ 1 \}$.
Any two such subgroups $\bL$ and $\bL'$ are conjugate by an
element of $ \ur(\bG)_{\bbQ}$.   Moreover, each $\bbQ$-defined 
reductive subgroup $\bH$ of $\bG$ is conjugate to a subgroup of $L$
by an element of $ \ur(\bG)_{\bbQ}$. 
\end{theorem}

From the theorem  corresponding rational splitting  
and conjugacy statements are  easily deduced for any
subfield $k \leq \bbC$. In fact, by the theorem, the variety $\bG$
is a product of varieties $\bL$ and $\ur(\bG)$, where $\bL$ is
isomorphic to $\bG \rmod \ur(\bG)$. The same statement holds 
then for the rational points of these varieties. In particular
we deduce that $(\bG /  \ur(\bG))_{\bbQ} =  \bG_{\bbQ} /  \ur(\bG)_{\bbQ}$,
and $$ \bG_{\bbQ} = \bL_{\bbQ}  \ur(\bG)_{\bbQ} \; \; \mbox{, with }  
\bL_{\bbQ} \cap \ur(\bG)_{\bbQ} = \{ 1 \} \; . $$
The application of this result replaces in our setting the explicit
use of {\em cohomological reasoning\/} which is pervasive in most treatments of 
the classical Bieberbach theory.


\section{Deformation spaces of finite extensions} \label{FExtdef}
Let $\Gamma$ be a \ftn-group which is isomorphic to an ACG of type $A$.
The extension results of the previous section are 
particular well suited to understand the deformation
spaces of the finite extensions of $\Gamma$. 
We continue to use the notational conventions of the previous section. \newline

Let $\beta: F \longrightarrow \Out(\Gamma)$ be an abstract kernel for $\Gamma$
which is realized by a finite extension $\Gamma \leq_f \Delta$.  
The restriction map 
$r_{\Gamma}: \Hom_c(\Delta,A) \longrightarrow \Hom_c(\Gamma,A)$
factorizes to a continuous map  
$$\bar{r}_{\Gamma}:  \; \Def(\Delta,A) \longrightarrow \Def(\Gamma,A) \; $$ 
on the deformation spaces.
We denote the set of fixed points for the kernel $\beta$
on $\Def(\Gamma,A)$ with  $\Def(\Gamma,A)^F$, and we put 
$\Hom_c(\Gamma,A)^F$ for its preimage in $\Hom_c(\Gamma,A)$. 
We topologize  $\Def(\Gamma,A)^F$ 
with the induced subspace topology from $\Def(\Gamma,A)$.
It follows from Proposition \ref{fundamental} that the image of 
$\bar{r}_{\Gamma}$ is contained in $\Def(\Gamma,A)^F$. 
The following 
is then a commutative diagram of (continuous) maps:
$$ \xymatrix{  \Hom_c(\Delta,A) \ar[d]\ar[r]^{r_{\Gamma}} &
 \Hom_c(\Gamma,A)^F   \ar[d] \\
 \Def(\Delta,A)  \ar[r]^ {\bar{r}_{\Gamma}}  & \Def(\Gamma,A)^F 
} \; .
$$

We apply
the results of the previous section to study the properties 
of the map $\bar{r}_{\Gamma}$.
Let  $\Gamma$ be an
\ftn-group and $\bU$ its Malcev completion. Recall from section \ref{Groupactions}
that $\Out(\bU_{\bbR})$ acts naturally on  $\Def(\Gamma,A)$. 
We consider first  the case of commensurable \ftn-groups. 
Their deformation spaces are essentially identical. 
\begin{proposition} \label{commensurable} 
Let $\Gamma$ and $\Gamma'$ be commensurable 
\ftn-groups. Then there exists an $\Out(\bU_{\bbR})$-equivariant
homeomorphism 
$h: \Def(\Gamma,A)\xrightarrow{\; \; {\scriptstyle \approx} \; \;} \Def(\Gamma',A)$. 
\end{proposition} 
\begin{proof} $\Gamma$ and $\Gamma'$ contain a common subgroup of
finite index which is a \ftn-group. It is therefore enough to prove the statement
in the case of a finite extension $\Gamma \leq \Gamma'$. Now by Proposition
\ref{nilrealization} and Proposition \ref{Extfunctor} the restriction map from
$\Gamma'$ to $\Gamma$ is a homeomorphism from $\Hom_c(\Gamma',A)$
to $\Hom_c(\Gamma,A)$. This homeomorphism is easily seen to be compatible with the
actions of $\Aut(\bU_{\bbR})$ on these spaces.  
\end{proof} 

However,  this phenomenon is not restricted to commensurability. More
generally the following holds:
\begin{proposition} Let $\Gamma$ and $\Gamma'$ be lattices in the 
real Lie group $\bU_{\bbR}$. Then there exists an $\Out(\bU_{\bbR})$-equivariant
homeomorphism 
$  h: \Def(\Gamma,A)\xrightarrow{\; \; {\scriptstyle \approx} \; \;} \Def(\Gamma',A)$.
\end{proposition} 
We leave the proof to the reader. In fact, Theorem \ref{hullequivalence} implies that
the properties of the deformation space $\Def(\Delta,A)$ depend only on the
algebraic hull $\bU_{\Delta}$ of $\Delta$.

\begin{corollary} \label{DefF}
Let $\Gamma$ be a \ftn-group and $\Gamma \leq_f \Delta$, where
$\Gamma$ is normal in $\Delta$, be a finite, effective extension, $F = \Delta \rmod \Gamma$.
Then there is a continuous bijection 
$$  \bar{r}_{\Gamma}: 
\Def(\Delta,A) \longrightarrow   \Def(\Gamma,A)^F \; . $$
\end{corollary} 

\begin{remark}
We do  not expect  in general that the map $\bar{r}_{\Gamma}$ will be a homeomorphism. 
Instead,  we expect that there is a stratification of  $\Def(\Delta,A)$ such
that $ \bar{r}_{\Gamma}$ will induce a homeomorphism on the strata.  
\end{remark}
Before we give the proof of the corollary, let us exemplify the previous remark by 
considering a particular important stratum in the deformation space. 
Recall from section \ref{Groupactions}, that $\rho \in \Hom_c(\Gamma,A)$
is called {\em symmetric\/} if the centralizer of $\rho(\Gamma)$ in $A$ 
acts transitively on $V$. Let $\Hom_c(\Gamma,A)_s$ denote the set of all symmetric
crystallographic homomorphisms of $\Gamma$. $\Hom_c(\Gamma,A)_s$
is a closed $A$-invariant subspace of $\Hom_c(\Gamma,A)$, and
the deformation space  $\Def(\Gamma,A)_s$ embeds as a
closed subspace into $\Def(\Gamma,A)$. Let us denote 
$\Def(\Delta,A)_s =  \bar{r}_{\Gamma}^{-1} (\Def(\Gamma,A)_s)$. 

\begin{theorem} \label{Defstratum}
The map  $\bar{r}_{\Gamma}$ induces a homeomorphism
$$  \Def(\Delta,A)_s \xrightarrow{\; \; {\scriptstyle \approx} \; \; \; }   \Def(\Gamma,A)_s^F \; . $$
 \end{theorem} 

\begin{prf}{Proof of Corollary \ref{DefF}}  We want to show that the map 
$\bar{r}_{\Gamma}:  \Def(\Delta,A) \longrightarrow \Def(\Gamma,A)^F$
is a bijection. 
Now $\Delta$ is a subgroup of $\Delta_f$, where $\Delta_f = \Gamma_f \cdot \tilde{F}$
is as in Proposition \ref{splitting}, and it is shown in the proof of  Theorem \ref{realizationA} that
there is  a commutative diagram of restriction maps 
$$ \xymatrix{ 
\Def(\Gamma,A)^F   &  \ar[l]^{\; \; \bar{r}_{\Gamma}} \Def(\Delta,A) \\
\Def(\Gamma_f,A)^{\tilde{F}} \ar[u]^{\cong} &  \ar[l]^{\; \; \bar{r}_{\Gamma_f}}_{\cong}   
\ar[u]  \Def(\Delta_f,A) 
}  \; .
$$
By Proposition \ref{nilrealization} the left upgoing arrow is a bijection. 
The map $ \bar{r}_{\Gamma_f}$ is a bijection by the conjugacy statement of 
Corollary \ref{srealization}, and the right upgoing arrow is onto by 
Proposition  \ref{fextension}. It follows that $\bar{r}_{\Gamma}$ is
a bijection.
\end{prf} 

For the proof of the theorem we need  a lemma.
\begin{lemma} \label{csection}
 Let $\tilde{F} \leq \Aut(\Gamma_f)$ be a finite subgroup, and
$\Delta_f  = \tilde{F} \cdot \Gamma_f$ be the corresponding split extension.  
Then the  
restriction map $r_{\Gamma_f}: \Hom_c(\Delta_f, A) \rightarrow  \Hom_c(\Gamma_f,A)$
admits a (semi-algebraic) continuous  
cross section $s:  \Hom_c(\Gamma_f,A)^F \rightarrow \Hom_c(\Delta_f, A)$.
\end{lemma}
\begin{proof} By Proposition \ref{s_x}, there exists a continuous map 
$s_x:  \Hom_c(\Gamma,A)_s^F \rightarrow \Hom(F, A_x)$
such that $c_{\rho}( s_x (\rho, g) ) = g $, for all $g \in F$. 
Let $\rho \in  \Hom_c(\Gamma,A)_s^F$. 
By Corollary \ref{srealization}, there exists a corresponding 
$\rho_{\Delta_f}=s(\rho) \in \Hom_c(\Delta_f,A)$ which restricts to
$\rho$.  
\end{proof}

\begin{prf}{Proof of Theorem \ref{Defstratum}}
We consider again the diagram in the proof of Corollary \ref{DefF}. 
By Proposition \ref{Extfunctor} the restriction map 
from $\Hom_c(\Gamma_f,A)$ to  $\Hom_c(\Gamma,A)$ is
a homeomorphism, and consequently the induced map 
from $\Def(\Gamma_f,A)_s$ to  $\Def(\Gamma,A)_s$
is a homeomorphism. The same is true for the bijection from
$\Def(\Delta_f,A)$ to $\Def(\Delta,A)$.  Now it follows
from the lemma that $\bar{r}_{\Gamma_f}:  \Def(\Delta_f,A)_s 
\rightarrow  \Def(\Gamma_f,A)$ is a homeomorphism. 
\end{prf} 

\subsection{Some applications} \label{applications}

\paragraph{Moduli spaces of finite extensions}
Let $\Gamma$ be a \ftn-group. Recall that the moduli space $\Mod(\Gamma,A)$ is just
the  space $ \Def(\Gamma,A) \rmod \Out(\Gamma)$
with the quotient topology. Let $\Delta$ be a finite effective extension group
of $\Gamma$. There is a natural induced restriction map 
$$\tilde{r}:  \Mod(\Delta ,A) \longrightarrow \Mod(\Gamma,A)$$
of moduli spaces. In general, this map is not injective and does
not give a precise picture for the space $\Mod(\Delta ,A)$.  

We give now a description of $\Mod(\Delta,A)$ in the case that the
extension is  radicably effective, i.e., we assume that $\Gamma= \Fitt(\Delta)$.
We may view $F= \Delta \rmod \Gamma$ as a subgroup of $\Out(\Gamma)$.
Since $\Gamma$ is an invariant normal subgroup of $\Delta$, there is a natural homomorphism 
$$\bar{\nu}: \, \Out(\Delta) \longrightarrow \N_{\Out(\Gamma)}(F) \rmod F$$
which comes associated with the extension.  Note
that the group $\N_{\Out(\Gamma)}(F) \rmod F$ acts as a group of 
homeomorphisms on $\Def(\Gamma,A)^F$. In fact, this action
desribes the action of $\Out(\Delta)$ on $\Def(\Delta,A)$, and
it is possible to recover $\Mod(\Delta,A)$ from $\Def(\Gamma,A)$
and the action of $F \leq \Out(\Gamma)$.

\begin{theorem} \label{moduli} 
Let $\Delta$ be an effective finite extension of the  \ftn-group $\Gamma$,  and
assume that $\Gamma = \Fitt(\Delta)$. Let $\Out(\Delta)$ act on $\Def(\Delta,A)$
with respect to the homomorphism $\nu$. 
Then the  embedding $\bar{r}_{\Gamma}:  \Def(\Delta,A) \longrightarrow  \Def(\Gamma,A)^F$
is $\Out(\Delta)$-equivariant. In particular, there is a continuous bijection from the moduli space 
$\Mod(\Delta,A)$ to the quotient space
$\Def(\Gamma,A)^F \rmod \bar{\nu}(\Out(\Delta))$. The bijection induces 
a homeomorphism on the moduli space $\Mod(\Delta,A)_s$ 
of symmetric crystallographic homomorphisms to $\Def(\Gamma,A)_s^F \rmod \bar{\nu}(\Out(\Delta))$.
\end{theorem} 
\begin{proof} That the map $\bar{r}_{\Gamma}$ is $\Out(\Delta)$-equivariant is a routine
calculation from the definitions. We omit it therefore. 
The remaining statements follow then from the
previous results on deformation spaces. 
\end{proof} 

\paragraph{Inheritance of the Hausdorff-property} 
Though the topology on $\Def(\Delta,A)$ is potentially
finer than the topology induced from $\Def(\Fitt(\Delta),A)^F$
we still can deduce an immediate useful consequence from 
Corollary \ref{DefF}.
\begin{corollary} 
Let $\Delta$ be a virtually nilpotent ACG, and let 
$\Gamma = \Fitt(\Delta)$ be the Fitting subgroup. 
If\/ $\Def(\Gamma,A)$ is a Hausdorff topological
space, then $\Def(\Delta,A)$ is a Hausdorff space too. 
\end{corollary} 

\paragraph{Inheritance of convexity to finite extensions}
If $\Delta$ is a virtually nilpotent ACG then the 
Fitting subgroup $\Gamma = \Fitt(\Delta)$ is a 
crystallographic \ftn-group,  and is also a characteristic
subgroup of $\Delta$. Therefore there is a natural
restriction homomorphism 
$$\nu: \Aut(\Delta) \rightarrow \Aut(\Gamma) \; ,$$
and, evidently, the map 
$r_{\Gamma}: \Hom_c(\Delta,A) \rightarrow \Hom_c(\Gamma,A)$
is equivariant with respect to $\nu$. It follows 

\begin{proposition} Let $\Delta$ be virtually nilpotent ACG, and let $\Gamma = \Fitt(\Delta)$.
Then the induced inclusion  on deformation spaces 
$$\Def(\Delta,A) \longrightarrow   \Def(\Gamma,A) $$ 
is equivariant with respect to the actions of $\Aut(\Delta)$ and $\Aut(\Gamma)$. 
\end{proposition}

As a particular consequence we see that convexity properties
of $\Def(\Fitt(\Delta),A)$ will be  inherited to $\Delta$. 

\begin{theorem} \label{inheritance1} 
Let $\Delta$ be a virtually nilpotent ACG, and let $\Gamma = \Fitt(\Delta)$.
Then the following hold:
\begin{itemize}
\item[i)] If $\Def(\Gamma,A)$ is fixed pointed then $\Def(\Delta,A)$ 
is fixed pointed too. 
\item[ii)] If $\Def(\Gamma,A)$ is convex then $\Def(\Delta,A)$ 
is convex too. 
\end{itemize}
\end{theorem}
\begin{proof} Let $F = \Delta \rmod \Gamma$.  
If  $\Def(\Gamma,A)$ has a fixed point $[\rho_o]$ for
$\Aut(\Gamma)$ then,  by  Corollary \ref{DefF},
$\bar{r}_{\Gamma}(\Def(\Delta,A))$ meets $[\rho_o]$.
By equivariance and the injectivity of the map $\bar{r}_{\Gamma}$, 
$[\rho_o]$ is a fixed point for $\Aut(\Delta)$ on $\Def(\Delta,A)$.

We assume now that $\Def(\Gamma,A)$ is convex. The restriction
homomorphism $\nu$ induces a map
$$ \bar{\nu}: \Out(\Delta) \longrightarrow \N_{\Out(\Gamma)}(F) \rmod F \; .$$ 
Let $\mu$ be
a finite subgroup of $\Out(\Delta)$, and let $\tilde{\mu} \leq N_{\Out(\Gamma)}(F)$ 
be the preimage of $\bar{\nu}(\mu)$. 
Since $F$ is finite so is $\tilde{\mu}$ which is a finite normal
extension of $F$. If $[\rho_o] \in \Def(\Gamma,A)$ is a fixed point for $\tilde{\mu}$
then, since $F \leq \tilde{\mu}$, $\bar{r}_{\Gamma}(\Def(\Delta,A))$ meets $[\rho_o]$. 
By the equivariance properties of the embedding of $\Def(\Delta,A)$ into $\Def(\Gamma,A)$, 
$[\rho_o]$ is a fixed point for $\mu$ too. 
\end{proof}

\paragraph{The realization theorem} 
We now come to the proof of the realization theorem 
for  finite extensions of a virtually nilpotent affine crystallographic group:
\newline

\begin{prf}{Proof of Theorem \ref{Affinerealization1}} 
Since the virtually nilpotent group $\Lambda$ is isomorphic to an affine 
crystallographic group, $\Theta= \Fitt(\Lambda)$ is
a \ftn-group, and the centralizer of  $\Theta$ is 
contained in $\Theta$. Hence the extension $\Theta \leq \Lambda$ 
is effective. It follows from Lemma \ref{eclosed} 
that $\Delta$ is an effective extension of $\Gamma = \Fitt(\Delta)$,
and also that $\Gamma= \Fitt(\Delta)$ is a \ftn-group. 
Since we already proved the realization theorem 
for effective extensions of \ftn-groups, to show that 
$\Delta$ is isomorphic to a crystallographic group, it is enough to
show that the kernel $F \leq \Out(\Gamma)$ which is associated to the 
extension $\Gamma \leq \Delta$ has a fixed point in $\Def(\Gamma,A)$. 

We remark that $\Theta = \Fitt(\Lambda)$  has 
finite index in $\Gamma = \Fitt(\Delta)$. 
Therefore (compare Proposition \ref{commensurable}) $\Out(\Theta)$ 
and the kernel $F_{\Theta} \leq \Out(\Theta)$ which
is associated to $\Theta \leq \Lambda$
act also  on $\Def(\Gamma,A)$, and the set of fixed points of $F_{\Theta}$ 
corresponds to $\Def(\Lambda, A)$. The action of $\Out(\Theta)$ on 
$\Def(\Gamma,A)$ factorizes over $\Out(\bU_{\bbR})$. The action of
$\Out(\Lambda)$ on  $\Def(\Lambda, A) \subseteq \Def(\Gamma,A)$ is then 
encoded in the homomorphism
$$ \bar{\nu}: \Out(\Lambda) \longrightarrow \N_{\Out(\bU_{\bbR})}(F_{\Theta}) \rmod F_{\Theta} \; .$$
Let $\mu \leq \Out(\Lambda)$ be the kernel
which is associated to the extension $\Lambda \leq \Delta$, and
let $\tilde{\mu} \leq \Out(\bU_{\bbR})$ be the 
preimage of $ \bar{\nu}(\mu)$
in $\N_{\Out(\bU_{\bbR})}(F_{\Theta})$.   
By our assumption,  $\mu$ has a fixed point in $\Def(\Lambda, A)$ 
and this implies that $\tilde{\mu}$ has a fixed point in $\Def(\Gamma,A)$. 
It is easy to see
that the kernel $F \leq \Out(\bU_{\bbR})$ 
is contained in $\tilde{\mu}$. Therefore $F$ has a
fixed point in $\Def(\Gamma,A)$. 
\end{prf}

\subsection{Manifolds of Euclidean type}  \label{applications1} 
Let $M$ be a closed manifold which admits a Riemannian metric of
constant curvature zero, i.e., a flat Riemannian metric. In the spirit
of the theory of geometrization,  $M$ is called a
{\em manifold of Euclidean type}. By Bieberbach
theory $M$ is finitely covered by a torus $T^n$. The fundamental
group $\pi_1(M)$ is isomorphic to a finite extension of $\bbZ^n$. It was remarked
previously in \cite{BauesV} that the deformation space 
$\Def(T^n, \Aff(\bbR^n))$ of affine strucutures on
$T^n$ is a Hausdorff-space and homeomorphic to a semi-algebraic
set. We obtain the following generalization of this result: 

\begin{corollary} Let $M$ be a closed manifold of Euclidean type. 
Then the deformation space of complete affine structures of type
$A$ is a Hausdorff space and homeomorphic to a semi-algebraic
set.
\end{corollary}
\begin{proof} First let us remark that, 
by Theorem \ref{tsemialg}, $\Def(\bbZ^n, A)$ is a
semi-algebraic set. 
Now $\pi_1(M)$ is a normal extension of $\bbZ^n$ with some
finite group $F$. We remark next 
that the crystallographic actions of $\bbZ^n$ are
all symmetric. 
By Corollary \ref{Defstratum} it follows that 
$\Def(M, A)$ is homeomorphic to $\Def(\bbZ^n,A)^F$. 
Therefore $\Def(M, A)$ is homeomorphic to a closed subspace 
of $\Def(\bbZ^n, A)$, and also it is Hausdorff. Since the action  
of $\GL(n,\bbZ)$ on $\Def(\bbZ^n,A)$ is algebraic $\Def(\bbZ^n,A)^F$
is a semi-algebraic set. 
\end{proof} 

\paragraph{The Burckhardt-Zassenhaus theorem}
We let $\GL_A$ act on $\GL(V)$ by left multiplication and put
$$  X_A \; = \;    {_{\displaystyle \GL_A}}\! \lmod \, \GL(V)    \; .$$
Let $F= {G}_{st}(V,A)$ be the algebraic 
variety of simply transitive abelian subgroups of $A$.
$\GL(A)$ acts by conjugation on $F$. 
It was proved in Proposition \ref{fproduct} that 
$$  \Def(\bbZ^n, A) \; = \; \GL(V) \underset{\GL_A}{\times}  F$$
is a fiber product, and in particular a bundle over the homogeneous
space $X_A$ with fiber $F$. The $\GL(V)$-action on $\Def(\bbZ^n, A)$
corresponds to the natural $\GL(V)$-action on the fiber product which is
induced by right multiplication of  $\GL(V)$ on itself. 

If $A$ contains the
subgroup $V$ of translations we note that $\GL_A$
fixes $V$ as an element of  $G_{st}(V,A)$.
Therefore, provided $A$ contains the
subgroup $V$ of translations, $\mu \leq \GL(V)$ fixes 
a point in $\Def(\bbZ^n, A)$ if and only if it
fixes a point on $X_A$.  
It follows that, in this case,  the  existence problem for fixed points
in the deformation spaces of tori is only a matter of
the action of $\GL(V)$ on $X_A$.
In particular we recover from Theorem \ref{realizationA} 
the Burckhardt-Zassenhaus theorem
(see  \cite{Zassenhaus}) on Euclidean crystallographic
groups:
\begin{corollary} 
Every finite effective extension of $\bbZ^n$ acts as
an Euclidean crystallographic  group on affine space 
$\bbA^n$. 
\end{corollary}
\begin{proof} Every finite subgroup $\mu$ of $\GL_n$ has
a fixed point in $X = O_n \lmod \GL_n$, as is well
known. Since ${G}_{st}(V, \E(n))$ contains the
subgroup of translations, $\mu$ has a fixed point  
on $\Def(\bbZ^n, E(n))$. 
\end{proof}

In general,  it seems to be a nontrivial problem to give
a precise description of the variety  ${G}_{st}(V,A)$, 
and therefore also of the deformation spaces for
complete affine structures on manifolds of Euclidean type. 
Next we cover  some accessible special cases.  

\paragraph{Complete affine surfaces}
It is well known that a closed affine surface is
diffeomorphic to the two-torus or the Klein bottle. 
In  \cite{BauesG} the following result is proved:
\begin{proposition}The deformation
space $\Def(T^2, \Aff(\bbR^2))$ of the two-torus is homeomorphic to $\bbR^2$.
In these coordinates the action of $\GL_2(\bbZ)$ 
corresponds to the natural linear action. The fixed point 
$0 \in \bbR^2$ represents the flat Euclidean structure on $T^2$.
\end{proposition}  

It is easy to see from \cite{BauesG} that also the 
natural $\GL_2(\bbR)$-action
is linear in the above model for  $\Def(T^2, \Aff(\bbR^2))$. 
In a sense, our model encodes all information 
on the affine crystallographic groups in dimension two.
We apply our method to study the two dimensional affine
crystallographic groups up to affine conjugacy: 
Recall that there exist up to 
isomorphism seventeen finite effective extensions
of $\bbZ^2$. All 
these groups act as Euclidean crystallographic groups on the plane. 
These are the famous {\em wall paper groups}.
As crystallographic subgroups of the Euclidean group, the wall paper groups are unique
up to conjugacy with affine transformations. 
Let us call an affine crystallographic group
{\em non-Euclidean\/} if it is not conjugated to an Euclidean crystallographic
group by an affine transformation.
So which of the finite effective extensions of $\bbZ^2$ admit non-Euclidean affine
crystallographic actions? 

\begin{proposition} Let\/  $\Gamma$ be a finite effective extension of $\bbZ^2$ which 
acts as a non-Euclidean affine crystallographic group.  If\/  $\Gamma$ preserves
orientation then $\Gamma$ is
isomorphic to  $\bbZ^2$. 
If  $\Gamma$ does not preserve orientation then $\Gamma$ contains $\bbZ^2$ as 
a translation subgroup of index two. There are three isomorphism classes of
the latter type and for these groups the deformation space
$\Def(\Gamma, \Aff(\bbR^2))$ is homeomorphic to the 
real line.
\end{proposition}
\begin{proof} Let $\mu$ be a finite subgroup of $\GL_2(\bbZ)$ which
fixes a line in $\bbR^2$. Every nonidentity element of $\mu$ has
eigenvalues $1$ and $-1$ and is of order two. It follows that $\mu$
is conjugate in  $\GL_2(\bbZ)$ to 
$$
     \mu_1= \;  < \begin{matrix}{cr} 1 & 0  \\ 0 & -1\end{matrix}> \: \mbox{ or }\; 
     \mu_2  = \; < \begin{matrix}{cc} 0 & 1 \\ 1 & 0 \end{matrix} > \;  \; . 
$$
From the classification (see for example \cite{Zieschang}[\S 22, \S 23])
it follows that there are two (effective) extensions of $\bbZ^2$
by $\mu_1$, and only one by  $\mu_2$.  The torsionfree extension of $\bbZ^2$ 
by $\mu_1$ is the fundamental group of the Klein bottle.
\end{proof} 

\begin{corollary} The deformation space $\Def({\cal K}, \Aff(\bbR^2))$
of complete affine structures on
the Klein bottle ${\cal K}$ is homeomorphic to the real line. The moduli space
$\Mod({\cal K}, \Aff(\bbR^2))$ coincides with the deformation
space. 
\end{corollary}
\begin{proof} 
Let $\Delta$ be the Klein-bottle group, $q \in \mu_1$ the non-identity
element of $\mu_1$, $a,b$ the standard generators for the lattice $\bbZ^2$.
Then $\Delta$ is generated by $q,a,b$ with the relations $qaq^{-1} = a$,
$ q b q^{-1} = -b$ and $q^2 =a$.   
The group $\mu_1$ fixes a line in $\bbR^2$ 
which is the deformation space for $\Delta$.
It is easy to see 
that the image of the map from $\Out(\Delta)$ to the
normalizer of $\mu_1$ in $\GL_2(\bbZ)$ is $\mu_1$ itself.
By Theorem  \ref{moduli} the fixed line in $\bbR^2$ is
therefore also the moduli space.   
\end{proof} 

\paragraph{Metric and symplectic affine structures on tori}
We turn our attention to some particular interesting subgeometries,
namely complete affine structures with a parallel inner product 
(flat pseudo-Riemannian manifolds) or a parallel nondegenerate 
alternating product (symplectic affine manifolds). In these 
cases we can give a precise description of the 
spaces ${G}_{st}(V,A)$, see section \ref{atori}.   
This amounts to a good understanding of the deformation spaces for these 
geometries on manifolds of Euclidean type. 
\newline 



\section{The realization problem for unipotent shadows} \label{srealizations}
Every torsionfree polycyclic group $\Gamma$ has an unipotent shadow $\Theta$ 
in its real algebraic hull $H_{\Gamma}$. 
Recall from Proposition \ref{cshadowdescent}, that if $\Gamma$ is an affine
crystallographic group then
the unipotent shadow $\Theta$ is an affine crystallographic group
as well. 
Therefore the following {\em realization problem\/} makes sense:  \newline

\noindent {\em Let $\Gamma$ be a torsionfree polycyclic group with
unipotent shadow $\Theta$. Assume that $\Theta$
is isomorphic to an affine crystallographic group. Is it true
that $\Gamma$ is isomorphic to an affine crystallographic group? } \newline

We need the following
\begin{definition} Let $H_{\Gamma}$ be a real algebraic hull for $\Gamma$. We choose a splitting
$H_{\Gamma} = T \cdot U_{\Gamma}$, where $T$ is a Levi subgroup of $H_{\Gamma}$
and put $U_{\Gamma} = \ur(H_{\Gamma})$. Then $\Theta \leq U_{\Gamma}$, and 
$U_{\Gamma}$ is a real Malcev hull for $\Theta$. Since $H_{\Gamma}$ has a strong unipotent
radical, conjugation embeds $T$ as a subgroup of $\Aut(U_{\Gamma})$.
We call the corresponding image $T_{\Gamma} \leq \Out(U_{\Gamma})$  of $T$ the 
{\em semisimple kernel associated to $\Gamma$}.  
\end{definition}
Recall (see 
section \ref{Groupactions}) that, since $U_{\Gamma}$ 
is a Malcev-hull for $\Theta$,
$\Out(U_{\Gamma})$, and hence also $T_{\Gamma}$, 
act on the deformation space  $\Def(\Theta,A)$
of $\Theta$.
The main result of ths section is the following theorem. 
It shows that the answer to the question above lies in 
the action of the semisimple kernel $T_{\Gamma}$ of $\Gamma$ on the 
deformation space of $\Theta$. 
\begin{theorem} \label{shadowrealization} 
Let $\Gamma$ be a torsionfree polycyclic group with
unipotent shadow $\Theta$. Assume that $\Theta$ is an ACG and
let $\tau \in \Hom_c(\Theta,A)$. If $[\tau] \in \Def(\Theta,A)$ is
a fixed point for the semisimple kernel $T_{\Gamma}$ of\/ $\Gamma$
then there exist $\rho \in \Hom_c(\Gamma, A)$ which descends 
to $\tau$. 
\end{theorem}  

An immediate consequence is the following realization
result. It gives an answer to the question which polycyclic
groups can act as affine crystallographic groups in terms 
of the set of crystallographic actions of the  unipotent shadow. 
\begin{corollary} Let $\Gamma$ be a torsionfree polycyclic group, 
$\Theta \leq H_{\Gamma}$ a unipotent shadow for $\Gamma$.  
Then $\Gamma$ may be realized as an affine crystallographic group of type $A$
if and only if there exists a fixed point for the action of the  semisimple kernel 
$T_{\Gamma}$ in the deformation space $\Def(\Theta, A)$. 
\end{corollary}
\begin{remark} 
Auslander (compare \cite{Auslander})
remarked that the unipotent shadow of a crystallographic
group is crystallographic. 
The realization criterion of the 
Corollary  may be seen as providing a converse to his result.    
\end{remark}
Applications to the existence problem for crystallographic actions on
large classes of polycyclic groups follow. (See section \ref{existence})

\subsection{The shadow map on deformation spaces}
Let $\Gamma$ be a torsionfree polycyclic group, $\Gamma \leq \bH_{\Gamma}$ and
$\Theta = \Gamma^u_C \leq \ur(\bH_{\Gamma})$ a unipotent shadow 
for $\Gamma$ with respect
to some (almost) nilpotent supplement $C$ for $\Fitt(\Gamma)$. Recall
from Proposition \ref{cshadowdescent}, that there exists a canonical  crystallographic
shadow map 
$$  s_u: \Hom_c(\Gamma, A) \longrightarrow \Hom_c(\Theta,A) \; . $$
We remark first (see Proposition \ref{inducedshadowmap}) that the crystallographic
shadow map induces a shadow map 
$$ \sigma_u: \Def(\Gamma, A) \longrightarrow \Def(\Theta,A)  $$
on the deformation spaces.  We study here the properties of the shadow map $\sigma_u$
in order to obtain information about the deformation spaces $\Def(\Gamma, A)$. 
The results and methods are analoguous to
those in section \ref{realproblems} on the realization of finite extension for
\ftn-groups. We can interpret them also as a solution to a suitable 
realization problem for polycyclic groups with a prescribed shadow.

\paragraph{The induced shadow map}
Let $H_{\Gamma}$ be a real algebraic hull for $\Gamma$, and let
$\Theta$ be a unipotent shadow. We choose a splitting
$H_{\Gamma} = T \cdot U_{\Gamma}$ of $H_{\Gamma}$.
We may view $T$ as a subgroup of $\Aut(U_{\Gamma})$, and $T$ acts
on  $\Def(\Theta,A)$ via the semisimple kernel $T_{\Gamma} \leq \Out(U_{\Gamma})$.

\begin{proposition} \label{inducedshadowmap}
The crystallographic shadow map  $s_u$ induces a continuous
map $$ \sigma_u:  \Def(\Gamma, A) \longrightarrow \Def(\Theta,A)^T  \; .
$$ 
\end{proposition}
\begin{proof} To see that 
there is an induced map $\sigma_u$
on deformation spaces which is given by 
$ \sigma_u([\rho \, ])  \; =  \; [ s_u(\rho) \,]$,
it is enough to verify that,  for all $g \in A$,  $\rho \in \Hom_c(\Gamma, A)$,
$$      s_u (\,  \rho^g \, ) \; = \;  s_u(\rho)^g  \; .  $$
To verify this,  we have to recall the construction of the
homomorphism $s_u(\rho) = \rho^u$. If $\Gamma = C \,  \Fitt(\Gamma)$
then, by Proposition \ref{shadowdescent},  
$\rho^u: \Gamma^u_C \rightarrow A$ is determined
by the conditions $\rho^u(\gamma_u) = \rho(\gamma)_u$,  for
all $\gamma \in C$, $\gamma \in \Fitt(\Gamma)$. Therefore the
above formula is immediate since the Jordan-decomposition in $A$
is preserved by conjugation. The map $\sigma_u$ is continuous 
since (Corollary \ref{shadowmapisc}) $s_u$ is continuous.   

Next we have to prove that  $\sigma_u(\rho)  \in  \Def(\Theta,A)^T$. 
But this is an immediate consequence of Theorem \ref{acextension}:
The homomorphism $\rho: \Gamma \rightarrow A$ 
extends to a $u$-simply transitive embedding
$\rho_{H_{\Gamma}}: H_{\Gamma} \rightarrow A$.  By Proposition 
\ref{shadowdescent}, $\rho_{H_{\Gamma}}$ restricts to
$s_u(\rho)= \rho^u$ on $\Theta = \Gamma^u_C \leq U_{\Gamma}$. 
Then  the action of 
$T \leq \Aut(U_{\Gamma})$ on $\Hom_c(\Theta, A)$
is induced by conjugation via $\rho_{H_{\Gamma}}(T)$.
Hence,  $T \leq \Aut(U_{\Gamma})_{\sigma_u([\rho]) }$. 
\end{proof}

We are going to prove here

\begin{theorem} \label{shadowembedding}
Let $\Gamma$ be a torsionfree polycyclic group,
$\Theta \leq U_{\Gamma}$ a unipotent shadow for $\Gamma$,
and $T$ the semisimple kernel associated to $\Gamma$. 
Then the  induced shadow map 
$$ \sigma_u: \Def(\Gamma, A) \longrightarrow \Def(\Theta,A)^T$$
is a continuous bijection onto $\Def(\Theta,A)^T$. 
\end{theorem}

It does not seem to be clear whether the map $\sigma_u$ is 
also a homeomorphism onto its image. However, this is the case
on a certain stratum of $\Def(\Gamma, A)$. Let 
$\Def(\Theta, A)_s$ be the closed
stratum of symmetric structures in the deformation space
$\Def(\Theta, A)$. (Compare Theorem \ref{Defstratum}). 
We define a stratum $\Def(\Gamma, A)_s = \sigma_u^{-1}(\Def(\Theta,A)_s)$ 
in $\Def(\Gamma, A)$. 

\begin{theorem} \label{shadowstratum}
The map  $\sigma_u$ induces a homeomorphism
$$  \Def(\Gamma,A)_s \xrightarrow{\; \; {\scriptstyle \approx} \; \; \; }   
     \Def(\Theta,A)_s^T \; . $$ 
\end{theorem}

We  postpone the proof of Theorem
\ref{shadowrealization}, Theorem \ref{shadowembedding}, and Theorem \ref{shadowstratum}  
for a moment in order to derive some consequences.  

\paragraph{Inheritance of the Hausdorff-property} 
Though the topology on $\Def(\Gamma,A)$ is potentially
finer than the topology induced from $\Def(\Theta,A)_s^T$
we still can deduce an immediate useful consequence.
\begin{corollary} 
Let $\Gamma$ be a torsionfree polycyclic ACG, and let $\Theta$ be a unipotent
shadow for $\Gamma$. If\/ $\Def(\Theta,A)$ is a Hausdorff topological
space, then $\Def(\Gamma,A)$ is Hausdorff too. 
\end{corollary} 
Some applications of this fact were already described 
in chapter \ref{defspaces}. If $\Theta$ is abelian 
then $\Def(\Theta, A) = \Def(\Theta, A)_s$. Therefore it 
follows from Theorem \ref{shadowstratum}: 
\begin{corollary} 
Let $\Gamma$ be a torsionfree polycyclic ACG, so that the unipotent
shadow  $\Theta$ is abelian. Then $\Def(\Gamma,A)$ is 
homeomorphic to a semi-algebraic set. 
\end{corollary}

\paragraph{Inheritance of convexity}
Let us recall that every automorphism $\phi \in \Aut(\Gamma)$ extends
to a unique algebraic automorphism $\Phi$ of $H_{\Gamma}$. Since
$\Phi$ is algebraic it restricts to an automorphism $\Phi^u$ of
the unipotent radical $U_{\Gamma}$ of $H_{\Gamma}$.  
Therefore there is a natural
induced homomorphism 
$$ \nu^u: \Aut(\Gamma) \longrightarrow \Aut(U_{\Gamma}) \; . $$
Since $U_{\Gamma} = U_{\Theta}$ is a Malcev-hull for $\Theta$,
the group $\Aut(U_{\Gamma})$ acts on the Deformation space 
$\Def(\Theta,A)$ and it is easy to see that the shadow map 
$s_u: \Hom_c(\Gamma,A) \rightarrow \Hom_c(\Theta,A)$
is equivariant with respect to $\nu^u$. We summarize:  

\begin{proposition} Let $\Gamma$ be a polycyclic ACG, $\Theta$ a unipotent
shadow for $\Gamma$.
Then the shadow map on  deformation spaces 
$$\sigma_u: \Def(\Gamma,A) \hookrightarrow   \Def(\Theta,A) $$ 
is equivariant with respect to the actions of $\Aut(\Gamma)$ and $\Aut(U_{\Gamma})$
and the homomorphism $\nu^u$.  
\end{proposition}

As a particular consequence we see that the convexity properties
of the unipotent shadow $\Theta$ will be  inherited.

\begin{theorem} \label{inheritance2}  
Let $\Gamma$ be a polycyclic ACG and $\Theta$ a unipotent
shadow for $\Gamma$.
Then the following hold:
\begin{itemize}
\item[i)] If $\Def(\Theta,A)$ is fixed pointed then $\Def(\Gamma,A)$ 
is fixed pointed too. 
\item[ii)] If $\Def(\Theta,A)$ is convex then $\Def(\Gamma,A)$ 
is convex too. 
\end{itemize}
\end{theorem}
\begin{proof}   
If  $\Def(\Theta,A)$ has a fixed point $[\rho_o]$ for
$\Aut(U_{\Gamma})$, Theorem \ref{shadowembedding}
implies that $\sigma_u(\Def(\Gamma,A))$ meets $[\rho_o]$.
By equivariance and injectivity of the map $\sigma_u$, 
$[\rho_o]$ is a fixed point for $\Aut(\Gamma)$ in $\Def(\Theta,A)$.

We assume now that $\Def(\Theta,A)$ is convex. 
Let $T_{\Gamma} \leq \Out(U_{\Gamma})$ be
the semisimple kernel for the shadow $\Theta$.
It is easy
to see that there is a well defined homomorphism 
$$ \bar{\nu}^u: \Out(\Delta) \longrightarrow \N_{\Out(U_{\Gamma})}(T_{\Gamma}) \rmod T_{\Gamma} \; . $$
Let $\mu \leq \Aut(\Gamma)$ be a finite
subgroup, and let $\tilde{\mu} \leq N_{\Out(\Gamma)}(T_{\Gamma})$ 
be the preimage of $\bar{\nu}^u(\mu)$.
Since $T_{\Gamma}$ is reductive so is $\tilde{\mu}$ which is a finite normal
extension of $T_{\Gamma}$. If $[\rho_o]$ is a fixed point for $\tilde{\mu}$
then, since $T_{\Gamma} \leq \tilde{\mu}$, $\sigma_u(\Def(\Gamma,A))$ meets $[\rho_o]$.
Moreover $[\rho_o]$ is a fixed point for $\mu$ too. Therefore  $\Def(\Gamma,A)$ is convex. 
\end{proof}

\paragraph{Surjectivity of the shadow map} 
We start now with the proofs of the main theorems.
Let $\Theta$ be a \ftn-group, $U_{\Theta}$ a real Malcev-hull for $\Theta$. 
For $\tau \in \Hom_c(\Theta,A)$ we
put $U= \ac{\tau(\Theta)}$ to denote the unipotent simply transitive hull
of $\tau(\Theta) \leq A$. Recall from Proposition \ref{Ustabilizer} that in this
situation conjugation in $A$ defines a natural map 
$$ c_{\tau}:  N_A(U) \longrightarrow \Aut( U_{\Theta} )_{[\tau]} \; . $$
The map $c_{\tau}$ is then onto with kernel $\C_A(U)$.

\begin{proposition} \label{splitT} 
Let $T \leq \Aut(U_{\Theta})$ be a reductive algebraic subgroup.
If\/  $[ \tau \, ] \in \Def(\Theta,A)$ is a fixed point for $T$, then there exists
an  embedding of (real-) linear algebraic groups  
$$    \jmath: \, T  \longrightarrow   \N_{A}(  U ) $$ 
which satisfies  $c(\,  \jmath(t)) =t $, for all $t \in T$.  Any two 
embeddings $\jmath, \jmath': T \rightarrow  \N_{A}(U)$ with
the latter property  are conjugate by an element of $\C_A(U)$. 
\end{proposition}
\begin{proof} 
By assumption $T  \leq  \Aut(U_{\Theta})_{[\tau]}$. 
Let $H= c^{-1}(T)$. $H$ is a Zariski-closed subgroup of $N_A(U)$, and
$H$ contains the unipotent group $\C_A(U)$ as a normal  subgroup.
In fact, since 
$H \rmod \C_A(U)$ is isomorphic to $T$, 
$\C_A(U)$ is the unipotent radical of $H$. Now by splitting of algebraic groups there
exists a subgroup $T_H \leq H$ such that $H = T_H \cdot\C_A(U)$. 
Since the restriction of $c$ to $T_H$
is an isomorphism of algebraic groups onto $T$, 
we can set $\jmath = {c}^{-1}: T \rightarrow T_H \leq H$. 
We thus proved the existence of $\jmath$.

We prove now the conjugacy statement.  Let us remark first that
every homomorphism $\jmath: T \rightarrow \N_{A}(U)) $
which satisfies the assumptions maps $T$ into the group $H$. 
Moreover  $\jmath$ is
uniquely determined by its image $\jmath(T) \leq H$. Since $\jmath(T)$
and $\jmath'(T)$ are Levi-subgroups in $H$, $\jmath'(T)$  is
conjugated to $\jmath(T)$ by an element of $\C_A(U)$.  
Hence also the homomorphisms $\jmath$ and $\jmath'$ are conjugate. 
\end{proof}

Let $\Gamma$ be a torsionfree polycyclic group with real 
algebraic hull $H_{\Gamma}$, $\Theta \leq U_{\Gamma}$
a unipotent shadow for $\Gamma$. 
We say that $\rho_{H_{\Gamma}}: H_{\Gamma}  \rightarrow  A$
extends a homomorphism $\tau \in \Hom_c(\Theta, A)$ if $\tau$ is the restriction of
$\rho_{H_{\Gamma}}$ to $\Theta$.

\begin{corollary} \label{hullrealization}
A homomorphism $\tau \in \Hom_c(\Theta,A)$ extends 
to a $u$-simply transitive embedding of linear algebraic groups
$\rho_{H_{\Gamma}}: H_{\Gamma}  \rightarrow  A$ 
if and only if\/ $[ \tau \, ] \in \Def(\Theta,A)$ is a fixed point 
for the semisimple kernel $T_{\Gamma}$.  
Any two such extensions $\rho_{H_{\Gamma}}, \rho'_{H_{\Gamma}}$
of $\tau$ are conjugate by an element of $\C_A( \ac{\tau(\Theta)} )$.
\end{corollary}
\begin{proof} Since $H_{\Gamma}= T \cdot U_{\Gamma}$ has a strong unipotent
radical, we may view the reductive group $T$ as a subgroup of 
$\Aut(U_{\Gamma})$. Let us assume that $[ \tau \, ] \in \Def(\Theta,A)$ is a fixed point 
for $T_{\Gamma}$.
By Proposition \ref{splitT}, there exists an embedding 
$\jmath: T \rightarrow \N_{A}(\,U)$, $U= \ac{\tau(\Theta)} \leq A$, 
such that $c(\,  \jmath(t)) =t $, for all $t \in T$.
We already used implicitly that, since $U_{\Gamma}$ is a real
Malcev hull for $\Theta$,  the homomorphism $\tau \in \Hom_c(\Theta,A)$
extends to a $u$-simply transitive homomorphism $\tau_{U_{\Gamma}}: U_{\Gamma}
\rightarrow A$. 
For $t \in T$, $u  \in U_{\Gamma}$
we can then define $$\rho_{H_{\Gamma}}(t \, u)  \; =  \; \jmath(t)  \, 
\tau_{U_{\Gamma}}(u) \; .$$
By the properties of $\jmath$, this defines a homomorphism 
$\rho_{H_{\Gamma}}: H_{\Gamma} \rightarrow A$ of algebraic groups
which is clearly injective and $u$-simply transitive. 
Thus we proved the ``if'' part of the first statement.
The ``only if'' we saw in Proposition \ref{inducedshadowmap}.

Now,  if $\rho_{H_{\Gamma}}$ is any $u$-simply transitive 
algebraic group embedding which extends $\tau$ to $H_{\Gamma}$, then
$\jmath(t) = \rho_{H_{\Gamma}}(t)$ defines a homomorphism 
$\jmath: T \rightarrow  \N_{A}(\, U)$, $U = \overline{\tau(\Theta)}$,  
which satisfies $c(\,  \jmath(t)) =t $, for all $t \in T$. Therefore the conjugacy
statement of Proposition \ref{splitT} implies that any two $u$-simply transitive 
extensions of  $\tau$ to $H_{\Gamma}$ are conjugate by an element 
$u \in \C_A(U)$. 
\end{proof} 

\begin{prf}{Proof of Theorem \ref{shadowrealization}}
By the previous corollary, $\tau \in \Hom_c(\Theta,A)$ extends
to a $u$-simply transitive algebraic group homomorphism 
$\rho_{H_{\Gamma}}:  H_{\Gamma} \rightarrow A$. Let
$\rho$ be the restriction of $\rho_{H_{\Gamma}}$ to $\Gamma$.
Let $\gamma \in \Gamma$ be so that $\gamma_u \in \Theta$.
Since $\rho_{H_{\Gamma}}$ preserves the Jordan-decomposition,
$\tau(\gamma_u) = \rho_{H_{\Gamma}}(\gamma_u)= \rho(\gamma)_u$.
Therefore  $\tau$ descends to $\rho$ in the sense of Definition
\ref{uextension}.  By Theorem \ref{acextension}, $\rho$ is
a crystallographic homomorphism for $\Gamma$.
\end{prf} 

\begin{prf}{Proof of Theorem \ref{shadowembedding}}
It follows from Theorem  \ref{shadowrealization} that the 
map $\sigma_u$ is onto $\Def(\Theta, A)^T$. Let 
$\tau \in \Hom_c(\Theta,A)$.  By 
Theorem \ref{acextension},  every extension $\rho \in \Hom_c(\Gamma,A)$ 
of $\tau$ extends {\em uniquely\/} to a $u$-simply transitive homomorphism
$\rho_{H_{\Gamma}}:  H_{\Gamma} \rightarrow A$ which then is also 
an extension of $\tau$. 
Therefore the conjugacy statement of Corollary \ref{hullrealization}
implies that $\sigma_u$ is injective. By Corollary \ref{shadowmapisc},
the shadow map $s_u: \Hom_c(\Gamma, A) \rightarrow \Hom_c(\Theta,A)$
is algebraic, in particular $s_u$ is continuous, hence also $\sigma_u$.
\end{prf} 

\begin{prf}{Proof of Theorem \ref{shadowstratum}}
To prove that $\sigma_u$ is an open mapping on $\Def(\Gamma,A)_s$
we show that $s_u: \Hom_c(\Gamma,A)_s \rightarrow  \Hom_c(\Theta,A)_s$
admits a continuous section. This is done in a manner completely 
analoguous to the proof of Lemma \ref{csection}. We leave the
details to the reader.  
\end{prf}


\section{Applications to the existence of affine crystallographic actions} 
\label{existence}
One of the main questions of the subject is to decide whether
a given virtually polycyclic group $\Delta$ may act as an
affine crystallographic group. Not many general results on
this question seem to be known. (But see 
\cite{GM} for the classification of 
torsionfree polycyclic groups which act crystallographically
by affine Lorentz-transformations.)
Our methods here apply when
some information on the unipotent shadow of $\Delta$ 
is available. In general, this suggests that the main difficulty
of the problem lies in the existence and structure
of crystallographic actions of torsionfree nilpotent 
groups. In section \ref{convexamples} we exhibited 
some classes of well understood $\ftn$-groups.
Since these examples 
have strongly convex deformation spaces our realization
results provide a positive answer for the existence of
crystallographic actions for their finite extensions and
also certain associated polycyclic groups. (See Corollary
\ref{convexftnbfgroups}, Corollary \ref{convexpcgroups})

\paragraph{Crystallographic actions of finite extensions}
Let $\Gamma$ be a \ftn-group. 
We say that $\Gamma$ admits an invariant grading 
if the Lie algebra of the real Malcev hull $U_{\Gamma}$
has a positive grading which is preserved by a Levi
subgroup of $\Aut(U_{\Gamma})$. 
(Compare Definition \ref{filtrations}.)

\begin{corollary} \label{nilexamples}
Let $\Gamma$ be a \ftn-group which satisfies one of the
following conditions 
\begin{itemize}
\item[i)]   $\Gamma$ is of nilpotency class $\leq 3$,       
\item[ii)]  $\Gamma$ has rank $\leq 5$,  
\item[iii)] $\Gamma$ admits an invariant grading.  
\end{itemize}
Then $\Gamma$ is isomorphic to an affine crystallographic group 
which satisfies the realization property. Moreover, 
if $\Delta$ is a finite effective extension
of\/ $\Gamma$ then $\Delta$ is isomorphic to an 
affine crystallographic group, and $\Delta$ has
the realization property.
\end{corollary} 
\begin{proof} By Theorem \ref{convexftngroups}, the deformation
space of $\Gamma$ is convex, in particular, it is not empty. 
If $\Delta$ is a finite effective extension of $\Gamma$, then the associated
kernel $\beta: \Delta/\Gamma \rightarrow \Out(\Gamma)$ has
a fixed point in $\Def(\Gamma,\Aff(V))$. 
By Theorem \ref{realizationA},  $\Delta$ admits an affine
crystallographic action. Since the convexity of the
deformation space is inherited to $\Delta$, also $\Delta$
has the realization property. 
\end{proof} 

\begin{remark} Our corollary includes the result of Lee (\cite{Lee}) who
proved that a finitely generated torsionfree virtually nilpotent group $\Gamma$
with  $\rank \Gamma \leq 3$ acts as an affine crystallographic group. 
\end{remark}

\paragraph{Examples of affine crystallographic polycyclic groups}
Let $\Gamma$ be a torsionfree polycyclic group. We associated 
to $\Gamma$ its unipotent shadow  $\Theta \leq U_{\Theta}$, 
and the semisimple kernel 
$T_{\Gamma} \leq \Out(U_{\Theta})$. 
Let
$\lu_{\Theta}$ be the Lie algebra of $U_{\Theta}$.
We say that $T_{\Gamma}$ 
centralizes a nonsingular derivation $D$ of $\lu_{\Theta}$
if there exists a subgroup $T \leq \Aut(U_{\Theta})$ which
projects onto $T_{\Gamma}$ and centralizes $D$.
\begin{corollary} 
Let $\Gamma$ be a torsionfree polycyclic group
with unipotent shadow $\Theta$. We assume that $\Theta$ 
satisfies one of the following conditions 
\begin{itemize}
\item[i)]   $\Theta$ is of nilpotency class $\leq 3$,   
\item[ii)]  $\Theta$ admits an invariant grading,
\item[iii)] The semisimple kernel $T_{\Gamma}$ of\/ $\Gamma$ centralizes 
a nonsingular derivation $D$ of $\lu_{\Theta}$. 
\end{itemize}
Then $\Gamma$ is isomorphic to an affine crystallographic group. 
\end{corollary} 
\begin{proof} By Theorem \ref{convexftngroups} conditions i) and ii)
imply that $\Def(\Theta,\Aff(V))$ is strongly convex. In
particular, the semisimple kernel $T_{\Gamma} \leq \Out(U_{\Theta})$
has a fixed point in $\Def(\Theta,\Aff(V))$. 
If condition iii) is satisfied it follows from Proposition \ref{derivation}
that there exists an affine crystallographic action $\rho$ of $\Theta$
so that
$T_{\Gamma} \leq \Out(U_{\Theta})_{[\rho\,]}$, that is, $T_{\Gamma}$ 
fixes the point $[\rho\,] \in \Def(\Theta, \Aff(V))$.  
Therefore, in all three cases, it follows from  Theorem \ref{shadowrealization} that 
$\Gamma$ admits an affine crystallographic action. 
\end{proof}

To illustrate the corollary, we consider the particular case 
that $\Gamma$ is a Zariski-dense
lattice in a semi-direct product. Then the 
conditions on the unipotent shadow $\Theta$ of $\Gamma$
become in particular transparent. 
\begin{example1}[semi-direct products] \rm
Let $G$ be
a connected simply connected solvable Lie group,
and $\Gamma \leq G$ a Zariski-dense lattice.
We assume that $G$ is a semisimple semi-direct product, that is, 
the Lie algebra $\lg$ of $G$ splits as a direct sum 
$$  \lg = \la \oplus \ln \; , $$
where $\ln$ is the nilpotent radical of $\lg$, 
and $\la$ is an abelian subalgebra which acts
on $\ln$ by semisimple transformations. The unipotent
shadow $U$ of $G$ is obtained by killing the torus 
action. That is, the Lie algebra $\lu$ of $U$ 
is just the direct product 
$$  \lu = \lb  \oplus \ln \; , $$  
of an abelian ideal $\lb$ and the ideal $\ln$. 
The Lie algebra $\la$ acts on $\lu$ by centralizing
$\lb$, and by the adjoint action on $\ln$. Moreover,
the Zariski closure of the semisimple abelian subgroup $T_A \leq \Aut(U)$,
which belongs to $\la$, is a torus $T \leq \Aut(U)$ 
which projects onto the semisimple kernel $T_{\Gamma} \leq \Out(U)$. 
\end{example1} 

Therefore we get:
\begin{corollary} Let $G$ be a simply connected solvable 
Lie group which is a semisimple semi-direct product 
$G= T_A N$, where $N$ is the nilradical of $G$. Assume
that $G$ has a Zariski-dense lattice $\Gamma$, and
assume further that one of the following conditions 
is satisfied:
\begin{itemize}
\item[i)]   $N$ has a positive grading which is invariant by $T_A$.
\item[ii)]  $T_A$ centralizes a nonsingular derivation of $N$. 
\item[iii)] $N$ is of nilpotency class $\leq 3$. 
\end{itemize}
Then $\Gamma$ is isomorphic to an affine
crystallographic group.    
\end{corollary} 

%
%
%
%


\section{Group actions on deformation spaces} \label{Groupactions}
Let $\Gamma$ be an \ftn-group, and $\Delta$ a finite effective 
extension group of $\Gamma$.  We let $\bU$ denote the Malcev
completion of $\Gamma$,  and $\bU_{\Delta}$ the algebraic hull
for $\Delta$.  
In this section  we view $\Gamma$ as a lattice in the connected 
simply connected (real) Lie group $U_{\Gamma} = \bU_{\bbR}$, and
$\Delta$ as a subgroup of $U_{\Delta}= \bU_{\Delta, \bbR}$. 
The purpose of this final section is to exploit the functorial 
properties of the real algebraic hull  $U_{\Delta}$ in the study of the  
deformation space $\Def(\Delta, A)$, and to provide some 
auxiliary results which we need at  various places in this work.

\paragraph{The hull functor on deformation spaces} 
Recall that the group $\Aut(\Delta)$ acts in 
a natural way on $\Hom_c(\Delta,A)$. In fact, the action
of $\phi \in \Aut(\Delta)$ is described for all
$\rho \in  \Hom_c(\Delta,A)$ by 
$$    \rho \longmapsto \rho^{\phi} = \rho  \circ \phi^{-1} \;. 
$$
The $\Aut(\Delta)$-action on $\Hom_c(\Delta,A)$ induces then an 
$\Out(\Delta)$-action on the deformation space $\Def(\Delta,A)$.
We fix an embedding $j: \Delta \hookrightarrow U_{\Delta}$ of
$\Delta$ as a lattice in its real algebraic hull $U_{\Delta}$. 
By Corollary \ref{hullrigidity}
the embedding $j$ induces an embedding 
$$ \epsilon_{j}: \Aut(\Delta) \longrightarrow  \Aut( U_{\Delta} ) \; . $$
We show now that the action of  $\Aut(\Delta)$ 
extends to an action of  $\Aut( U_{\Delta})$. 
By Proposition \ref{hullextension},  there exists,  for every
 $\rho \in  \Hom_c(\Delta,A)$, a unique homomorphism 
$ \bar{\rho}=\rho_{U_{\Delta}}:  U_{\Delta}  \rightarrow A$ which
satisifies $ \bar{\rho} \circ  j =  \rho$.
For $\rho \in \Hom_c(\Delta, A)$, $ \Phi \in \Aut(U_{\Delta})$,  
we put then  
$$   \rho^{\Phi} = \bar{\rho} \circ \Phi^{-1} \! \circ j \; .$$
Since  $\rho^{\Phi}(\Gamma)$ is a lattice in the unipotent simply
transitive group $\bar{\rho}(U_{\Gamma})$ it follows that $\rho^{\Phi}(\Delta)$ is 
an ACG, hence  $\rho^{\Phi} \in \Hom_c(\Delta,A)$.

\begin{proposition} \label{extaction} 
The correspondence $\rho \longmapsto \rho^{\Phi^{-1}}$ extends the
$\Aut(\Delta)$-action on $\Hom_c(\Delta,A)$ to an 
action of $\Aut(U_{\Delta})$ on $\Hom_c(\Delta,A)$. 
The $\Aut(U_{\Delta})$-action on $\Hom_c(\Delta,A)$ is free, and induces an
action of $\Out(U_{\Delta})$ on $\Def(\Delta,A)$. In particular the
action of $\Out(\Delta)$ factorizes over $\Out(U_{\Delta})$.
\end{proposition}
\begin{proof} Let $\Phi, \Psi \in \Aut(U_{\Delta})$.  
We calculate  
$ (\rho^{\Phi})^{\Psi} = \overline{\rho^{\Phi}} \circ \Psi^{-1} \circ j$, where
$\overline{\rho^{\Phi}}: U_{\Delta} \rightarrow A$ is the extension of 
$\rho^{\Phi}: \Delta \rightarrow A$.  
Therefore $\overline{\rho^{\Phi}} = \bar{\rho} \circ \Phi^{-1}$, and we get
$ (\rho^{\Phi})^{\Psi} =  \bar{\rho} \circ ( \Phi^{-1} \circ \Psi^{-1} ) \circ j 
= \rho^{\Psi \circ \Phi}$.
So, in fact, $\Aut(U_{\Delta})$ acts on $\Hom_c(\Delta,A)$. 
Obviously, the $\Aut(U_{\Delta})$-action extends the $ \Aut(\Delta)$-action, 
in the sense that, for $\phi \in \Aut(\Delta)$,  
$\rho^{\phi} = \rho^{\epsilon_j( {\phi}) }$.
Since $\Delta$ 
is Zariski-dense in $U_{\Delta}$, $\Aut(U_{\Delta})$ acts freely. 

For $h \in U_{\Delta}$, $\gamma \in \Delta$, let $c_h: U_{\Delta} \rightarrow U_{\Delta}$
denote conjugation with $h$. We get 
$ \rho^{c_h} (\delta)  = \bar{\rho}(h) \rho(\delta)  \bar{\rho}(h)^{-1}$.
Therefore  $\Out(U_{\Delta})$ acts on $\Def(\Delta,A)$. 
\end{proof} 

\begin{remark} 
The $\Aut(U_{\Delta})$-action depends on the embedding 
$j: \Delta \rightarrow U_{\Delta}$. 
Nevertheless, the image of  $\Aut(U_{\Delta})$ as a group of morphisms of  
$\Hom_c(\Delta, A)$
is independent from the choice of embedding. In particular, $\Out(U_{\Delta})$ maps to 
a well defined subgroup of homeomorphisms on the deformation space $\Def(\Delta,A)$.  
\end{remark}

The presence of the $\Aut(U_{\Delta})$-action on the deformation
space $\Def(\Delta,A)$ comes from  a correspondence of the
affine crystallographic representations of $\Delta$ with certain representations
of its real hull $U_{\Delta}$.  We want to briefly explain this now: \newline 

\noindent 
Let $H$ be a real algebraic group which contains a 
normal unipotent subgroup $U$ of finite index, and 
assume that the centralizer of $U$ is contained in $U$. 
Note that the assumptions on $H$ are satisfied by the 
algebraic hull $U_{\Delta}$.

\begin{definition} 
A homomorphism $\bar{\rho}: H \rightarrow A$
is called {\em crystallographic\/} if  $\bar{\rho}$ is injective 
and if $U$ acts simply transitively on $V$.  
\end{definition}
We then define the space $\Hom_c(H, A)$ of crystallographic
homomorphisms of $H$ and the  corresponding deformation 
space $\Def(H, A) =  \Hom_c(H, A) /A$. 
The space 
$\Hom(H, A)$ is a topological space with
the compact open topology. Also by the results
of chapter \ref{choms} there is a natural
structure of real algebraic variety on $\Hom_c(H, A)$.

We fix an embedding $j: \Delta \hookrightarrow U_{\Delta}$ of
$\Delta$ in its real algebraic hull $U_{\Delta}$. Recall that 
every $\rho \in \Hom_c(\Delta,A)$ is unipotent on the Fitting subgroup $\Gamma$
of $\Delta$. 
By Proposition \ref{hullextension},  $\rho \in \Hom_c(\Delta,A)$ extends to a 
representation $\rho_{U_{\Delta}}$ of $U_{\Delta}$. Without a proof
we state:

\begin{theorem} \label{hullequivalence}
The correspondence $\rho \longmapsto \rho_{U_{\Delta}}$ defines a
homeomorphism  
$$    h_{j}:  \Hom_c(\Delta,A) \xrightarrow{\; \; {\scriptstyle \approx} \; \;}  \Hom_c(U_{\Delta}, A)$$
which commutes with the natural actions of automorphism groups
on both spaces. In particular, $h_{j}$ induces a homeomorphism 
$\bar{h}_{j}:  \Def(\Delta,A) \xrightarrow{\; \; {\scriptstyle \approx} \; \;}  \Def(U_{\Delta},A)$ which
commutes with the natural outer automorphism actions.  
\end{theorem}

The theorem implies that the topology and structure of the deformation
space $\Def(\Delta,A)$ for $\Delta$, and also of the moduli space $\Mod(\Delta,A)$,
depend only on the properties of the  algebraic hull for $\Delta$. \\

Here is another application. Look at the variety 
${\cal C} = \Hom_c(U_{\Delta}, A) \times V$. Then $A$ 
acts on $\cal C$, where $g (\rho, v) = (\rho^g, g v)$.
Moreover, $\Aut(U)$ acts on the first factor of ${\cal C}$,
commuting with the action of $A$. Define ${\cal LC}(U) = {\cal C} / A$.
 
\begin{proposition} \label{Matsushima}
The deformation space $\Def(U,A)$ is a quotient of the variety 
${\cal LC}(U)$ by the induced unipotent action of ${\rm Inn}(U)$ on ${\cal LC}(U)$
\end{proposition}
\begin{proof}  Since $A$ acts transitively on $V$, we can deduce
that $$\Def(U,A) =  {\cal LC}(U) / {\rm Inn}(U) \; .$$
Furthermore,  ${\cal LC}(U)$ is the set of  $A$-conjugacy classes
of \'etale representation with basepoint, and this corresponds  to
the algebraic variety of left-symmetric algebra products on the
Lie algebra $\lu$. (See, for
example,  \cite{BauesCortes}, for discussion.) 
\end{proof}

\paragraph{Stabilizers in deformation space}
For a moment we restrict our considerations
to \ftn-groups. Let $\Gamma \leq U_{\Gamma}$ be a \ftn-group 
and assume that $\Gamma$ is an ACG. Let $\rho \in \Hom_c(\Gamma, A)$
be an affine crystallographic representation, and 
$[\rho \, ] \in \Def(\Gamma,A)$ the corresponding point in
the deformation space. We let $U= \bar{\rho}(U_{\Gamma}) \leq A$ 
be the unipotent simply transitive hull for $\rho(\Gamma)$. The group
$ \Aut(\U_{\Gamma})$ acts on  $\Def(\Gamma,A)$.
The stabilizer $\Aut(U_{\Gamma})_{[\rho \,]} \leq \Aut(U_{\Gamma})$ 
of $[\rho \, ]$ may be desribed in terms of the normalizer $\N_A(U)$
of the unipotent simply transitive hull $U$. 

\begin{proposition} \label{Ustabilizer} 
There are  canonical isomorphisms 
\begin{eqnarray*} 
        \N_A(U) \rmod \C_{A} (U) & \stackrel{\cong}{\longrightarrow} & 
        \Aut(U_{\Gamma})_{[\rho \,]} \;\; , \\
        \N_A(U) \rmod (\C_A(U) U)&  \stackrel{\cong}{\longrightarrow} & 
         \Out(U_{\Gamma})_{[\rho \, ]}  \; . 
\end{eqnarray*} 
\end{proposition}
\begin{proof} 
$\N_A(U) \rmod \C_{A}(U)$ acts by conjugation as a group of automorphisms 
of $U$.  Since $\bar{\rho}$ is an isomorphism, there exists 
for each $g \in \N_A(U)$ an element $\Phi_g \in \Aut(U_{\Gamma})$ such that
$  \rho^g =  \rho^{\Phi_g}$.  This is,  precisely, the condition that 
$[ \rho^{\Phi_g} ] =   [\, \rho\, ]$. Therefore the map $g \mapsto \Phi_g$ 
defines a homomorphism of $\N_A(U)$ onto $\Aut(U_{\Gamma})_{[\rho\,]}$ with
kernel $\C_A(U)$. 
This  factorizes to an isomorphism
$ \N_A(U) \rmod (\C_A(U) \, U) \longrightarrow 
\Out(U_{\Gamma})_{[\rho \, ]}$.
\end{proof}

We consider next the $\Aut(\Gamma)$-action on $\Def(\Gamma,A)$.
Let us put 
$$ \Out_{A,\rho}(\Gamma)  = \N_A(\rho(\Gamma)) \rmod (\C_A(U) \rho(\Gamma)) \; . $$
We describe the stabilizer of $[\rho \,]$
in terms of the normalizer $\N_{A}(\rho(\Gamma)) \leq A$.

\begin{proposition} \label{Gstabilizer} 
There are canonical isomorphisms 
\begin{eqnarray*} 
          \N_{A}(\rho(\Gamma)) \rmod \C_{A} (U) & \stackrel{\cong}{\longrightarrow} &
           \Aut(\Gamma)_{[\rho \, ]} \;\; , \\
           \Out_{A,\rho}(\Gamma)  &  \stackrel{\cong}{\longrightarrow} 
                 & \Out(\Gamma)_{[\rho \, ]}  \; .  \\
\end{eqnarray*} 
\end{proposition}
\begin{proof} As in  Proposition \ref{extaction} we view $\Aut(\Gamma)$ as
a subgroup of $\Aut(U_{\Gamma})$. With this identification
$\Aut(\Gamma)_{[\rho \, ]} =  \Aut(\Gamma) \cap  \Aut(U_{\Gamma})_{[\rho \, ]}$. 
By Proposition \ref{Ustabilizer} conjugation on $U$ induces
a surjective homomorphism $c: \N_A(U) \longrightarrow \Aut(U_{\Gamma})_{[\rho \,]}$.
Therefore
$$ \Aut(\Gamma)_{[\rho \, ]}  =  \Aut(\Gamma) \cap c(\N_A(U)) 
= c( \N_A(\rho(\Gamma))) \; . $$
This shows that $\N_A(\rho(\Gamma))$ is mapped onto 
$\Aut(\Gamma)_{[\rho \, ]}$. The first isomorphism follows. 

Now $\Out(\Gamma)_{[\rho]}$ is just the image 
of $ \Aut(\Gamma)_{[\rho]\, }$ in $\Out(\Gamma)$. Hence, 
$\Out(\Gamma)_{[\rho]}$ is a quotient of  $\N_A(\rho(\Gamma))$,
and clearly the kernel is $\C_A(U) \rho(\Gamma)$.
\end{proof}

\paragraph{The normalizer of a unipotent simply transitive group}
Let $A \leq \Aff(V)$ be a Zariski-closed subgroup, and
$U \leq A$ a simply transitive unipotent group. We
discuss some of the structure of the normalizer $\N_{A}(U)$ of $U$ in $A$. 

\begin{lemma} \label{centralizer}  The centralizer
$\C_A(U) \leq A$ is a unipotent normal subgroup of $\N_{A}(U)$  and
acts freely on $V$. The group  $\C_A(U) \, U$ is a unipotent 
normal subgroup in $\N_A(U)$. 
\end{lemma}   
\begin{proof} $\C_A(U)$ acts freely  since $U$ is simply
transitive, and clearly  $\C_A(U)$ is a normal Zariski-closed
subgroup of $A$. Since any reductive subgroup of $\C_A(U)$ 
has fixed points on $V$,  $\C_A(U)$ is unipotent. 
\end{proof}

There are some subgroups of $\Aut(U)$ which
are associated to the normalizer $N_A(U)$, and 
which are important in our context. (These groups
also carry some (differential-) geometric interpretations 
in terms of the simply transitive
group action of $U$ and the associated flat left invariant 
connection on the real Lie group $U$. We do not go into these
details here, but see a related discussion in \cite{FriedGoldman}[\S 3.11].)   
We consider the homomorphism 
$$    c_U :   \N_A(U) \longrightarrow \Aut(U) \;\, ,    $$
where $c_U(g): U \rightarrow U$ is given by conjugation with $g \in  \N_A(U)$.   
Let us define subgroups $\Aff_A(U)$, 
$\Aut_A(U)$, $\Inn_A(U) \leq \Aut(U)$ as follows  
$$
\begin{array}{lclrr}
 \Aff_A(U)  & = &  c(\N_A(U))  &  (\cong \,  \N_{A}(U) \rmod \C_A(U) ) \\
 \Aut_A(U)  & = &  c(\N_{\GL_A}(U))   &  (\cong \, \N_{A}(U) \rmod  U)   \\
 \Inn_A(U)  & = &  c(\N_{\GL_A}(U) \cap \C_A(U) U)  &  (\cong \,  (\C_A(U)U)  \rmod U) \\
\end{array} 
$$ 
Finally, we put $\Out_A(U)$ for the
image of $\N_A(U)$ in $\Out(U)$. Note that the 
group $\Aut_A(U)$ projects onto $\Out_A(U)$. 
Proposition \ref{Ustabilizer} implies now  
\begin{proposition} \label{Autstabilizer} 
Let $\Gamma$ be a \ftn-group and $\rho \in \Hom_c(\Gamma,A)$. 
If $U \leq A$ is the simply transitive hull for $\rho(\Gamma)$ then there are
natural isomorphisms 
\begin{eqnarray*} 
\Aff_A(U)  & \stackrel{\cong}{\longrightarrow} &  \Aut(U_{\Gamma})_{[\rho \,]} \\
\Out_A(U)  & \stackrel{\cong}{\longrightarrow} &  \Out(U_{\Gamma})_{[\rho \, ]}
\; \; \; \; .
\end{eqnarray*} 
\end{proposition} 
\begin{proof} In fact, the representation $\bar{\rho}$ defines
an isomorphism $\varrho: U_{\Gamma} \rightarrow U$, and
the proof of Proposition \ref{Ustabilizer} shows that
$\Aut(U_{\Gamma})_{[\rho \,]} = \varrho^{-1} \Aff_A(U) \varrho$. 
The statement for $\Out_A(U)$ follows as well. 
\end{proof}

The following are geometrically interesting special cases:  
\begin{definition} Let $U \leq A$ be a simply transitive unipotent subgroup. 
$U$ is called {\em $A$-symmetric\/} if $\C_A(U)$ is transitive on $V$. 
We call $U$ {\em fully $A$-symmetric\/} if $U$ is $A$-symmetric and
$c_U(\Out_A(U))= \Out(U)$.
Finally,  we call $U$  {\em $A$-invariant} if $c_U(\Out_A(U))= \Out(U)$,
and $U$ is called  {\em $A$-convex} if $c_U(\Out_A(U))$ contains 
a Levi subgroup of  $\Out_A(U)$
\end{definition}
In our context these conditions are interpreted 
in terms of group actions on the deformation spaces.

\paragraph{The normalizer of a unipotent ACG}
Let us consider now the normalizer $\N_A(\, \rho(\Gamma))$
of $\rho(\Gamma)$ in $A$. The Zariski-closure of 
$U = \overline{\rho(\Gamma)}$ is a simply transitive unipotent
subgroup of type $A$. By Zariski-denseness of $\rho(\Gamma)$ in $U$, 
we have  $\N_A(\, \rho(\Gamma)) \leq \N_A(U)$, as well as 
$\C_A(\rho(\Gamma)) = \C_A(U)$.  
%
For $x \in V$,  we put  $\N_{A,x}(\, \rho(\Gamma))=   \N_A(\, \rho(\Gamma))\,  \cap A_{x}$. 
Our interest will be in the image of $\N_{A,x}$ in $\Out_{A,\rho}(\Gamma)$. 
In a special case, the projection from $\N_{A,x}$ to  $\Out_{A,\rho}(\Gamma)$
is surjective, for all $x \in V$. 
\begin{lemma} \label{symmetric}
If $U$ is an $A$-symmetric simply transitive group then, for all $x \in V$,
the natural map  
$$ 
     \N_{A,x}(\, \rho(\Gamma)) \longrightarrow  \N_A(\, \rho(\Gamma)) \rmod \C_A(U)
$$
is an isomorphism. 
\end{lemma}
\begin{proof} Put $H=\N_A(\, \rho(\Gamma))/  \C_A(U)$. $H$ acts on the orbit space 
$\C_A(U) \lmod V$. It is easy to see that $\N_{A,x}(\, \rho(\Gamma))$ projects 
onto $H_{[x]}$,  where $H_{[x]}$ is the stabilizer of $[x] = \C_A(U)x$. 
If $U$ is fully symmetric then $ \Z_A(U)$ acts transitively on $V$.  Hence $H= H_{[x]}$,
$\N_{A,x}(\, \rho) \cong \N_A(\, \rho(\Gamma))/\C_A(U)$. 
\end{proof}

\begin{definition} $\rho \in \Hom_c(\Gamma,A)$ is called symmetric, 
if $U = \overline{\rho(\Gamma)}$ is $A$-symmetric.  
\end{definition} 
Let $c_{\rho}: N_A(\rho(\Gamma)) \rightarrow \Aut(\Gamma)$ denote
the conjugation homomorphism, and $\Hom_c(\Gamma,A)_s$ the set
of symmetric crystallographic homomorphisms. The next fact was required
in section \ref{FExtdef}.
\begin{proposition} \label{s_x} 
Let $F \leq \Aut(\Gamma)$ be a finite subgroup. 
There exists a continuous map  
$s_x:  \Hom_c(\Gamma,A)_s^F \longrightarrow \Hom(F, A_x)$
such that $c_{\rho}( s_x (\rho, g) ) = g $, for all $g \in F$. 
\end{proposition}
\begin{proof} 
For every $\rho \in \Hom_c(\Gamma,A)$, $\bar{\rho}(U_{\Gamma})$ 
is a simply transitive subgroup of $A$. Therefore there exists, for every $g \in \Aut(\Gamma)$,
a unique polynomial diffeomorphism $\Phi_g$ of $V$ which satisfies
$\Phi_g \bar{\rho}(u) x = \bar{\rho} (g u) x$. In particular, 
$\Phi_g  \bar{\rho}(u) \Phi_g^{-1}  =  \bar{\rho} (g u)$ holds, for all
$u \in U_{\Gamma}$.     
Now, if $\rho$ is fixed by $F$,  and $\rho \in \Hom_c(\Gamma, A)_s$ 
then by Proposition \ref{Gstabilizer}, 
and Lemma \ref{symmetric} there exists a $\phi \in N_{A,x}(\rho(\Gamma))$
such that $\phi  \bar{\rho}(u) \phi^{-1} = \bar{\rho} (g u)$. We remark that,
since $\C_A(U)$ is transitive on $V$, it coincides with the centralizer of $U$ 
in the group of diffeomorphisms of $V$.
We conclude that $\Phi_g \phi^{-1} \in \C_A(U)$, and in particular
$\Phi_g \in N_{A,x}(U)$. So we set $s_x (\rho, g)=  \Phi_g$. The map
$s_x$ is easily seen to be continuous, since it is the restriction
of the continuous map on $\Hom_c(\Gamma,A)$ which assigns to
$\Phi_g(\rho)$ its first derivative in $x$. 
\end{proof}    

We remark that symmetric simply transitive
actions may be constructed for many unipotent
Lie groups $U$.   
\begin{example} Let $\cal A$ be a nilpotent associative algebra, finite dimensional 
over $\bbR$ and
$\lu$ the commutation Lie algebra of $\cal A$. Then it is known,
compare \cite{Auslander,Milnor}, that the corresponding simply
connected Lie group $U$ admits a simply transitive 
unipotent representation. It is possible to show that this 
representation is also symmetric, and that every symmetric
simply transitive representation arises from this construction.  
\end{example}




\end{document}